\documentclass[a4paper]{article}

\usepackage[english]{babel}
\usepackage[utf8x]{inputenc}
\usepackage[T1]{fontenc}
\usepackage{booktabs}
\usepackage{subfig}
\usepackage[titletoc,title]{appendix}
\usepackage{amsfonts}
\usepackage{float}
\usepackage{array}
\usepackage{multirow}
\usepackage[title]{appendix}
\usepackage[a4paper,top=3cm,bottom=2cm,left=3cm,right=3cm,marginparwidth=1.75cm]{geometry}
\usepackage{ulem}
\usepackage{amsmath}
\usepackage{graphicx}
\usepackage[colorinlistoftodos]{todonotes}
\usepackage{authblk}
\usepackage{soul}
\usepackage{multirow}

\newtheorem{lemma}{Lemma}[section]

\newcommand{\mcO}{\mathcal{O}}

\newcommand{\mcL}{\mathcal{L}}

\newcommand{\mbR}{\mathbb{R}}
\newcommand{\mbRd}{{\mathbb{R}^d}}

\newcommand{\omg}{{\Omega}}
\newcommand{\omgd}{{\Omega_\delta}}
\newcommand{\unn}{{u_{N\!N}}}

\def \mub{{\boldsymbol \mu}}

\def \xb{{\bf x}}
\def \yb{{\bf y}}

\def \thetab{{\boldsymbol\theta}}

\title{nPINNs: nonlocal Physics-Informed Neural Networks for a parametrized nonlocal universal Laplacian operator. Algorithms and Applications}

\author[1]{G. Pang}
\author[2]{M. D'Elia}
\author[2]{M. Parks}
\author[1]{G. E. Karniadakis}
\affil[1]{Division of Applied Mathematics, Brown University, RI}
\affil[2]{Center for Computing Research, Sandia National Laboratories, NM}
\date{}
\begin{document}
\maketitle
\begin{abstract}
    Physics-informed neural networks (PINNs) are effective in solving inverse problems based on  differential and integro-differential equations with sparse, noisy, unstructured, and multi-fidelity data. PINNs incorporate all available information, including governing equations (reflecting physical laws), initial-boundary conditions, and observations of quantities of interest, into a loss function to be minimized, thus recasting the original problem into an optimization problem. In this paper, we extend PINNs to parameter and function inference for integral equations such as nonlocal Poisson and nonlocal turbulence models, and we refer to them as nonlocal PINNs (nPINNs). The contribution of the paper is three-fold. First, we propose a unified nonlocal Laplace operator, which converges to the classical Laplacian as one of the operator parameters, the nonlocal interaction radius $\delta$ goes to zero, and to the fractional Laplacian as $\delta$ goes to infinity. This universal operator forms a super-set of classical Laplacian and fractional Laplacian operators and, thus, has the potential to fit a broad spectrum of data sets. We provide theoretical convergence rates with respect to $\delta$ and verify them via numerical experiments. Second, we use nPINNs to estimate the two parameters, $\delta$ and $\alpha$, characterizing the kernel of the unified operator. The strong non-convexity of the loss function yielding multiple (good) local minima reveals the occurrence of the \textit{operator mimicking} phenomenon, that is, different pairs of estimated parameters could produce multiple solutions of comparable accuracy. Third, we propose another nonlocal operator with spatially variable order $\alpha(y)$, which is more suitable for modeling turbulent Couette flow. Our results show that nPINNs can jointly infer this function as well as $\delta$. More importantly, these parameters exhibit a universal behavior with respect to the Reynolds number, a finding that contributes to our understanding of nonlocal interactions in wall-bounded turbulence.
\end{abstract}


\section{Introduction}

Nonlocal models can capture effects that standard partial differential equations (PDEs) fail to capture thanks to their ability to describe long-range interactions.
As an example, they can better model physical problems that exhibit multi-scale behavior (e.g., fracture and failure in solid materials) as well as diffusion phenomena whose mean square displacement is not a linear function of time (e.g., super- and sub-diffusion occurring in solute transport in ground-water). In the former case, we mention the peridynamics model for continuum mechanics introduced and extended in \cite{silling2000reformulation,silling2003deformation,silling2005peridynamic}; in the latter case, we refer to anomalous diffusion models described by fractional differential equations \cite{meerschaert1999multidimensional,benson2000application}, and more general nonlocal diffusion models such as those introduced in \cite{du2012analysis,d2013fractional,d2017nonlocal}. 

However, the increased accuracy of these models comes at the price of several modeling and numerical challenges that hinder their usability. Among others, we mention the uncertain nature of model parameters, which may be non-measurable, sparse, and subject to noise, and of the nonlocal kernels characterizing the functional form of the operator. In this work, we focus on the former concern and investigate how to discover proper nonlocal models based on a limited number of measurements of quantities of interest (QoIs). This task is, in general, much harder than parameter identification in a PDE setting, as in the nonlocal case the operator itself (e.g., its spectral properties) may change significantly as its parameters change.

Previous attempts to model identification in the context of nonlocal Laplace operators relied on optimal control strategies \cite{Antil2018SpectralControl,Antil2019OptControl,DElia2014DistControl,d2019priori,d2016identification,turner2015inverse,Zheng2019}. Here, the derivation of the adjoint equations can be problematic when the nonlocal kernel is not an affine function of the parameters to be identified, the so-called control parameters. In those cases taking the variational derivative of the objective functional with respect to the parameters is analytically intractable and an affine approximation of the kernel \cite{burkovska2019approximation} is needed to facilitate the solution of the optimization problem. Other works proposed surrogate-based inversion techniques \cite{garcia2006using,pang2017discovering,yan2012stochastic}, in which a large number of forward simulations act as the training set to construct a surrogate model that maps the parameters into QoIs or QoI-related quantities. These techniques usually work well for large training sets.  

We propose a new approach to model the learning that is in stark contrast with previously developed techniques. This paradigm shift is the combination of 1) {\it machine learning} and {\it physical principles}, and 2) {\it universal operators} and {\it versatile surrogates}, such as neural networks.
The outcome is a data-driven physics-informed framework for learning new complex nonlocal phenomena. While a physics-based approach to machine learning is now recognized as a very effective approach, augmenting neural networks with the structure provided by a unified operator is new and makes our technique more rigorous and physically consistent.

Moreover, a \textit{physics-informed learning} framework can considerably reduce the size of the training set by using the governing equations, which encode physical laws, as an implicit regularization term in an objective functional to be minimized. As a result, this implicit regularizer will quickly guide the optimization towards a good local minimum. The difference between this setting and the aforementioned optimal control techniques is that the latter uses the governing equations as constraints to the optimization problem, in a PDE-constrained optimization fashion. We are interested in physics-informed learning algorithms that adopt machine learning surrogates in the objective functional and do not require deriving adjoint equations. In this family of algorithms we mention physics-informed Gaussian processes (PIGPs) \cite{raissi2018hidden} and physics-informed neural networks (PINNs) \cite{raissi2019physics}. PIGPs can quantify the uncertainty of their predictions; however, when applying them to integral equations, analytically intractable and computationally expensive operations are involved. As an example, we have to evaluate the action of the composition of two nonlocal operators on a covariance function, i.e., $\mathcal{L}(\mathcal{L}\,k(\xb,\yb))$, where $\mathcal{L}$ is a nonlocal operator and $k(\cdot,\cdot)$ is the covariance function of the GP. On the other hand, PINNs for integral equations are much easier to implement, see paper \cite{pang2019fpinns}, where the authors propose a fractional version of PINNs, fPINNs. However, no uncertainty quantification is provided. We propose a generalized nonlocal version of PINNs, nPINNs, that reduces to PINNs when the nonlocality vanishes and to fPINNs for infinite nonlocal interactions. 

More specifically, the contribution of the paper is three fold. 

{\bf 1.} We introduce a universal nonlocal Laplace operator that is parameterized by two scalars: the nonlocal interaction radius $\delta$ and the decay rate, $\alpha$, of the nonlocal kernel. For finite $\delta$, this operator represents a superset of classical Laplacian and fractional Laplacian operators; in fact, as $\delta\rightarrow 0$ the operator reduces to the classical Laplacian $-\Delta$, whereas as $\delta\rightarrow \infty$ the operator converges to the fractional Laplacian $(-\Delta)^{\alpha/2}$. We provide theoretical convergence rates with respect to $\delta$ and verify them via numerical experiments. 

{\bf 2.} We use nPINNs to estimate  $\delta$ and $\alpha$; the strong non-convexity of the loss function may yield multiple good local minima -- this gives rise to the \textit{operator mimicking} phenomenon. By operator mimicking, we refer to the existence of multiple ($\delta$, $\alpha$) pairs that generate distinct operators, which are equally effective in reproducing the training data.

{\bf 3.} We propose a new nonlocal operator to model the total shear stress in the turbulent Couette flow. This operator is parameterized by the nonlocal interaction radius $\delta$ and a spatially variable decay exponent $\alpha(y^+)$, where $y^+$ is the wall coordinate. We successfully use nPINNs to jointly estimate this function as well as $\delta$ and discover a universal behavior with respect to the Reynolds number. This fact contributes to better understanding non-local interactions in wall-turbulence.

The rest of the paper is organized as follows. In Section \ref{sec:nonlocal-diffusion} we define the nonlocal Poisson problem. In Section \ref{sec:unified-operator} we introduce the unified nonlocal operator and prove the convergence rates of the operator to Laplacian and fractional Laplacian with respect to $\delta$. In Section \ref{sec:nPINNs-formulation} we describe the basic idea and implementation of nPINNs. In Section \ref{sec:numerical-tests}, we first numerically illustrate the convergence rates of the unified operator, then demonstrate the applicability and consistency of nPINNs by solving forward benchmark problems, and finally show results of model discovery with nPINNs. In Section \ref{sec:turbulence} we introduce a new nonlocal model for the turbulent Couette flow and estimate its parameters via nPINNs. We conclude the paper in Section \ref{sec:conclusion}. In Appendix \ref{sec:AppA} we report the proof of Lemma \ref{lem:limiting-behavior} on the limiting behavior of the unified operator, in Appendix \ref{sec:AppB} we describe how we compute singular integrals, and in Appendix \ref{sec:AppD} we report the results of numerical tests on the accuracy of quadrature rules for nonlocal operators.

\section{A nonlocal Poisson problem}\label{sec:nonlocal-diffusion}
In this section we introduce nonlocal Laplace operators and the corresponding equations. Given a nonnegative kernel function $\gamma(\xb,\yb;\mub)$ such that $\gamma(\xb,\yb;\boldsymbol{\mu})=\gamma(\yb,\xb;\boldsymbol{\mu})$ and a set of parameters $\mub$, we define the action of a nonlocal Laplace operator $\mcL^\mub$ on a scalar function $u:\mbRd\to\mbR$ as
\begin{equation}\label{original_operator}
\mathcal{L}^{\boldsymbol{\mu}}u(\xb):=\int_{\mathbb{R}^d}(u(\yb)-u(\xb))\gamma(\xb,\yb;\boldsymbol{\mu})d\yb.
\end{equation}
The set of kernel parameters $\mub$, usually unknown and non-measurable, determines regularity properties of the solution and is application dependent; its identification is the subject of our work. In \eqref{original_operator}, the integral form allows us to represent long-range interactions and reduces the regularity requirements on the solutions that, as opposed to the PDE setting, do not have to be differentiable.

The presence of long-range interactions has consequences on the way problems on bounded domains are solved, or, more specifically, on how nonlocal boundary conditions are prescribed. To clarify this concept, given an open and bounded domain $\Omega\subset\mbRd$, we introduce its interaction domain, or nonlocal boundary, as the set of those points outside of $\Omega$ interacting with points inside, i.e. 
\begin{equation}
\omg_I = \{\yb\in\mbRd\setminus\Omega: \gamma(\xb,\yb)\neq 0, \; 
{\rm for}\; \xb\in\Omega\}.
\end{equation}
As we point out later on, this is the set where conditions on the solution must be prescribed to guarantee the existence and uniqueness of the solution of a nonlocal Poisson problem in a bounded domain.

\begin{figure}[H] 
\centering
\includegraphics[width=.3\textwidth]{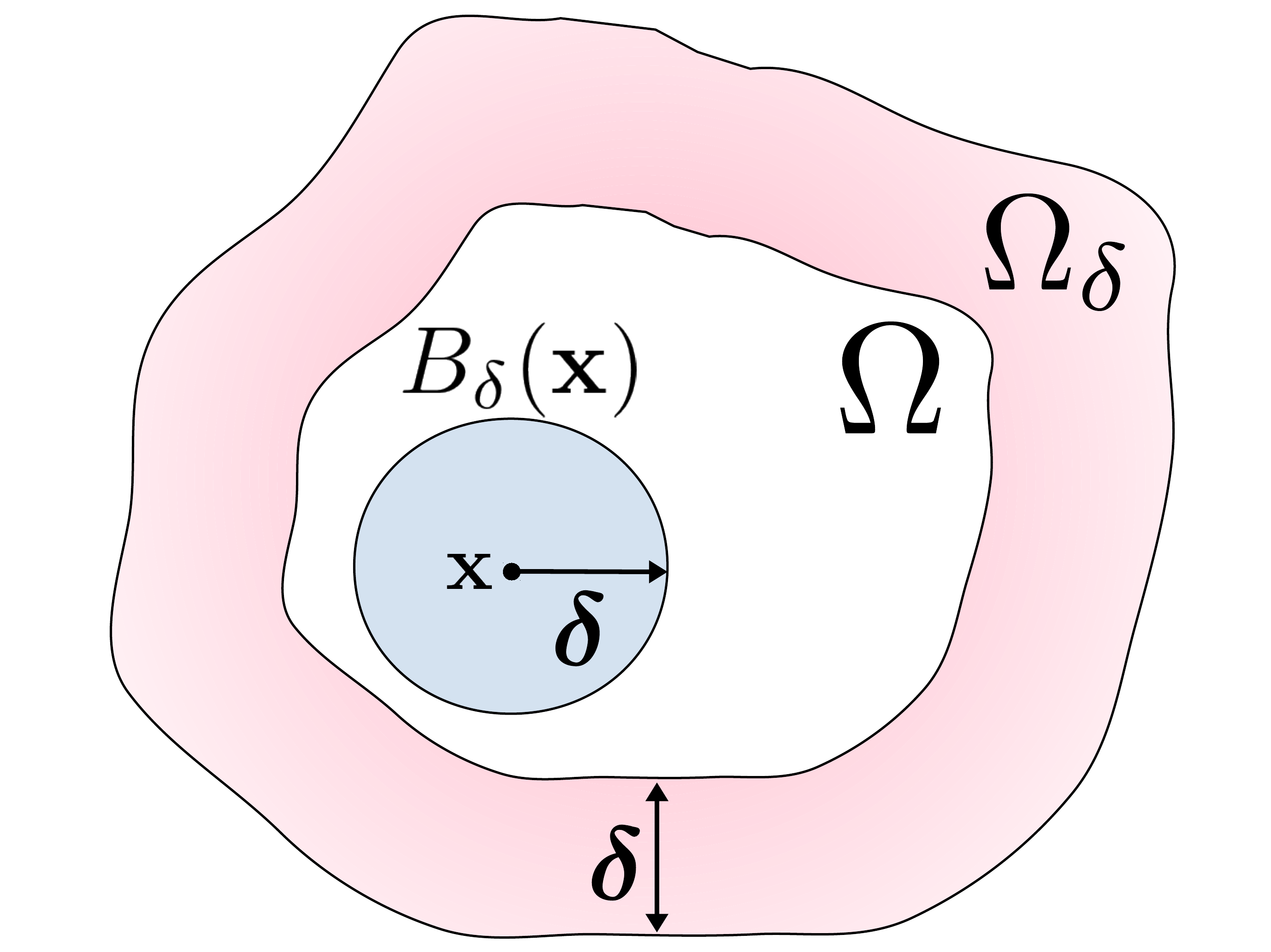}
\caption{\label{interaction_domain}Domain $\Omega$, interaction domain $\Omega_{\delta}$ and support of the kernel function around a point, $B_{\delta}(\mathbf{x})$.}
\label{fig:domain}
\end{figure}

In scientific and engineering applications, it is often the case that interactions have finite length; hence, we limit the action of the kernel to a {\it nonlocal neighborhood} $B_\delta(\xb)$, the Euclidean ball of radius $\delta$ centered at a point $\xb$. We refer to $\delta$ as the interaction radius. Thus, we require the kernel $\gamma$ to satisfy the following properties. Given positive constants $\delta$ and $\gamma_0$, 
\begin{equation}\label{eq:kernel-support}
\left\{\begin{array}{ll}
      \gamma(\xb,\yb;\delta)\ge 0, & \yb \in B_{\delta}(\xb), \\
      \gamma(\xb,\yb;\delta) = 0, & \forall \yb\in \mbRd \setminus B_{\delta}(\xb),\\
      \gamma(\xb,\yb;\delta)\ge \gamma_0 > 0, & \forall \yb\in B_{\delta/2}(\xb).\\
\end{array} \right.
\end{equation}
For kernels that satisfy \eqref{eq:kernel-support} the interaction domain takes the form
\begin{equation}
\omg_I=\omgd = \{\yb\in\mbRd\setminus\Omega: 
\|\xb-\yb\|\leq \delta, \; {\rm for}\; \xb\in\Omega\}.
\end{equation}
In Figure \ref{fig:domain} we present an illustration of $\omg$, $\omgd$ and the interaction neighborhood $B_\delta(\xb)$.

In this work we are interested in solving the nonlocal counterpart of a local Poisson problem in $\omg$, i.e., given $f:\omg\to\mbR$ and $g:\omgd\to\mbR$ we want to find $u:\omg\cup\omgd\to\mbR$ such that
\begin{equation}\label{eq:nonlocal-diffusion}
\begin{aligned}
-\mcL^\mub u(\xb)=-\int_{B_\delta(\xb)}(u(\yb)-u(\xb))
\gamma(\xb,\yb;\mub)d\yb = f(\xb), & \quad\xb\in\omg \\
u(\xb)=g(\xb) & \quad\xb\in\omgd,
\end{aligned}
\end{equation}
where we restricted the domain of integration to the nonlocal neighborhood. Here, the condition on $\omgd$ is a nonlocal analogue of a Dirichlet boundary condition and, hence, we refer to it as Dirichlet volume constraint.

Operators as in \eqref{original_operator} have been extensively analyzed in \cite{du2013nonlocal}, in which the authors introduce a nonlocal vector calculus that allows to express the operator $\mcL^\mub$ as a composition of a nonlocal divergence and a nonlocal gradient (just like in the local setting). They also introduce the concept of nonlocal {\em curl} and provide the nonlocal counterpart of several variational results of the classical local calculus (e.g., integration by parts, Green's identities, etc). In \cite{du2012analysis} the same authors use the nonlocal calculus to recast problems such as \eqref{eq:nonlocal-diffusion} in a variational setting and they study their well-posedness in relation to kernel properties. More specifically, they show that when sufficient conditions on the kernel function hold and nonlocal volume constraints are prescribed in $\omgd$, problem \eqref{eq:nonlocal-diffusion} is well-posed in a weak sense. We refer to \cite{du2012analysis} for more details on the functional space setting and properties of the solution in relation to the nonlocal kernel. Note that all kernels used in this paper are such that sufficient conditions for the well-posedness are satisfied.

\section{A unified nonlocal operator and its relationship with classical and fractional Laplacian}\label{sec:unified-operator}

In this section we introduce the unified nonlocal Laplace operator and discuss its behavior at the limit of vanishing or infinite interactions.

With the purpose of defining an operator that bridges the local and the fractional Laplacian, we consider the following kernel function parameterized by the interaction radius $\delta$ and the scalar $\alpha\in(0,2)$, we have
\begin{equation}\label{kernel}
\gamma(\xb,\yb;\mub:=[\delta,\alpha])=
\left\{\begin{array}{ll}
 \dfrac{C_{\delta,\alpha}}{\|\yb-\xb\|_2^{d+\alpha}}
& \|\yb-\xb\|_2\le \delta,\\ 
0, & \|\yb-\xb\|_2> \delta,
\end{array} \right. 
\end{equation}
where $\alpha\in (0,2)$ is a measure of how fast the nonlocal interaction decays as $\yb$ moves away from $\xb$. With this kernel, the nonlocal Laplace operator can be rewritten as
\begin{equation}\label{eq:param-L}
\mathcal{L}^{\delta,\alpha}u(\xb)=
C_{\delta,\alpha}\int_{B_{\delta}(\xb)}
\frac{u(\yb)-u(\xb)}{\|\yb-\xb\|_2^{d+\alpha}}d\yb,
\quad\forall\;\xb\in\Omega.
\end{equation}
The goal of the paper is to jointly estimate the $\delta$ and $\alpha$, given a limited number of observations for $u$ in $\omg\cup\omgd$ and for $f$ in $\omg$.

\smallskip
First, we recall some well-established results on the behavior of the operator in \eqref{eq:param-L} at the limit of vanishing and infinite interactions, for different kernel scaling constants. Let the scaling constant in \eqref{kernel} be defined as
\begin{equation*}
C_{\delta,\alpha}:=C'_{\delta,\alpha}=
\frac{2(2-\alpha)\Gamma(\frac{d}{2}+1)}{\pi^{d/2}\delta^{2-\alpha}}.
\end{equation*}
Then, by Taylor expansion, we have the following limit
\begin{equation}\label{lim1}
\lim\limits_{\delta\rightarrow 0}(-\mathcal{L}^{\delta,\alpha})
u(\xb)=-\Delta u(\xb), \quad \forall\;\alpha \in (0,2).
\end{equation}
Let the fractional Laplacian operator be defined as \cite{lischke2019fractional}
\begin{equation}\label{eq:frac-lapl}
(-\Delta)^{\alpha/2}u(\xb)=C''_{\delta,\alpha}\;\mbox{p.v.}\;\int_{\mbRd}\frac{u(\xb)-u(\yb)}{\|\xb-\yb\|_2^{d+\alpha}}d\yb,
\end{equation}
where ``p.v.'' denotes the principle value of the integral and $C''_{\delta,\alpha}$ is defined as
\begin{equation*}
C''_{\alpha}=\frac{2^{\alpha}\Gamma(\frac{d}{2}+\frac{\alpha}{2})}{\pi^{d/2}|\Gamma(-\frac{\alpha}{2})|}.
\end{equation*}
For $C_{\delta,\alpha}:=C''_{\delta,\alpha}$ in \eqref{kernel}, by simply applying the limit as $\delta\to\infty$ to the domain of integration, we have
\begin{equation}\label{lim2}
\lim\limits_{\delta\rightarrow\infty}
(-\mathcal{L}^{\delta,\alpha})u(\xb)=(-\Delta)^{\alpha/2} u(\xb), 
\quad \forall\; \alpha\in(0,2).
\end{equation}
These limits establish that when properly scaled and for two different scaling constants, the operator defined in \eqref{eq:param-L} converges to the classical Laplacian as the nonlocality vanishes ($\delta\to 0$) and to the fractional Laplacian as the interactions become infinite ($\delta\to\infty$). We propose a new single scaling constant that defines a unified nonlocal operator for which those limits hold true, as shown in the following lemma.
\begin{lemma}\label{lem:limiting-behavior}
For $d=1,2,3$ and 
$C_{\delta,\alpha}=C'_{\delta,\alpha}+C''_{\alpha}$
the nonlocal Laplace operator $\mathcal{L}^{\delta,\alpha}$ satisfies both (\ref{lim1}) and (\ref{lim2}).
\end{lemma}

\noindent{\it Proof.}
We split the proof in two parts. We focus on the case $d=1$, and prove the two- and three-dimensional cases in Appendix \ref{sec:AppA}.

\smallskip\noindent
I. The operator $\mathcal{L}^{\delta,\alpha}$ satisfies \eqref{lim1}.\\
We assume that $u(x)$ is $C^2(\omg\cup\omgd)$. 
We consider the following integral
\begin{equation}
    \int_{x-\delta}^{x+\delta}\frac{u(y)-u(x)}{|y-x|^{1+\alpha}}dy=\int_{-\delta}^{\delta}\frac{u(x+z)-u(x)}{|z|^{1+\alpha}}dz.
\end{equation}
Using Taylor series expansion at $x$ yields
\begin{equation}\label{ss}
  \begin{split}
   \int_{-\delta}^{\delta}\frac{u(x+z)-u(x)}{|z|^{1+\alpha}}dz & = \int_{-\delta}^{\delta}\frac{zu'(x)}{|z|^{1+\alpha}}dz + \int_{-\delta}^{\delta}\frac{z^2u''(x)}{2}\cdot\frac{1}{|z|^{1+\alpha}}dz + \int_{-\delta}^{\delta}\frac{o(|z|^2)}{|z|^{1+\alpha}}dz\\
   & = 0 + u''(x)\frac{\delta^{2-\alpha}}{2-\alpha}+\int_{-\delta}^{\delta}\frac{o(|z|^2)}{|z|^{1+\alpha}}dz  \\
   & = u''(x)\frac{\delta^{2-\alpha}}{2-\alpha}+\int_{-\delta}^{\delta}\frac{o(|z|^2)}{|z|^{1+\alpha}}dz \\
   & = u''(x)\frac{\delta^{2-\alpha}}{2-\alpha} + \mcO(\delta^\beta),
 \end{split}  
\end{equation}
where $\beta>2-\alpha$ and $u'(x)$ and $u''(x)$ are first and second derivatives. Here, the second equality follows from the fact that the integrand $\frac{z}{|z|^{1+\alpha}}$ is an odd function. Thus, we have
\begin{equation*}
\begin{split}
\lim\limits_{\delta\rightarrow0}C_{\delta,\alpha}
\int_{-\delta}^{\delta}\frac{u(x+z)- u(x)}{|z|^{1+\alpha}}dz & =\lim\limits_{\delta\rightarrow0} \left(\frac{2-\alpha}{\delta^{2-\alpha}}+\frac{2^{\alpha}\Gamma(\frac{1}{2}+\frac{\alpha}{2})}{\pi^{1/2}|\Gamma(-\frac{\alpha}{2})|}\right)
\left(u''(x)\frac{\delta^{2-\alpha}}{2-\alpha}+
\mcO(\delta^\beta)\right)\\
& = u''(x) + \lim\limits_{\delta\rightarrow 0} \frac{2^{\alpha}\Gamma(\frac{1}{2}+\frac{\alpha}{2})}{\pi^{1/2}|\Gamma(-\frac{\alpha}{2})|}u''(x)\frac{\delta^{2-\alpha}}{2-\alpha} = u''(x).
   \end{split}    
\end{equation*}

\smallskip\noindent  II. The operator $\mathcal{L}^{\delta,\alpha}$ satisfies \eqref{lim2}.\\
In one dimension, the fractional Laplacian is defined as
\begin{equation*}
    (-\Delta)^{\alpha/2}u(x)=C''_{\alpha}\int_{-\infty}^{\infty}\frac{u(x)-u(y)}{|y-x|^{1+\alpha}}dy.
\end{equation*}
We assume for simplicity that $u\equiv 0$ in $\mbR\setminus\omg$ (the non-homogeneous case can be treated in the same manner). By subtracting the nonlocal and fractional Laplacian we have
\begin{equation*}
\begin{split}
|-\mathcal{L}^{\delta,\alpha}u(x)&-(-\Delta)^{\alpha/2}u(x)|  =\left|-\frac{2-\alpha}{\delta^{2-\alpha}}\int_{-\delta}^{\delta}\frac{u(x+z)-u(x)}{|z|^{1+\alpha}}dz\right. \\
& \left. + C''_{\delta,\alpha}\int_{\delta}^{+\infty}\frac{u(x+z)-u(x)}{|z|^{1+\alpha}}dz 
+ C''_{\delta,\alpha}\int_{-\infty}^{-\delta}\frac{u(x+z)-u(x)}{|z|^{1+\alpha}}dz \right|.
\end{split}    
\end{equation*}
The first integral on the right-hand side is finite, and, thus, the integral converges to zero with rate $\alpha-2$ as $\delta\rightarrow\infty$. For $z\in(\delta,+\infty)\cup (-\infty,-\delta)$, we know that for large\footnote{At least larger than the size of the domain.} $\delta$, $u(x+z)=0$ because of $u\equiv 0$ in $\mbR\setminus\omg$, so that the second and third terms on the right-hand side become
\begin{equation}
C''_{\alpha}\int_{\delta}^{+\infty}
\frac{u(x+z)-u(x)}{|z|^{1+\alpha}}dz = 
- \frac{C''_{\delta,\alpha}u(x)}{\alpha}\delta^{-\alpha},
\end{equation}
and
\begin{equation}
C''_{\delta,\alpha}\int_{-\infty}^{-\delta}
\frac{u(x+z)-u(x)}{|z|^{1+\alpha}}dz = 
-\frac{C''_{\delta,\alpha}u(x)}{\alpha}\delta^{-\alpha},
\end{equation}
which converge to $0$ as $\delta\to\infty$. $\square$

Based on the proof above and the proofs in Appendix \ref{sec:AppA} for the multivariate case, we summarize the asymptotic behavior as follows:
\begin{equation}\label{asmp1}
|(-\mathcal{L}^{\delta,\alpha})u(\xb)-(-\Delta) u(\xb)|
\sim C_1(\alpha,\xb) \; \delta^{2-\alpha}
\quad {\rm as}\; \delta\rightarrow0,
\end{equation}
and
\begin{equation*}\label{asmp2}
|-\mathcal{L}^{\delta,\alpha}u(\xb)-(-\Delta)^{\alpha/2}u(\xb)|
\sim C_2(\alpha,\xb)\;\delta^{\alpha-2}
+    C_3(\alpha,\xb)\; \delta^{-\alpha}
\quad {\rm as}\; \delta\rightarrow\infty,
\end{equation*}
or, equivalently,
\begin{equation}\label{asmp22}
|-\mathcal{L}^{\delta,\alpha}u(\xb)-(-\Delta)^{\alpha/2}u(\xb)|
\sim C_4(\alpha,\xb)\;\delta^{\max\{\alpha-2,-\alpha\}}
\quad {\rm as}\; \delta\rightarrow\infty,
\end{equation}
where $C_i$, $i=1,\ldots 4$, are positive constants depending on $\alpha$ and $\xb$ only. 

These results show that the unified nonlocal operator in \eqref{eq:param-L} is a superset of the classical and fractional Laplacian, and, therefore, we expect the operator to be able to fit a broad spectrum of experimental data. As an example, in a stochastic context, $\mcL^{\delta,\alpha}$ acts as a bridge between the generator of a Gaussian process (the classical Laplacian) and the one of an isotropic $\alpha$- stable L\'{e}vy process (the fractional Laplacian) \cite{lischke2019fractional}. 

Note that tempered fractional operators, also related to stochastic processes \cite{sabzikar2015tempered}, do not belong to the superset spanned by the unified operator in \eqref{eq:param-L}. However, truncating has a similar effect as tempering, i.e., the second moment of truncated operators, such as ours, is finite, as for the tempered case. This is very important when modeling engineering applications such as turbulence, see Section \ref{sec:turbulence}. Also, the extension to generalized tempered operators only requires multiplication of the kernel function in \eqref{kernel} by the factor $\exp\{-\lambda\|\xb-\yb\|\}$, for some $\lambda>0$. Hence, adapting nPINNs to this class of operators can also be implemented.

\section{Nonlocal Physics-Informed Neural Networks (nPINNs)} \label{sec:nPINNs-formulation}

In this section we first provide an abstract formulation of the nPINNs algorithm and then describe in detail specific components, namely fully connected NNs, and discretization, evaluation, and minimization of the loss function.

The main idea of nPINNs can be summarized in three simple steps.
\begin{itemize}
\item[\bf 1.] Collect observations or high fidelity simulations of the solution, $u_{obs}$;
\item[\bf 2.] Approximate the solution with a fully-connected NN: $u(x) \cong \unn(\xb;\thetab)$;
\item[\bf 3.] Minimize the loss function
    \begin{equation}\label{loss_fun}
    \mbox{Loss}(\thetab,\delta,\alpha)  = \frac{\|-\mathcal{L}^{\delta,\alpha}
    \unn(\xb;\thetab)-f(\xb)\|^2}{\|f(\xb)\|^2} 
   + \frac{\|\unn(\xb;\thetab)
   -u_{obs}(\xb)\|^2}{\|u_{obs}(\xb)\|^2}, 
   \end{equation}
with respect to the unknown parameters $(\alpha,\delta)$ and the NN parameters $\thetab$, which are, in turn, the outcome of the minimization.
\end{itemize}

Being this an abstract formulation, we do not specify the norms used in \eqref{loss_fun} and how they are computed. Note that ``Loss'' has a physics-driven and a data-driven component: the first term controls the residual of the nonlocal equation, whereas the second the mismatch between solution and data. The residual term acts as a regularizer and improves the conditioning of the minimization problem. We point out that algorithm {\bf 1.}--{\bf 3.} can be used for both the solution of forward and inverse problems; in the first case, $u_{obs}$ is only needed in the interaction domain, whereas, in the second case, additional values of the solution inside the domain improve the learning process, see Section \ref{sec:eval-min-loss} for a detailed explanation.

As described in \cite{wang2020understanding}, both terms in ``Loss'' can be multiplied by penalization constants, which can be dynamically updated throughout the optimization process to speed up the convergence but we do not employ this technique in the current work.

\begin{figure}[H] 
\centering
\includegraphics[width=.9\textwidth]{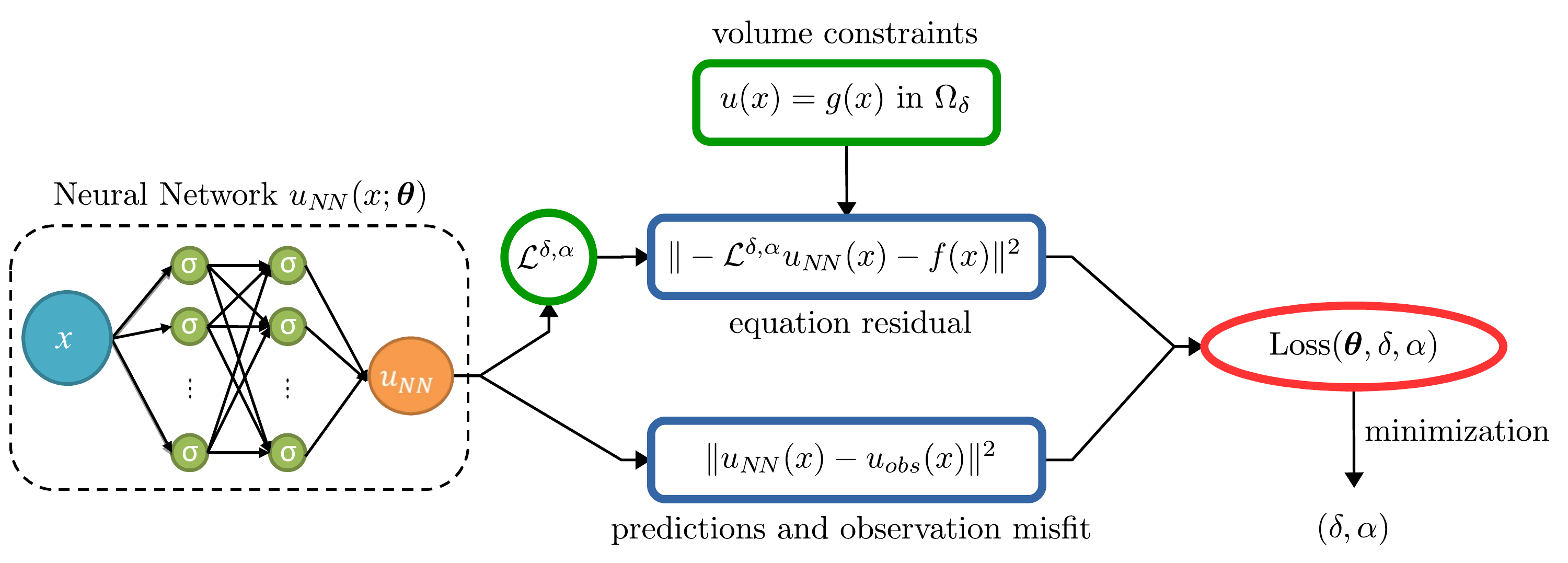}
\caption{\label{nPINN_illustration} Schematic of nPINNs for estimating parameters $(\delta,\alpha)$.}
\label{domain}
\end{figure}

In Figure \ref{nPINN_illustration} we present an illustration of {\bf 1.}--{\bf 3.}: on the left we sketch  a fully-connected NN; in the center, the blue boxes are the residual and the misfit terms. Their combination results in the loss function ``Loss'', whose minimization delivers the optimal parameters $(\delta,\alpha)$. The green boxes represent the model and the nonlocal boundary conditions.

\subsection{Fully-connected neural networks}

We approximate the solution of \eqref{eq:nonlocal-diffusion} as a fully-connected NN, $\unn$, consisting of a series of linear and nonlinear transformations. Let the spatial coordinate $\xb$ be the input of the NN and denote by $\thetab$ the collection of transformation parameters, so that $\unn= \unn(\xb;\thetab)$. A NN of depth $n+1$ consists of an input and an output layer and $n$ hidden layers; furthermore, each layer has a certain number of neurons, which we call ``width". The input layer has width $d_0$, the output layer $1$ and the hidden layers $d_i$, $i=1,\ldots n$; see Figure \ref{nPINN_illustration} on the left for an illustration.
The output of the $i$-th layer is denoted by $\xb_i\in \mbR^{d_i}$. The weight matrix $\mathbf{W}_i\in \mathbb{R}^{d_{i+1}\times d_i}$ and the bias vector $\mathbf{b}_i \in \mathbb{R}^{d_{i+1}}$ link the outputs of the $i$- and $(i+1)$- layers by the recurrence relation
\begin{equation}
\left\{
\begin{array}{ll}
\xb_0 = \xb, &   \\
\xb_{i+1} = \sigma(\mathbf{W}_i \xb_i + \mathbf{b}_i), & i=0, 1, \cdots, n-1,\\
\unn(\xb;\thetab) = \mathbf{W}_n \xb_n + \mathbf{b}_n, 
& \thetab=\{\{\mathbf{W}_i\},\{\mathbf{b}_i\}\}_{i=0}^n.
\end{array} \right.
\end{equation}
In other words, the output of each layer is first subject to a linear transformation parameterized by $\mathbf{W}_i$ and $\mathbf{b}_i$, and then a nonlinear transformation $\sigma(\cdot)$, which is an element-wise nonlinear function also known as \textit{activation function}\footnote{In this work we consider bounded and infinitely differentiable functions $\sigma$, namely the hyperbolic tangent. Sometimes we need to restrict the magnitude of the output $\unn$. One approach is to consider an extra activation function $\kappa(\cdot)$ in the output layer, namely, $\unn(\mathbf{x},\boldsymbol{\theta})=\kappa(\mathbf{W}_n\mathbf{x}_n+\mathbf{b}_n)$. For example, in Section \ref{sec:turbulence} we take $\kappa(\cdot)$ in $\alpha_{N\!N}$ to a sigmoid function in order to ensure the fractional order ranges from 0 to 1. Another example is to take $\kappa(\cdot)$ as a softmax function in the context of image classification. } Note that the output of the last hidden layer enters the output layer without nonlinear transformation $\sigma(\cdot)$. The vector $\thetab$ collects all the undetermined NN parameters $\mathbf{W}_i$ and $\mathbf{b}_i$; as such, these parameters are part of the outcome of the minimization algorithm.
The depth $n+1$ and the width $d_i$ are called hyper-parameters of the NN; note that the performance of nPINNs, specifically the convergence of the minimization algorithm could be sensitive to these hyper-parameters. Although there has been recent work on the automatic selection of these parameters \cite{zoph2016neural,finn2017model}, in this work we pick them heuristically.

Note that the NN is a global function defined over $\mbRd$; however, we only evaluate it within $\omg$ as its values outside the domain are determined by nonlocal boundary data.

\subsection{Evaluation and minimization of the loss function}\label{sec:eval-min-loss}

We describe how to evaluate the Loss function at given values of $\thetab$, $\delta$ and $\alpha$. While a standard choice for the norm of any function $v$ defined over $\omg$ would be  $\|v\|_{L^2(\omg)}^2:=\int_{\Omega}v^2(\xb)d\xb$, in practice, it is almost impossible to know the values of $v(\xb)$ everywhere in the domain. Thus, we consider the $\ell^2$ norm of the vector of values of $v$ at specific points $\xb_j\in\omg$, $j=1,\ldots J$, i.e. $\|v\|_{\ell^2}^2=\sum_{j=1}^J v^2(\xb_j)$.

It follows that the discrete version of the loss function in \eqref{loss_fun}, the empirical Loss, which for convenience we still denote by ``Loss'', is given by
\begin{equation}\label{loss_fun_discrete}
\mbox{Loss}(\thetab,\delta,\alpha)  =
\frac{\sum_{k=1}^N\left(-\mcL^{\delta,\alpha}\unn(\xb_k;\thetab)
- f(\xb_k)\right)^2}{\sum_{k=1}^N\left(f(\xb_k)\right)^2} 
+ \frac{\sum_{k=1}^{N_{obs}}(\unn\left(\hat{\xb}_k;\thetab)
- u_{obs}(\hat{\xb}_k)\right)^2}{\sum_{k=1}^{N_{obs}}
  (u_{obs}(\hat{\xb}_k))^2}.
\end{equation}
Here, we introduce two sets of points: the {\it residual} points $\{\xb_k\}_{k=1}^N$ and the {\it observation} points $\{\hat{\xb}_k\}_{k=1}^{N_{obs}}$, which do not necessarily coincide. Both sets can either be equally spaced or scattered in the computational domain $\Omega\cup\omgd$. However, the location of the observation points is generally dependent on the experimental setup\footnote{As an example, for groundwater solute transport, the location depends on where the monitoring wells are dug in the field experiment.}. 
Furthermore, depending on whether we are solving a forward or inverse problem, observations of $u$ may not be needed inside the domain $\omg$. Specifically, if we are solving a forward problem, it is enough to have observation points in $\omgd$; instead, if we are solving an inverse problem, additional observation points within $\omg$ improve the learning process. In the case where the nonlocal volume constraint $g$ is available everywhere in $\omgd$, it can be automatically prescribed while evaluating the nonlocal operator and observations in $\omgd$ are not needed. In this work, we assume that $g$ is a known function; thus, when using nPINNs for the solution of forward problems, we do not need any observations (neither in $\omg$ nor in $\omgd$) and the loss function reduces to
\begin{equation}\label{loss_fun_forward}
\mbox{Loss}(\thetab)  = \frac{\sum_{k=1}^N\left(-\mathcal{L}^{\delta,\alpha}\unn(\xb_k;\thetab)-f(\xb_k)\right)^2}{\sum_{k=1}^N\left(f(\xb_k)\right)^2}. 
\end{equation}
Note that, in this case, the outcome of the minimization of \eqref{loss_fun_forward} are the NN parameters $\thetab$ only.

\smallskip
The evaluation of the nonlocal operator in \eqref{loss_fun_discrete} at residual points can be performed by using any quadrature technique\footnote{Note that in PINNs for PDEs, we can adopt the chain rule when doing automatic differentiation to analytically compute partial derivatives of $\unn$ with respect to input parameters. However, in nPINNs (as well as fPINNs), the chain rule does not apply and we have to first approximate the integration, hence introducing approximation error.}; in this work we adopt the composite Gauss-Legendre quadrature combined with singularity subtraction, see Appendix \ref{sec:AppB}, where we describe quadrature rules for the one- and multi-dimensional cases.

\smallskip
Several algorithms can be employed for the minimization of the loss function; a popular strategy is to use gradient-based optimization algorithms such as the conjugate gradient method, the limited memory BFGS algorithm, and the Adam method \cite{kingma2014adam}. We adopt the last one, since it is extensively utilized for machine learning tasks and it needs less tuning on its controlling parameters. It should be noted that the loss function is strongly non-convex and, thus, locating its global minimum is very hard. However, as illustrated in Section \ref{sec:numerical-tests}, for several local minima, values of the loss function are comparable and the accuracy of the corresponding solutions (when a reference solution is available) is similar. These cases are instances of the so-called \textit{operator mimicking} phenomenon, where the strong non-convexity of the loss function may lead to multiple good local minima; i.e., there exist multiple ($\delta$, $\alpha$) pairs that generate distinct operators which are equally effective in reproducing the training data.

\section{Computational results}\label{sec:numerical-tests}

This section consists of three parts. In Section \ref{limit-sec}, we illustrate the behavior of the unified operator \eqref{eq:param-L} and the solution of \eqref{eq:nonlocal-diffusion} with respect to $\delta$ and $\alpha$, and then present a numerical study of the limit behavior as the nonlocal interaction radius $\delta$ goes to zero and infinity. In Section \ref{forward-sec}, we demonstrate the accuracy of nPINNs through the solution of one-, two-, and three-dimensional forward problems. In Section \ref{inverse-sec}, we report results of the estimation of $\delta$ and $\alpha$ via nPINNs.

Unless otherwise stated, $\unn$ has $n=4$ hidden layers with constant width $d_i \equiv 10$. The Xavier \cite{glorot2010understanding} and zero initializations are utilized to initialize weights $\mathbf{W}_i$ and biases  $\mathbf{b}_i$, respectively. The set of residual points $\{\xb_k\}_{k=1}^N$ is generated via Sobol sequences \cite{sobol1967distribution} and, when solving inverse problems, observation points are uniformly distributed in $\omg$. When a manufactured solution $u$ is available, we quantify the accuracy of $\unn$ with respect to $u$ by using the following metric: given $N_t$ test points $\{\xb^t_j\}_{j=1}^{N_t}$, also generated via Sobol sequences, we define the relative error $\epsilon$ as
\begin{equation}\label{def_err}
\epsilon^2=\frac{\sum_{j=1}^{N_t}
\left(u(\xb^t_j)-\unn(\xb^t_j;\thetab)\right)^2}
{\sum_{j=1}^{N_t}u(\xb^t_j)^2}.
\end{equation}

\subsection{Limit behavior of the parametrized operator}\label{limit-sec}

We illustrate the behavior of the unified nonlocal operator with respect to $\delta$ and $\alpha$ and its relation to the classical and fractional Laplacian.

We consider problem \eqref{eq:nonlocal-diffusion} in $\omg=(0,1)$ with a smooth source term $f(x)=\sin(2\pi x)$ and a zero volume constraint $g(x)=0$. We solve the forward problem using nPINNs for different values of interaction radius $\delta$ and decay rate $\alpha$. The quadrature parameters are set to $m=5$ and $M=10$ for all numerical examples (see Appendix \ref{sec:AppB}). Note that the purpose of the experiments presented in this section is not to test the performance of nPINNS; as such, these computational tests could be performed with any forward nonlocal discretization technique.

\begin{figure}[H]
\centering
\subfloat[]{
\includegraphics[width=.48\textwidth]{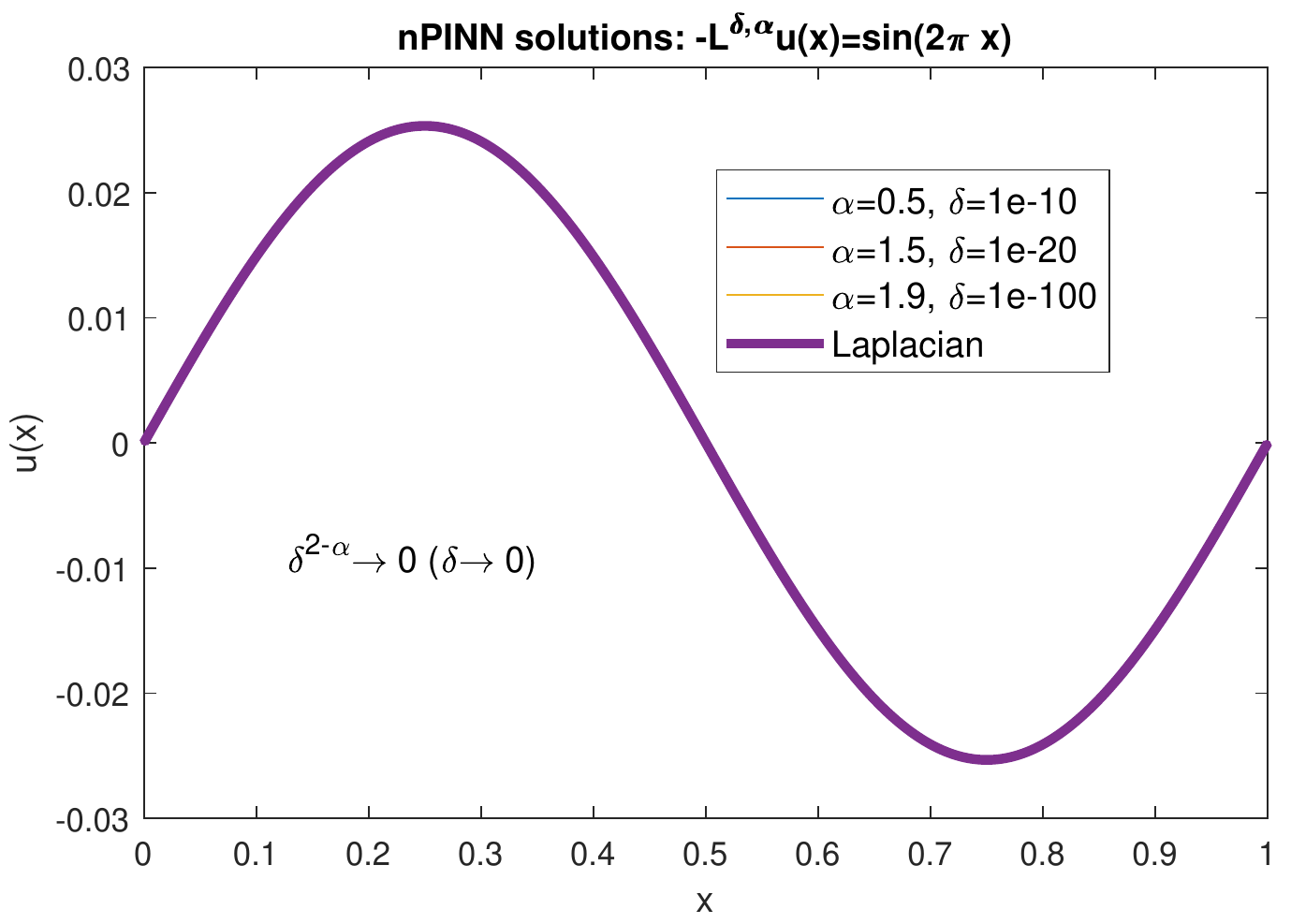}}
\subfloat[]{
\includegraphics[width=.48\textwidth]{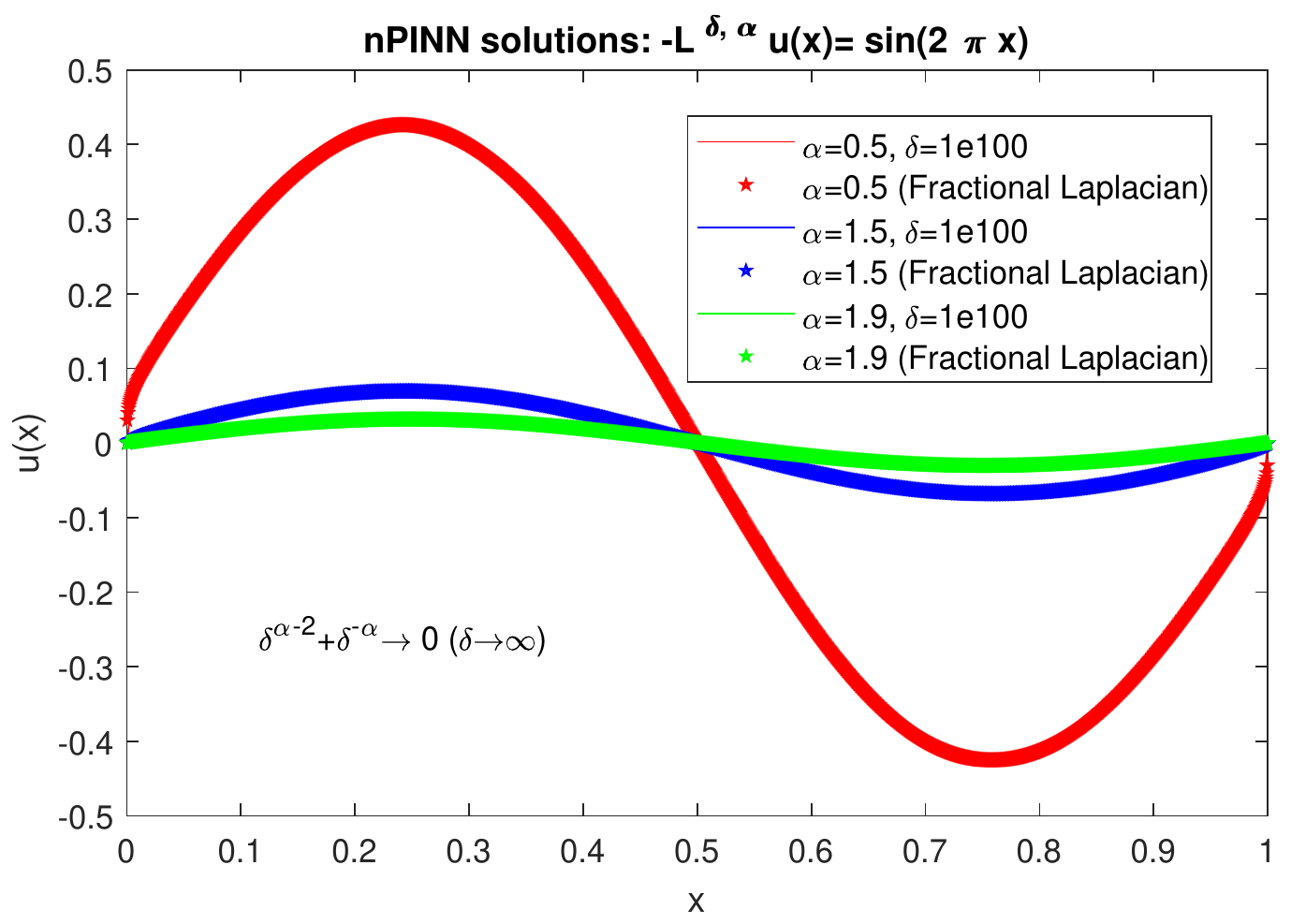}}
\caption{Limit behavior of the nonlocal operator $-\mathcal{L}^{\delta,\alpha}$ with respect to the nonlocal interaction radius $\delta$. (a) $\mcL^{\delta,\alpha}$ reduces to the classical Laplacian for small $\delta$ regardless of $\alpha$; (b) $\mcL^{\delta,\alpha}$ reduces to the fractional Laplacian for large $\delta$.  }
\label{delta_lim}
\end{figure}

\begin{figure}[H] 
\centering
\includegraphics[width=.5\textwidth]{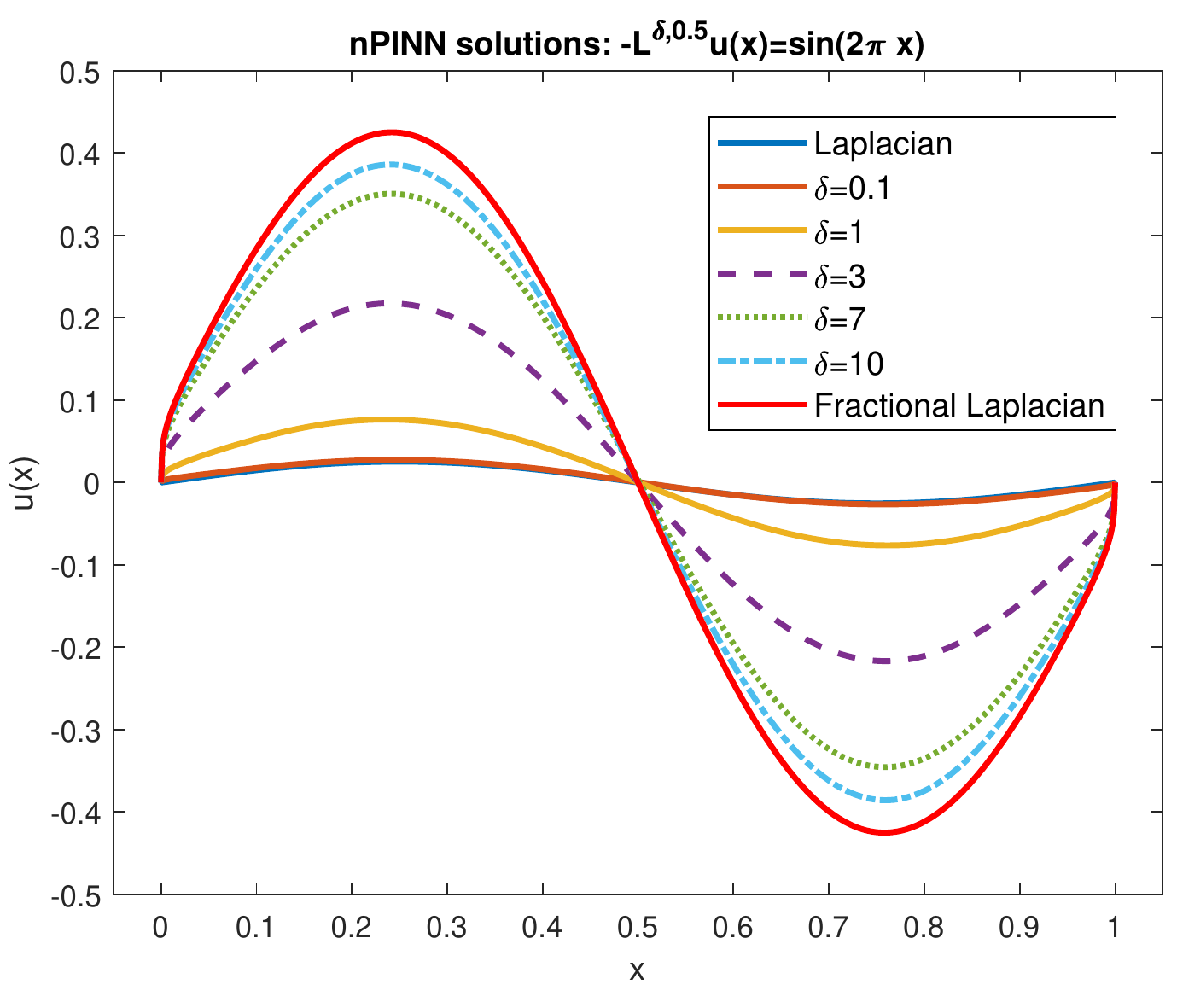}
\caption{Transition from classical to fractional limits for the solution of \eqref{eq:nonlocal-diffusion} as $\delta$ varies from zero to infinity for fixed $\alpha=0.5$.}
\label{var_delta}
\end{figure}

We report the results in Figure \ref{delta_lim}. On the left, our results confirm that the unified nonlocal operator reduces to the classical Laplacian for small $\delta$ regardless of $\alpha$. The curve corresponding to the classical Laplacian is computed analytically; we have $u(x)=\sin(2\pi x)/(4\pi^2)$. On the right, our tests show that the unified nonlocal operator reduces to the fractional Laplacian for large $\delta$ and different values of $\alpha$. The curves corresponding to the fractional Laplacian are generated with the third-order Gr\"unwald-Letnikov scheme (a finite difference scheme, see \cite{zhao2015series}). These results confirm our statements in Lemma \ref{lem:limiting-behavior}.
In Figure \ref{var_delta}, for the same nonlocal Poisson problem and fixed $\alpha=0.5$ we show that the unified operator is indeed a bridge between the classical and the fractional Laplacian. As the interaction radius transitions from zero to infinity, solutions of problem \eqref{eq:nonlocal-diffusion} span the whole spectrum of operators between the classical and fractional limit. Additionally, by letting the decay parameter change, we can model an even broader range of solutions; this makes our operator suitable for the identification of new complex nonlocal Laplace operators.

\paragraph{Convergence tests}

We illustrate the convergence rates determined in Section \ref{sec:unified-operator}. We introduce the relative error between the unified operator and classical and fractional Laplacian operators, respectively:
\begin{equation}\label{err_f}
\begin{split}
\epsilon_{clas}^2& =\frac{\sum_{j=1}^{N_t}(-\mathcal{L}^{\delta,\alpha}
u(x_j^t)-(-\Delta u(x_j^t)))^2}
{\sum_{j=1}^{N_t}(-\Delta u(x_j^t))^2}, \\[2mm]
\epsilon_{frac}^2& =\frac{\sum_{j=1}^{N_t}(-\mathcal{L}^{\delta,\alpha}
u(x_j^t)-(-\Delta)^{\alpha/2}u(x_j^t))^2}
{\sum_{j=1}^{N_t}((-\Delta)^{\alpha/2}u(x_j^t))^2},
\end{split}   
\end{equation}
where the test points $x_j^t$ are the same as those in \eqref{def_err}. We consider manufactured solutions $u(\cdot)$ and evaluate $\epsilon_{clas}$ and $\epsilon_{frac}$ assuming $u$ is known beforehand; for the classical case, we use $u(x)=\sin(2\pi x)/(4\pi^2)$ whereas for the fractional case we use $u(x)=x(1-x^2)^{1+\alpha/2}$. The action of $(-\Delta)$ and $(-\Delta)^{\alpha/2}$ on $u$ is computed analytically, while the one of $-\mathcal{L}^{\delta,\alpha}$ is computed numerically using the composite Gauss quadrature with $m=M=50$.

\begin{figure}[H]
\centering
\subfloat[$\epsilon_{clas}\sim \delta^{2-\alpha}$]{
\includegraphics[width=.48\textwidth]{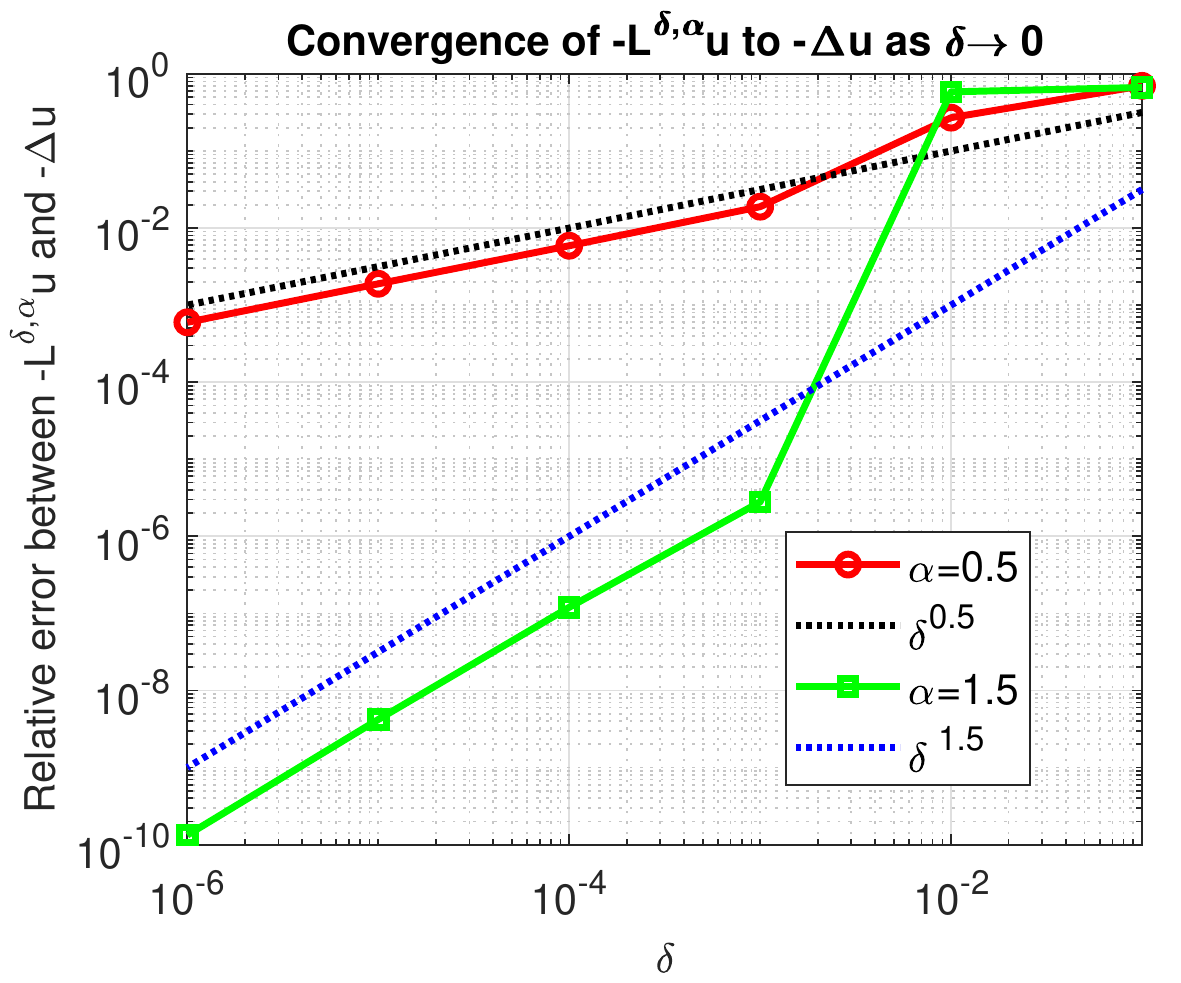}}
\subfloat[$\epsilon_{frac}\sim \delta^{\max\{-\alpha,\alpha-2\}}$]{
\includegraphics[width=.48\textwidth]{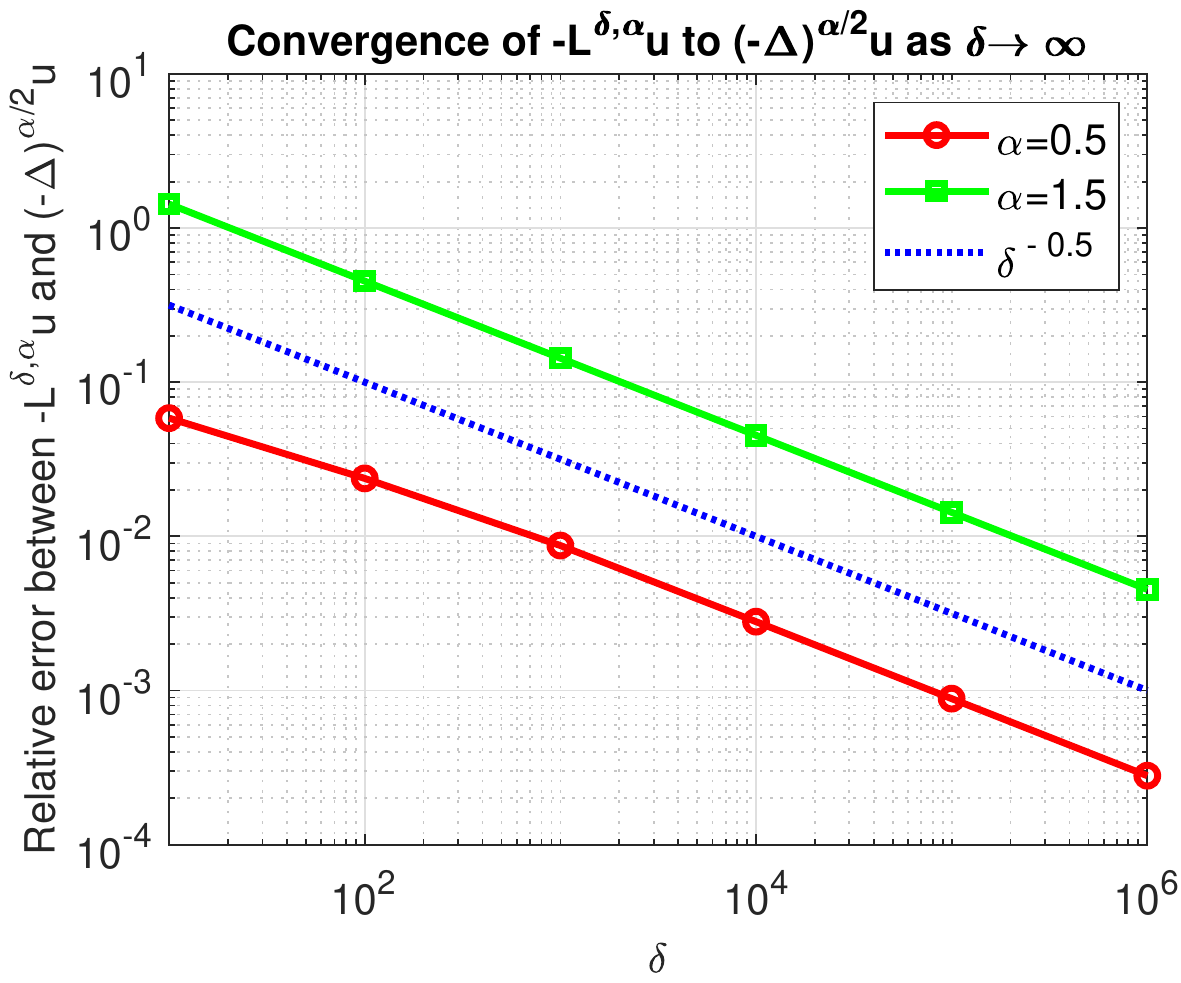}}
\caption{Convergence rates of the relative errors $\epsilon_{clas}$ (a) and $\epsilon_{frac}$ (b)  with respect to $\delta$: the numerical tests confirm our theoretical results.}
\label{delta_lim1}
\end{figure}

\begin{figure}[H]
\centering
\subfloat[$\epsilon_{clas,u} \sim \delta^{2-\alpha}$]{
\includegraphics[width=.48\textwidth]{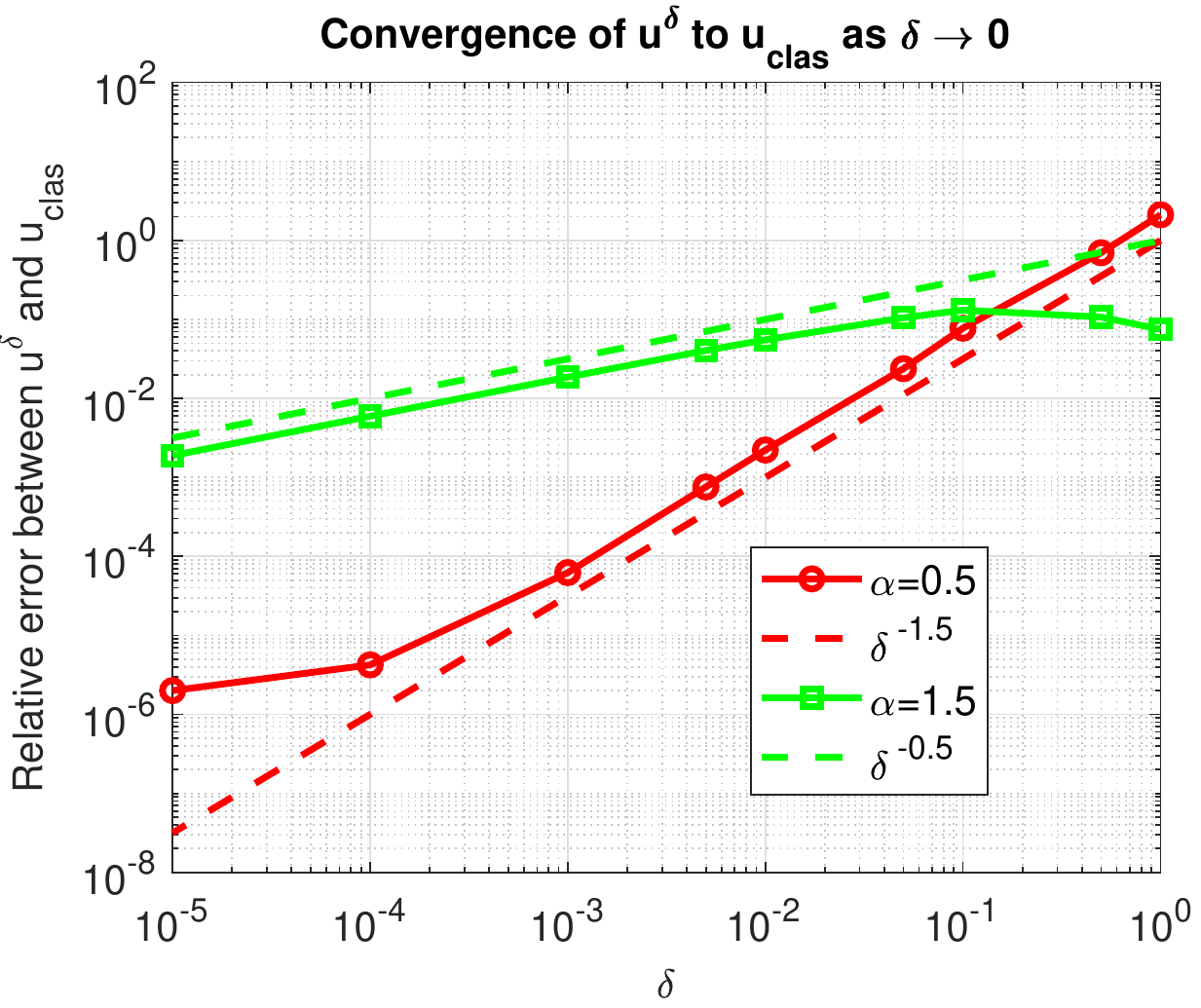}}
\subfloat[$\epsilon_{frac,u} \sim \delta^{\max\{-\alpha,-\alpha+2\}}$]{
\includegraphics[width=.48\textwidth]{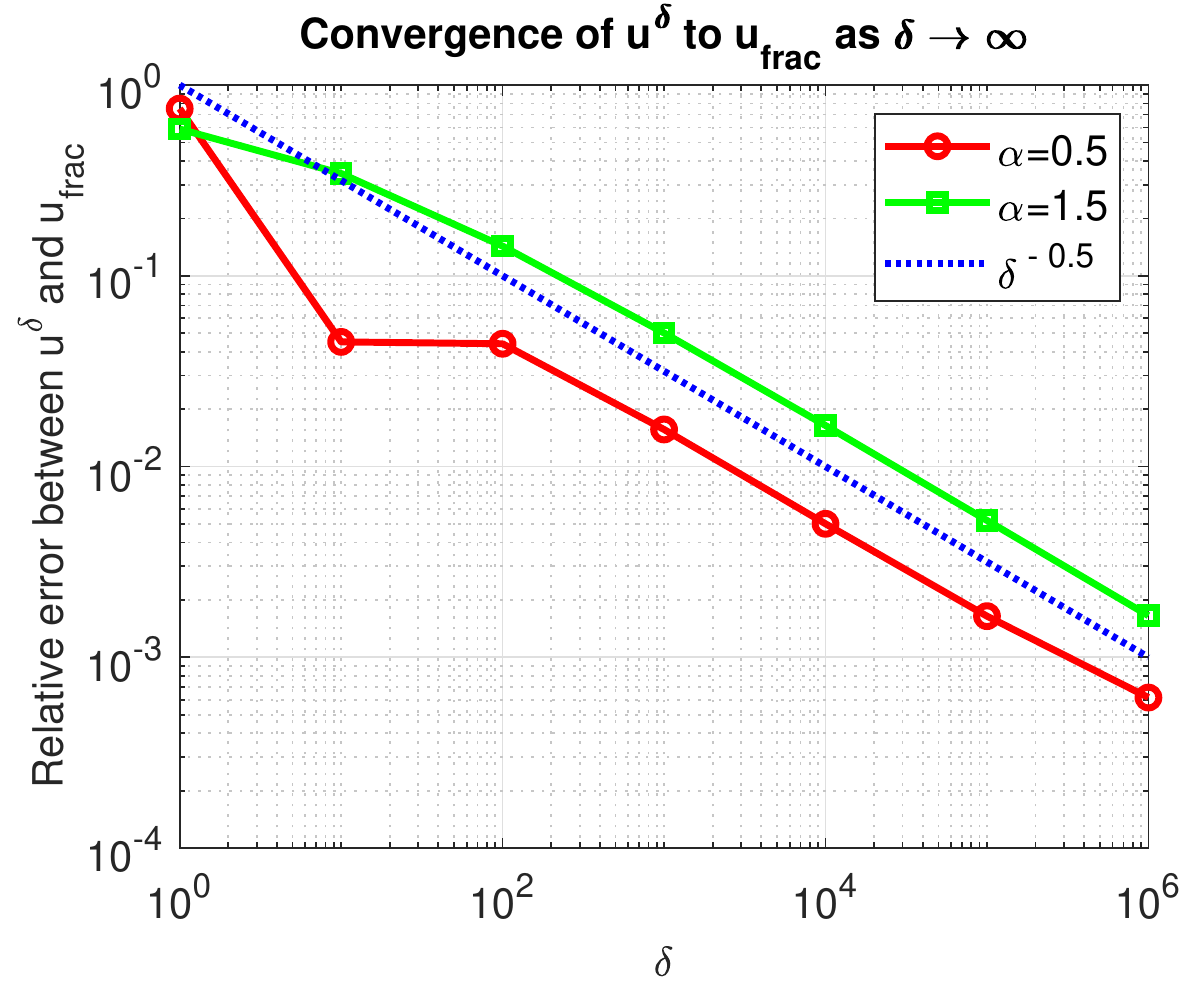}}
\caption{Convergence of nPINN solutions with respect to $\delta$ to the solution of the classical Laplacian as $\delta\to 0$ (a) and the fractional Laplacian as $\delta\to \infty$ (b). }
\label{delta_lim2}
\end{figure}

Figure \ref{delta_lim1} displays the convergence rates in log-log plots. The graph on the left shows that, as $\delta\rightarrow 0$, the relative error decays with the expected rate $2\!-\!\alpha$. From the graph on the right, we observe that when $\delta\rightarrow \infty$, the relative error decays with the expected rate $\max\{\alpha\!-\!2, -\!\alpha\}$.  
Even though not discussed in our theoretical section, we investigate experimentally the convergence of (forward) nPINNs solutions to the solution of the classical and fractional Laplacian equations. We denote the solutions of the classical and fractional equations as $u_{clas}$ and $u_{frac}$ respectively, and introduce the relative errors
\begin{equation}
\begin{split}
\epsilon^2_{clas,u}&=\frac{\sum_{j=1}^{N_t}
(\unn(x_j^t)-u_{clas}(x_j^t))^2}
{\sum_{j=1}^{N_t}(u_{clas}(x_j^t))^2} \\[2mm]
\epsilon^2_{frac,u}&=\frac{\sum_{j=1}^{N_t}
(\unn(x_j^t)-u_{frac}(x_j^t))^2}
{\sum_{j=1}^{N_t}(u_{frac}(x_j^t))^2}.
\end{split}    
\end{equation}
We consider two data sets. For the classical case: $\alpha=0.5$ and 1.5, $f(x)=\sin(2\pi x)$, and $g(x)=0$, for which $u_{clas}(x)=\sin(2\pi x)/(4\pi^2)$; and for the fractional case: $\alpha=0.5$ and 1.5, $f(x)=\frac{\Gamma(\alpha+3)}{6}\left(3-(3+\alpha)x^2\right)x$, and $g(x)=0$, for which $u_{frac}(x)=x(1-x^2)^{1+\alpha/2}$. Figure \ref{delta_lim2} shows the behavior of the relative errors of the solution with respect to $\delta$. We observe that the convergence rates are the same as those for the continuous operators, namely $\epsilon_{clas,u}\cong \mcO(\delta^{2-\alpha})$ and $\epsilon_{frac,u}\cong \mcO(\delta^{\max\{\alpha-2,-\alpha\}})$.

\subsection{Accuracy of nPINNs for forward problems} \label{forward-sec}
In this subsection, we first show that nPINNs for the solution of forward problems are as accurate as standard discretization methods for nonlocal equations when the number of nPINNs {\it residual} points, $N$, equals the number of {\it discretization} points. This is a very important property that shows consistency of our algorithm. Second, we illustrate the behavior of forward nPINNs in case of discontinuous solutions. Finally, we apply the forward nPINNs to two- and three-dimensional problems. These results show the applicability of our algorithm in higher dimensions.

\paragraph{Convergence with respect to the number of residual points}
For this task we focus on the pure fractional case, i.e. $\delta=\infty$. We consider the following one-dimensional manufactured solution in $\omg=(-1,1)$ for $\alpha=1.5$,
\begin{equation}\label{solu11}
   u(x)=\left\{
   \begin{array}{cc}
       x(1-x^2)^{1+\alpha/2},& \quad |x|\le 1 \\
       0, & \quad |x|>1.
   \end{array}  \right.
\end{equation}
The associated right-hand side is given by
\begin{equation}
f(x)  = (-\Delta)^{\alpha/2}u_{1D}(x)=\frac{\Gamma(\alpha+3)}{6}(3-(3+\alpha)x^2)x. 
\end{equation}
We solve this problem using three methods: (1) nPINNs with $\alpha=1.5$ and $\delta=10^{100}$ (to mimic the fractional limit), (2) fPINNs, see \cite{pang2019fpinns}, with $\alpha=1.5$, and (3) second-order Gr\"unwald-Letnikov, a standard finite difference discretization of the strong form of fractional equations. We take the same depth and width of NNs for both the nPINNs and fPINNs, i.e. $n=4$ and $d_i \equiv 10$. To minimize the quadrature error in nPINNs, we use the quadrature parameters $m=10$ and $M=50$.

\begin{figure}[H] 
\centering
\includegraphics[width=.6\textwidth]{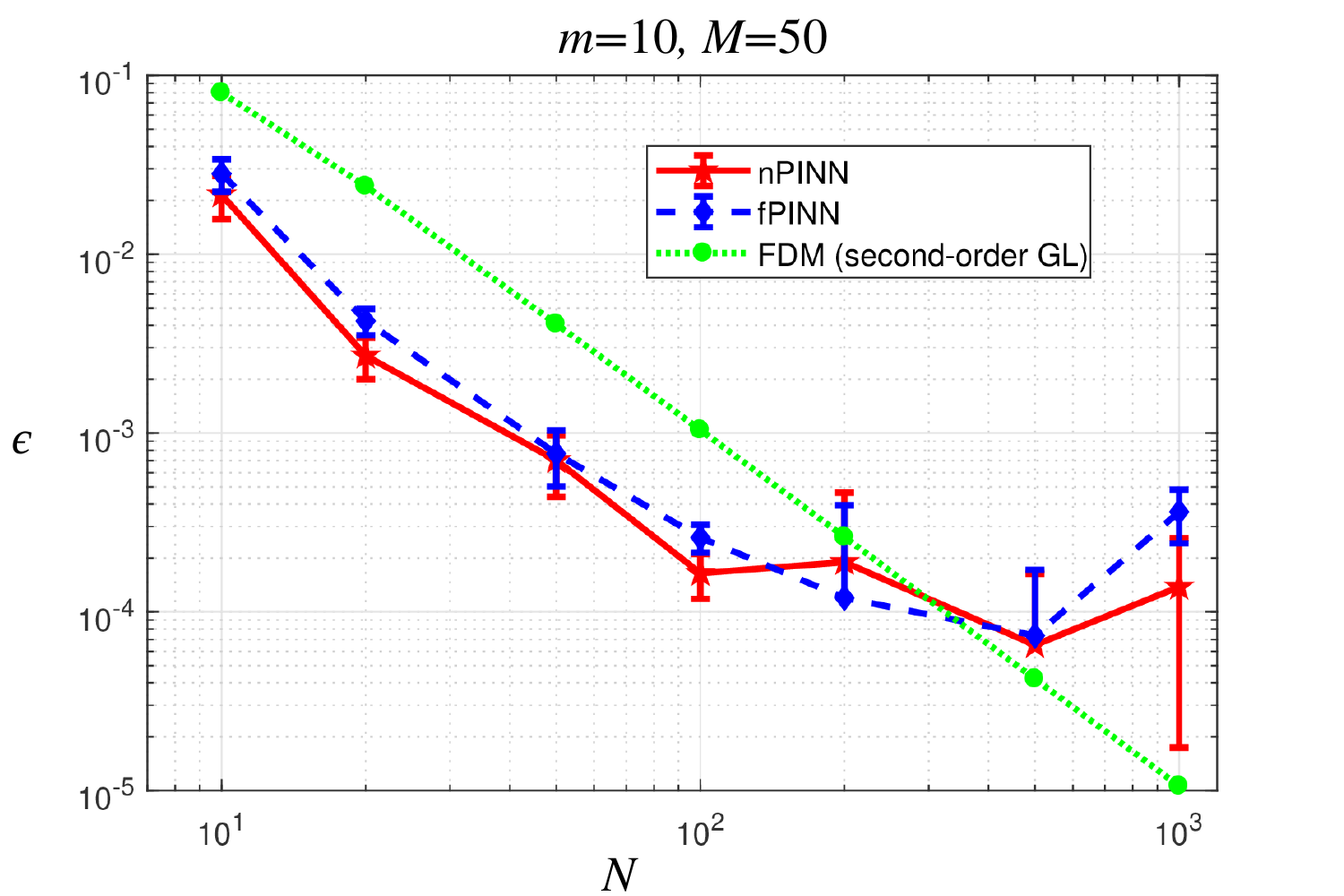}
\caption{Convergence of the relative error $\epsilon$ for nPINNs, fPINNs, and the finite difference scheme with respect to $N$ (number of residual points for nPINNs and fPINNs and number of discretization points for the finite difference scheme, the two sets coincide). Due to the sensitivity to the initial guess for $\thetab$, nPINNs and fPINNs curves are the result of ten simulations. The markers indicate the mean value of $\epsilon$ and the bars the standard deviation $\sigma$ (the length of the bar is $4\sigma$).}
\label{cmp-N}
\end{figure}

We compare the convergence of the relative error $\epsilon$ defined in \eqref{def_err} with respect to number, $N$, of residual points; these are uniformly distributed in $\Omega = (-1,1)$. To have a fair comparison, we let the discretization points of the finite difference scheme and the residual points of the nPINNs algorithm coincide. Results are reported in Figure \ref{cmp-N}: we observe that nPINNs and fPINNs have a very similar convergence behavior, which confirms that nPINNs reduces to fPINNs for $\delta=\infty$. Furthermore, for small $N$, the convergence behavior is the same as the one of the finite difference scheme and the error values are even lower. 
On one hand, increasing the number of residual points allows for a more accurate evaluation of the residual $\|-\mathcal{L}^{\delta,\alpha}\unn(x;\thetab)-f(x)\|^2$; on the other hand, increasing $N$ after a threshold makes optimization more difficult and, hence, less effective, resulting in error stagnation.

\begin{figure}[H] 
\centering
\includegraphics[width=.6\textwidth]{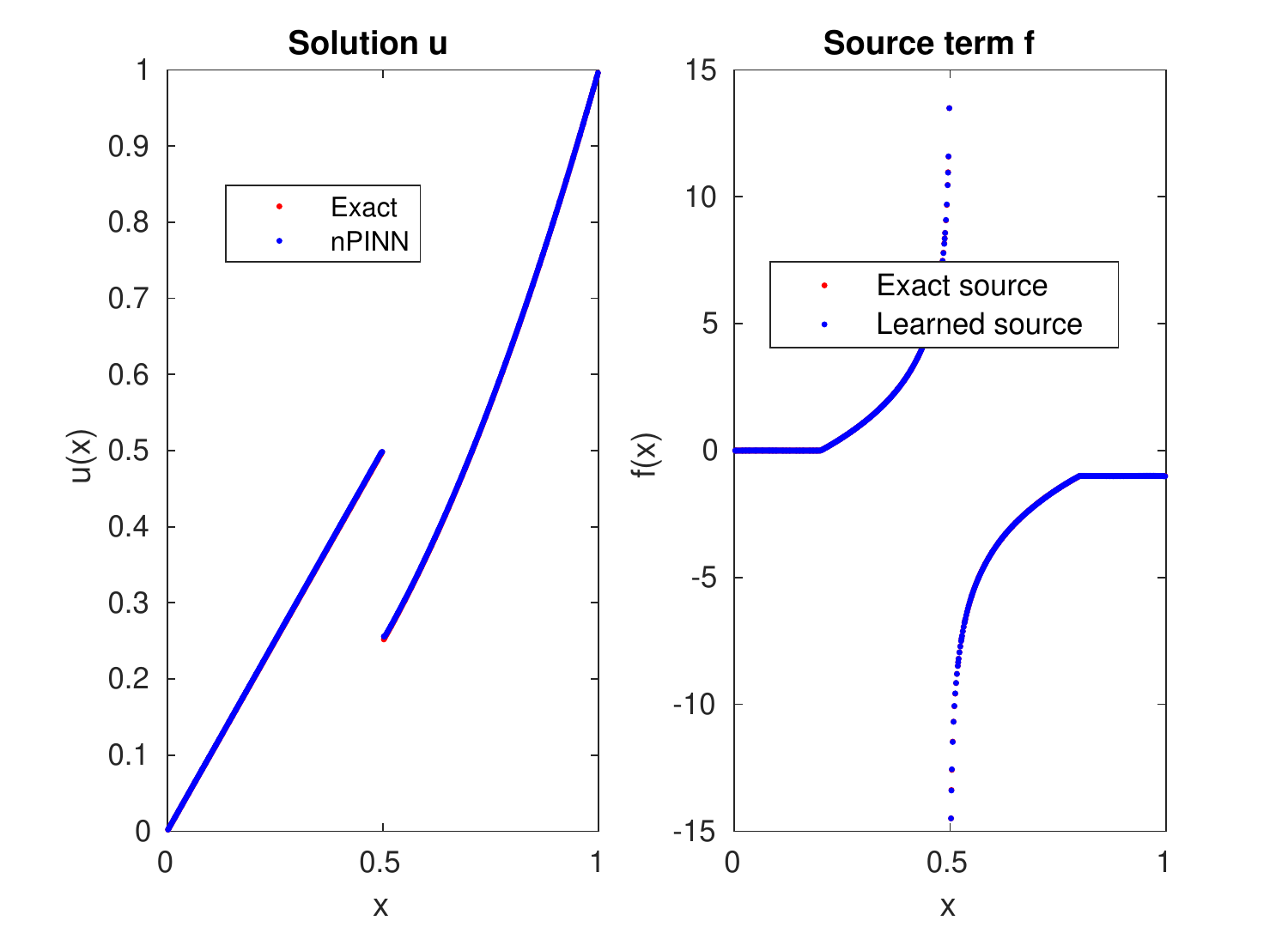}
\caption{For problem (I): nPINNs solution and exact solution on the left; trained $-\mcL^{\delta,\alpha}\unn$ and exact $f$ on the right. }
\label{1d-dis}
\end{figure}

\begin{figure}[H] 
\centering
\includegraphics[width=.6\textwidth]{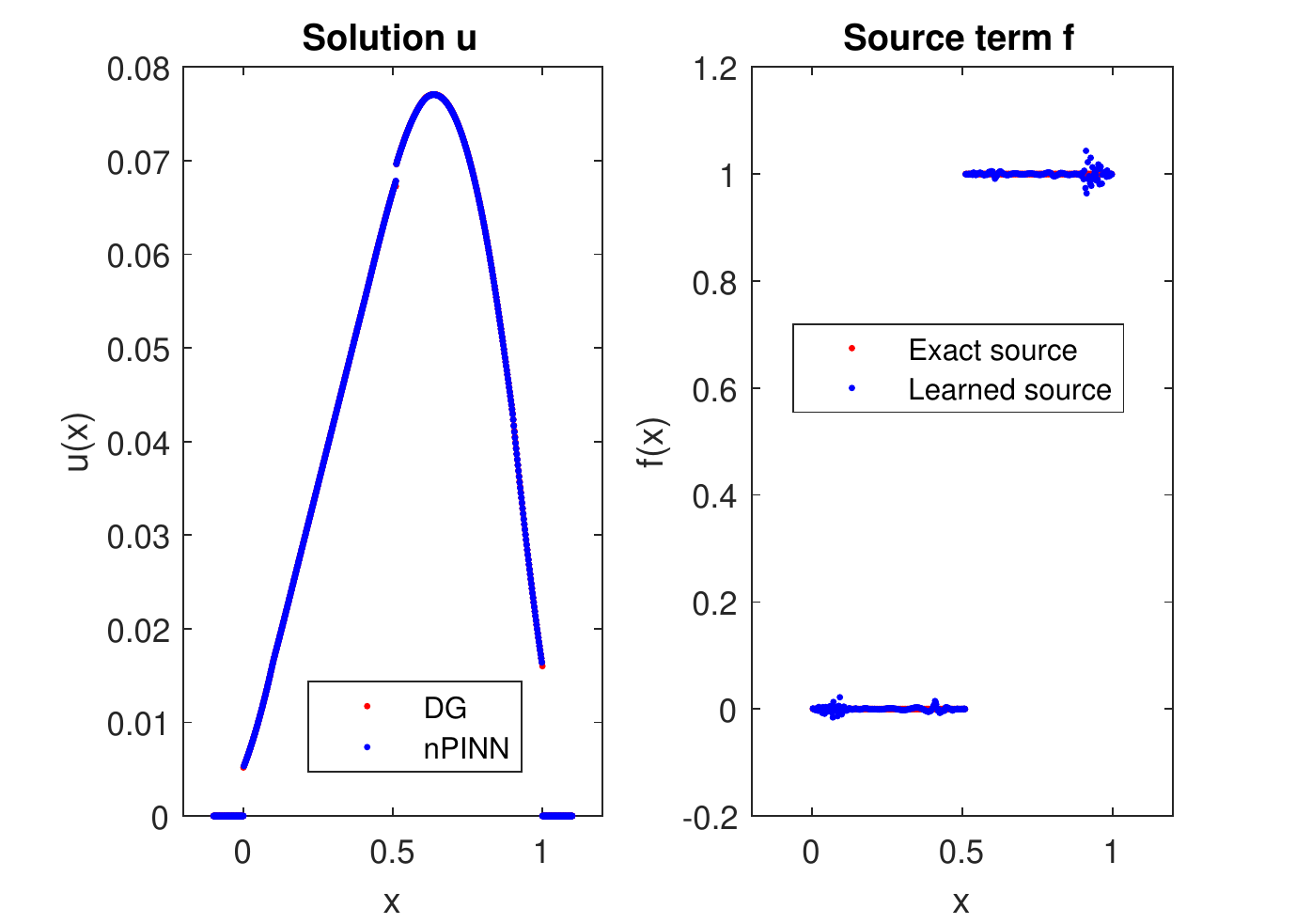}
\caption{For problem (II): nPINNs and finite element solutions on the left; trained $-\mcL^{\delta,\alpha}\unn$ and exact $f$ on the right.}
\label{1d-dis2}
\end{figure}

\paragraph{Treating solutions with discontinuities}
We consider two problem settings for $C_{\delta,\alpha}=C'_{\delta,\alpha}$ and $\omg=(0,1)$. In problem (I), taken from \cite{chen2011continuous}, we set $\alpha=0$, $\delta=0.3$ and
\begin{equation}
f(x)=\left\{\begin{array}{ll}
\;\; 0 & x\in[0,0.5-\delta)\\[2mm]
\begin{array}{l}
-\frac{2}{\delta^2}\left[\frac{1}{2}\delta^2-\delta+\frac{3}{8}+(2\delta-\frac{3}{2}-\ln{\delta})x  \right.\\
\left.+(\frac{3}{2}+x^2\ln{\delta}-(x^2-x)\ln{\frac{1}{2}-x} \right]
\end{array} & x\in[0.5-\delta, 0.5)\\[6mm]
\begin{array}{l}
-\frac{2}{\delta^2}\left[\frac{1}{2}\delta^2
-\delta-\frac{3}{8}+(2\delta+\frac{3}{2}+x\ln{\delta}) \right.\\
\left.-(\frac{3}{2}+x^2\ln{\delta}
+(x^2-x)\ln{x-\frac{1}{2}} \right]
\end{array} & x\in(0.5, 0.5+\delta)\\[6mm]    
\;\; -2 & x\in[0.5+\delta,1.0],
    \end{array} \right.  
\end{equation}
for which the corresponding analytic nonlocal solution is
\begin{equation}
u(x)=\left\{
\begin{array}{ll}
x & x\in [-\delta,0.5)\\
x^2 & x \in (0.5,1+\delta].
\end{array} \right.  
\end{equation}
In problem (II), we set $\alpha=-1$, $\delta=0.1$ and
\begin{equation}
f(x)=\left\{
\begin{array}{ll}
0 & x \in [0,0.51111]\\
1 & x \in  (0.51111, 1.0],
    \end{array} \right.  
\end{equation}
for which the corresponding reference solution is computed by using discontinuous piecewise linear finite elements \cite{chen2011continuous} on a grid with discretization size $2^{-12}$. Results in Figures \ref{1d-dis} and \ref{1d-dis2} show that nPINNs can accurately match the reference solutions. In fact, we have $\epsilon$=1.7e-03 for problem (I) using ~600 residual points and ~2000 testing points; and $\epsilon$=4.9e-04 for problem (II) using ~600 residual points and ~2000 testing points. These results demonstrate the flexibility of the algorithm and its robustness with respect to rough solutions. We take $m=40$ and $M=30$ in composite Gauss quadrature.


\begin{figure}[H]
\centering
\subfloat[2D]{
\includegraphics[width=.7\textwidth]{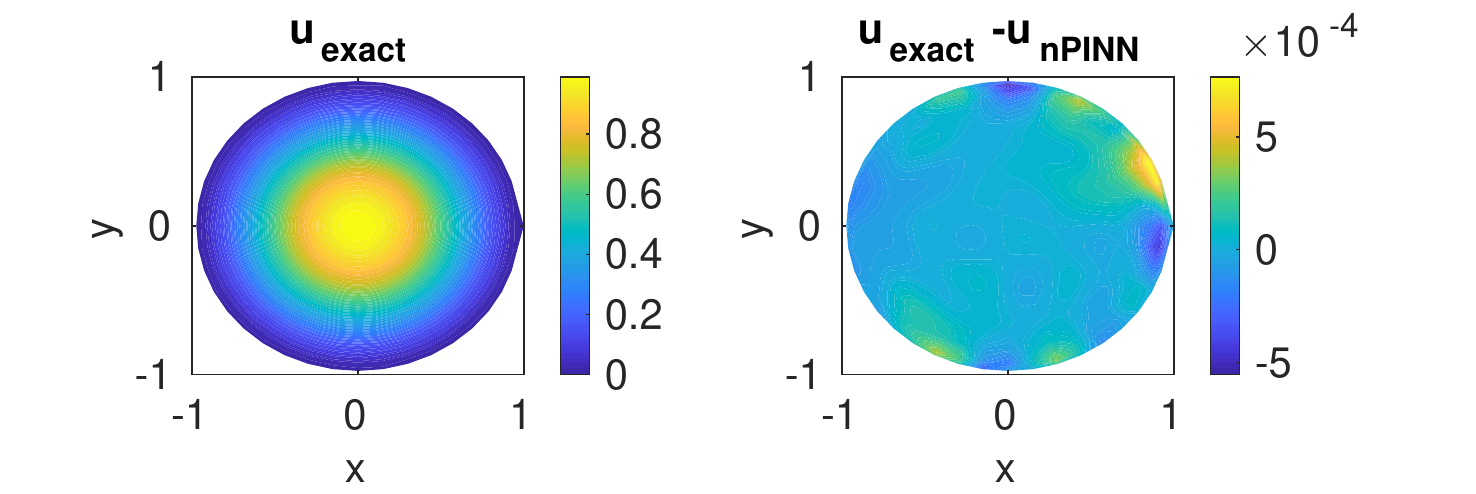}}\vfill
\subfloat[3D]{
\includegraphics[width=.7\textwidth]{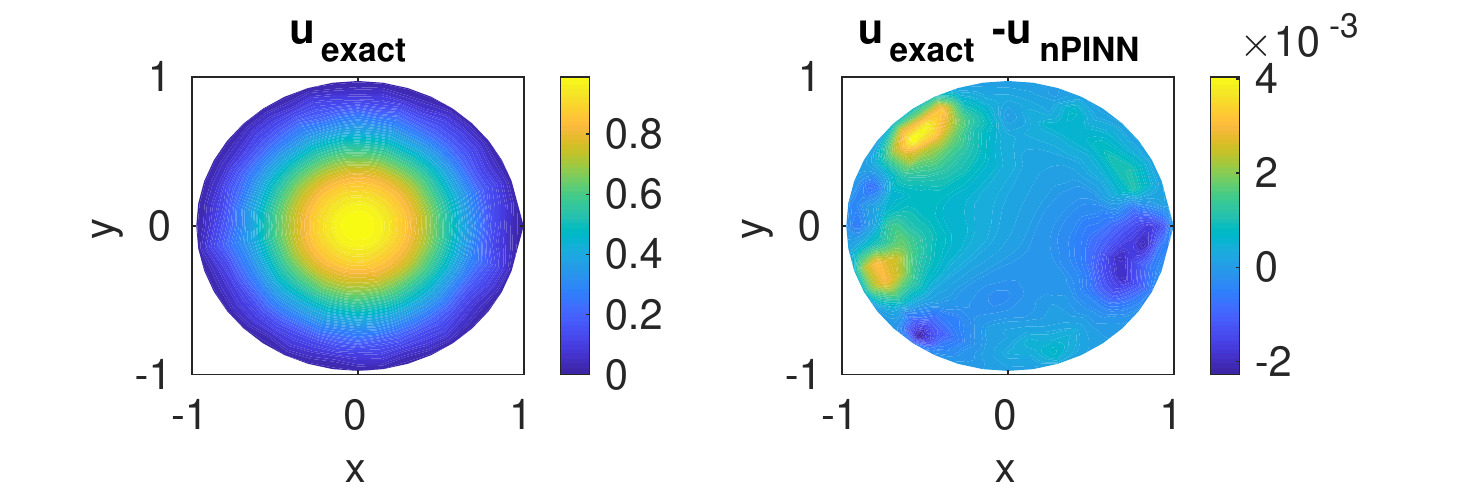}}
\caption{\label{2d-3d}Manufactured solutions on the left and their difference with nPINN solutions on the right for the two-dimensional (top) and the three-dimensional (bottom) problems.}
\end{figure}

\paragraph{Treating two- and three-dimensional problems}
We demonstrate the applicability of nPINNs to two- and three-dimensional (forward) problems with $\delta$=1.0e100; in both cases, we solve equation \eqref{eq:nonlocal-diffusion} in the unit ball. We consider the manufactured solution $u(\xb)=(1-\|\xb\|_2^2)^{1+\alpha/2}$ for which the corresponding forcing term is $f(\xb)=(-\Delta)^{\alpha/2}u(\xb)= 2^{\alpha}\Gamma(\frac{\alpha}{2}+2)\Gamma(\frac{d+\alpha}{2})\Gamma(\frac{d}{2})^{-1}(1-(1+\frac{\alpha}{d})||\xb||_2^2)$~\cite{pang2019fpinns}. Note that, because $\delta\sim\infty$, we can obtain the forcing term analytically by computing the action of the fractional Laplacian on $u$. Manufactured solutions and their difference with nPINNs solutions are reported in Figure \ref{2d-3d}; in the three-dimensional case the solution is plotted along the plane $z=0$. For $d=2$, 300 residual points, and 2000 test points the relative error is $\epsilon$=1.7e-04; and for $d=3$, 400 residual points, and 2000 testing points $\epsilon$=1.1e-03.\\

\subsection{Parameter estimation for determining the operator}\label{inverse-sec}
We use nPINNs to jointly estimate $\delta$ and $\alpha$ for a one-dimensional nonlocal Laplace operator. In $\omg=(0,1)$ we consider $N^{obs}\!=\!100$ equally spaced observation points. Recall that in the simulations presented in the previous sections observations were not needed because we were considering a forward problem with a known volume constraint and prescribed $\delta$ and $\alpha$. Here, we still assume $g$ to be known everywhere in $\omgd$ so that observation points belong to $\omg$ only. We select $N\!=\!200$ residual points scattered over $\omg$ via Sobol's sequences, and take $m\!=\!M\!=\!10$ in the composite Gauss quadrature for evaluating $\mathcal{L}^{\delta,\alpha}\unn$.

We generate four manufactured solutions by solving problem \eqref{eq:nonlocal-diffusion} for $f(x)=\sin(2\pi x)$, $g(x)=0$ and four pairs of parameters $(\delta^*,\alpha^*)$, which we refer to as {\it true parameters}. Specifically, we consider: (I) $\delta^*=1.4$ and $\alpha^*=0.8$; (II) $\delta^*=14$ and $\alpha=0.8$; (III) $\delta^*=1400$ and $\alpha^*=0.8$; (IV) $\delta^*=1.4$ and $\alpha^*=1.8$. The corresponding solutions are used to sample $u_{obs}$ at observation points. In order to investigate the sensitivity of the algorithm to the initial guess, we consider three pairs of initializations $(\delta_0, \alpha_0)$, specifically (1) $\delta_0 = 1$ and $\alpha_0 = 0.5$, (2) $\delta_0 = 10$ and $\alpha_0 = 0.5$, and (3) $\delta_0 = 1000$ and $\alpha_0=0.5$. 

For the solution of the optimization problem we run four million Adam optimization iterations and decrease the learning rate\footnote{The learning rate $\eta$ measures how fast we update the parameters $\mathbf{w}$, being optimized in (stochastic) gradient descent algorithms: $\mathbf{w}^{k+1}=\mathbf{w}^k-\eta \nabla Loss(\mathbf{w})$. We allow it to change with the iteration index $k$.} $\eta$ as the number of iterations increases:
\begin{equation*}
\eta = \left\{
\begin{array}{ll}
    10^{-3}, &  \mbox{iter} \in [0, 2\times 10^6], \\
    10^{-4}, &  \mbox{iter} \in (2 \times 10^6, 3\times 10^6],\\
    10^{-5}, &  \mbox{iter} \in (3\times 10^6, 4\times 10^6],
\end{array} \right.
\end{equation*}
where ``iter'' is the current iteration number.

\begin{figure}[H] 
\centering
\includegraphics[width=.7\textwidth]{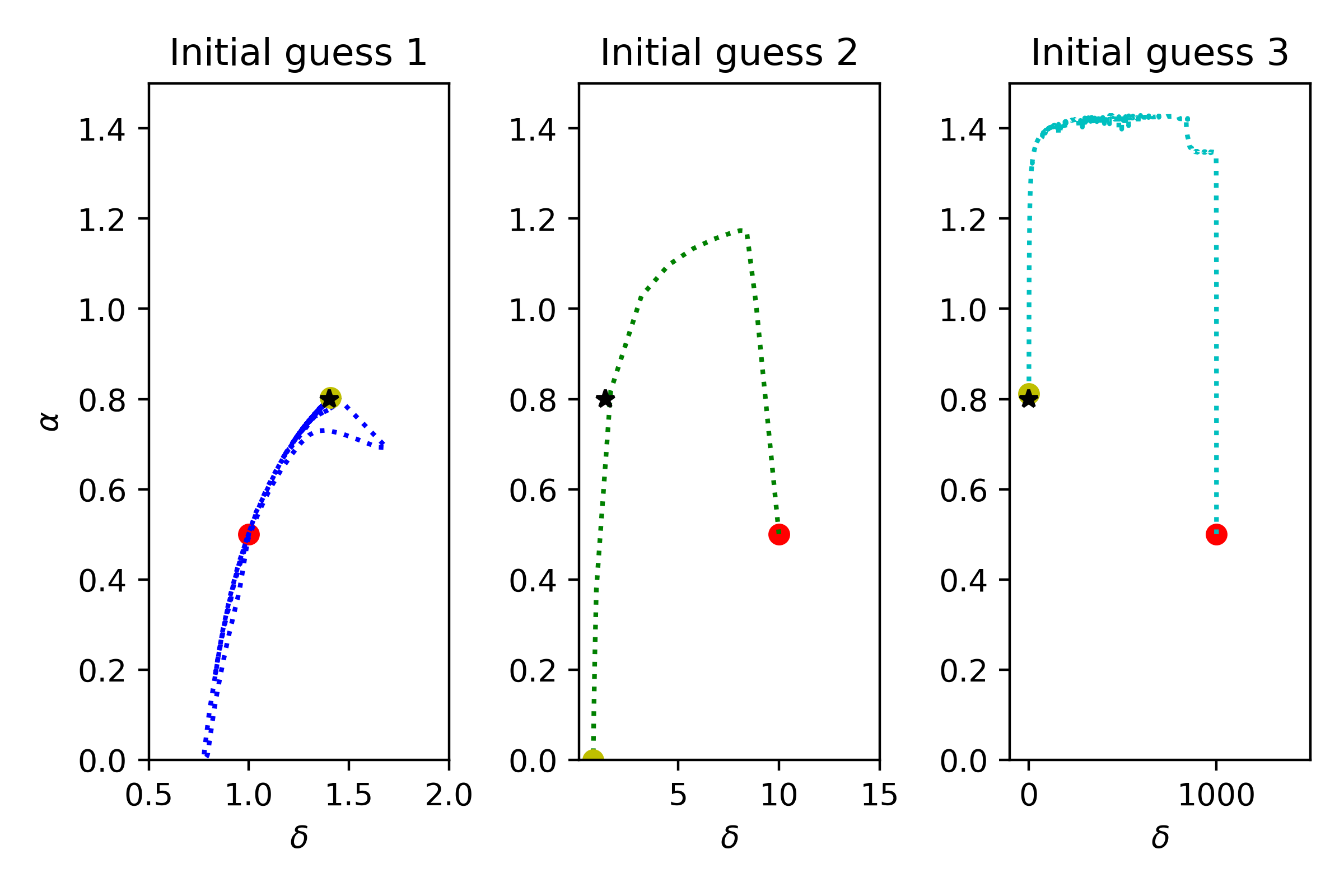}
\caption{Optimization trajectories in the parameter space $(\delta, \alpha$) with true parameters $(\delta^*,\alpha^*)=(1.4,0.8)$ and initial guesses: (1) $(\delta_0,\alpha_0)=(1,0.5)$, (2) $(\delta_0,\alpha_0)=(10,0.5)$, and (3) $(\delta_0,\alpha_0)=(1000,0.5)$. The red and yellow dots represent the initial guess and the optimal pair, the black stars represent the true parameters. Final loss values are (1) \textbf{2.68e-6}, (2) 9.89e-5, and (3) 3.56e-6. The relative errors $\epsilon$ are (1) \textbf{2.64e-4}, (2) 1.43e-2, and (3) 2.69e-4. Initial guess (1) yields the lowest final loss and relative error. Initial guess (2) leads to a bad local minimum, which will be neglected due to the large loss value compared to initial guess (1). The local minimum reached by initial guess (3) is acceptable since its final loss value is of the same order of (1).}
\label{fig1}
\end{figure}

\begin{figure}[H] 
\centering
\includegraphics[width=.7\textwidth]{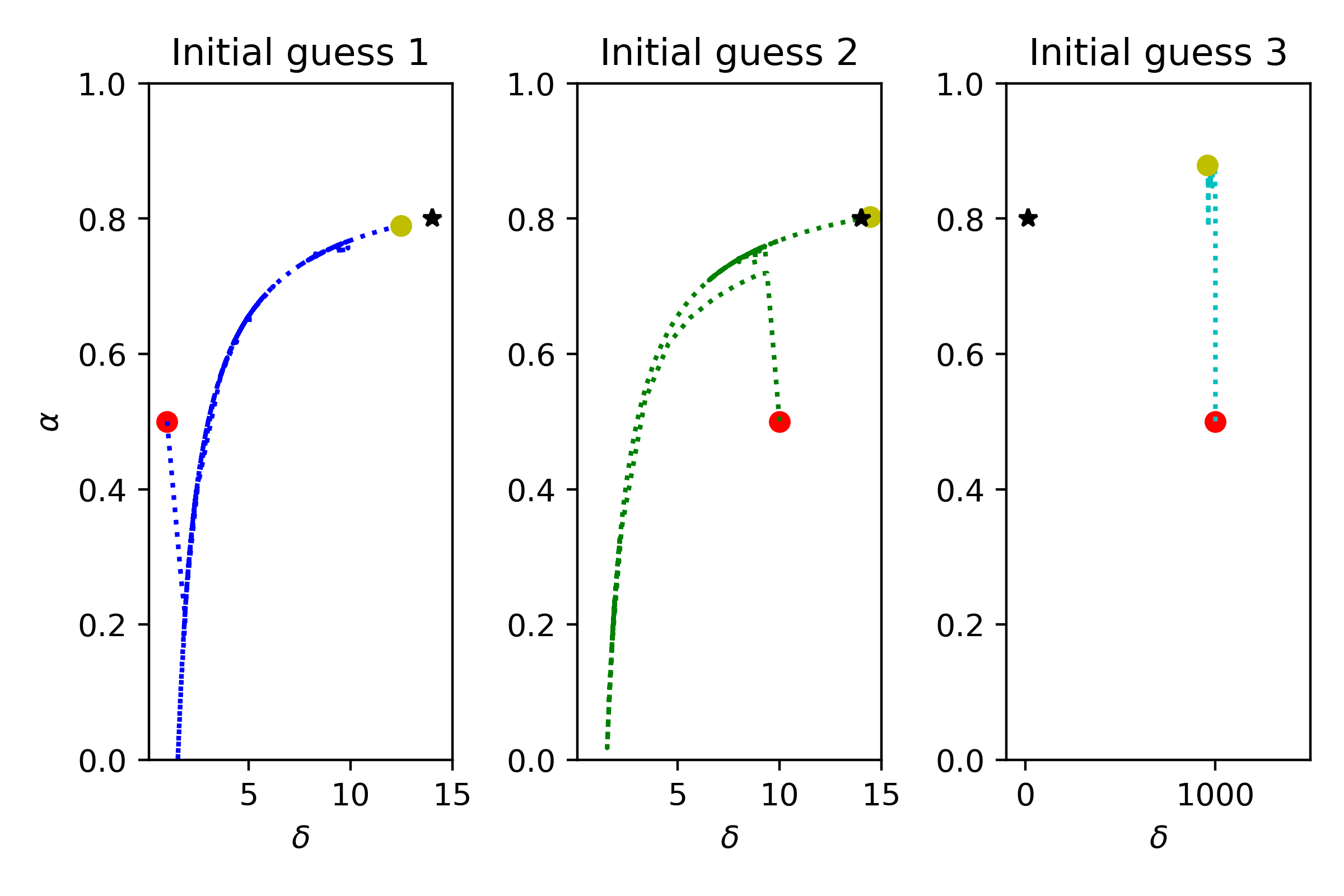}
\caption{Optimization trajectories in the parameter space $(\delta, \alpha$) with true parameters $(\delta^*,\alpha^*)=(14,0.8)$ and initial guesses: (1) $(\delta_0,\alpha_0)=(1,0.5)$, (2) $(\delta_0,\alpha_0)=(10,0.5)$, and (3) $(\delta_0,\alpha_0)=(1000,0.5)$. Final loss values are (1) 2.53e-6, (2) \textbf{2.47e-6}, and (3) 4.08e-6. The relative errors $\epsilon$ are (1) 5.46e-4, (2) \textbf{1.06e-4}, and (3) 8.17e-4. Initial guess (2) yields the lowest final loss and relative error.. Initial guess (1) is acceptable because the value of the loss function is comparable to (2). Initial guess (3) converges to the local minimum $(\delta, \alpha)=(958, 0.879)$ with a loss value of the same order of (1) and (2). Thus, the nonlocal operator $\mathcal{L}^{958,0.879}$ {\it mimics} the operator $\mathcal{L}^{14,0.8}$ for the source term $f$ and the volume constraint $g$, due to the comparable loss values. Note that $\epsilon$ for initial guess (3) is also comparable to (1) and (2).}
\label{fig2}
\end{figure}

\begin{figure}[H] 
\centering
\includegraphics[width=.7\textwidth]{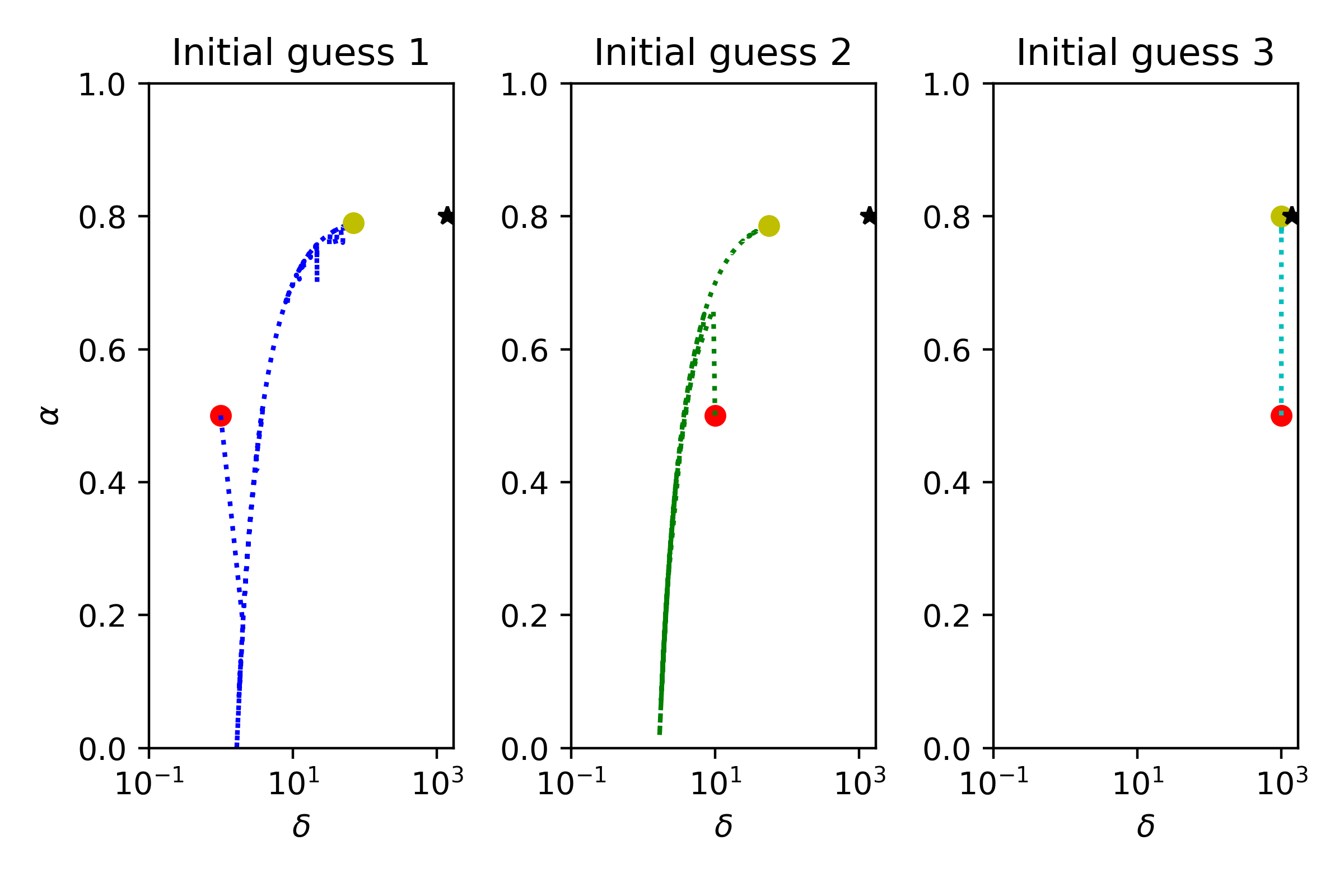}
\caption{Optimization trajectories in the parameter space $(\delta, \alpha$) with true parameters $(\delta^*,\alpha^*)=(1400,0.8)$ and initial guesses: (1) $(\delta_0,\alpha_0)=(1,0.5)$, (2) $(\delta_0,\alpha_0)=(10,0.5)$, and (3) $(\delta_0,\alpha_0)=(1000,0.5)$. Final loss values are (1) 2.72e-6, (2) \textbf{2.54e-6}, and (3) 3.56e-6. The relative errors $\epsilon$ are (1) 3.56e-4, (2) \textbf{2.40e-4}, and (3) 2.69e-4. While the estimated $\alpha$'s for (1) and (2) are $\approx\!0.8$ (true value), values of $\delta$ are much smaller than the true value. However, the final loss values are comparable with (3), which does reach the true parameters. This indicates that (1) and (2) yield two {\it mimicking operators} for $\mathcal{L}^{1400,0.8}$; it also implies that the farther the initial guess is from the true parameters, the more likely the operator mimicking phenomenon occurs.}
\label{fig3}
\end{figure}

\begin{figure}[H] 
\centering
\includegraphics[width=.7\textwidth]{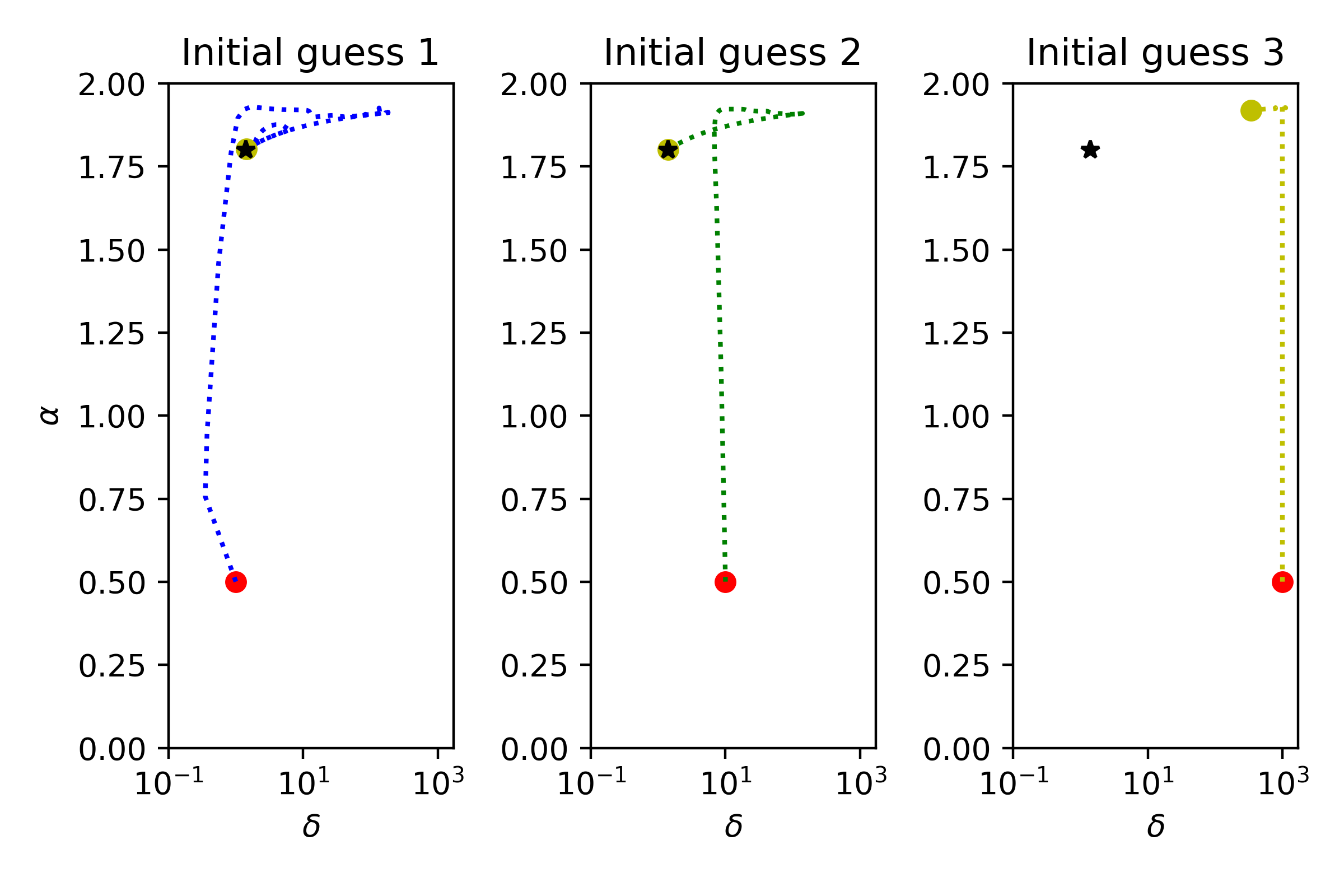}
\caption{Optimization trajectories in the parameter space $(\delta, \alpha$) with true parameters $(\delta^*,\alpha^*)=(1.4,1.8)$ and initial guesses: (1) $(\delta_0,\alpha_0)=(1,0.5)$, (2) $(\delta_0,\alpha_0)=(10,0.5)$, and (3) $(\delta_0,\alpha_0)=(1000,0.5)$. Final loss values are (1) 2.31e-7, (2) \textbf{1.63e-7}, and (3) 8.22e-5. The relative errors $\epsilon$ are (1) 1.95e-4, (2) \textbf{4.05e-5}, and (3) 8.41e-3. Initial guesses (1) and (2) reach the true parameters and have comparable final values of the loss function; initial guess (3) is disregarded because of higher values of the loss function. There is no operator mimicking for this pair of true parameters.}
\label{fig4}
\end{figure}

In Figures \ref{fig1}--\ref{fig4} we report the optimization trajectories in the parameter space $(\delta,\alpha)\in (0,+\infty)\times (0,2]$ in correspondence of the initial guesses (1)--(3) for every pair of true parameters (I)--(IV). In all figures the red disk corresponds to the initial guess, the yellow disc to the optimal (estimated) pair and the black star to the true pair. Note that the trajectories we plot are the projection of the trajectories in whole parameter space $(\delta,\alpha,\boldsymbol{\theta})$ (whose dimension is of the order of 1e+03 or 1e+04) onto a two-dimensional subspace $(\delta,\alpha)$. 

From the figures, we see that both the locations of true parameters $(\delta^*,\alpha^*)$ and initial guess $(\delta_0,\alpha_0)$ affect the outcome of the optimization. First, the shape of the loss function \eqref{loss_fun_discrete} varies for different true parameters. Second, due to the non-convexity of loss function, different initial guesses may lead to either the global minimum or other (possibly good) local minima corresponding to comparable loss values. The former case indicates that using nPINNs allows us to identify true parameters, the latter case indicates the occurrence of the operator mimicking phenomenon, see Figures \ref{fig2} and \ref{fig3}. Third, due to the non-convexity, optimal parameters may correspond to a local minimum for which the value of the loss function is much higher than the optimal one (global minimum). This can be observed in the second subplot of Figure \ref{fig1} and in third subplot of Figure \ref{fig4}.
Next, we illustrate the sensitivity of the algorithm with respect to the number of observation points $N^{obs}$. Even though, intuitively, more observation points may lead to more accurate parameters, this is not necessarily the case, as results in Table~\ref{table} confirm due to optimization errors. Here, for different true parameters and initial guesses, we report the estimated parameters for increasing $N^{obs}$. From the first two rows, we see that predicted values are closer to the true ones as we increase $N^{obs}$; however, the last row shows that only $10$ observation points are enough to obtain accurate predictions. This indicates that $N^{obs}$, $(\delta^*,\alpha^*)$ and $(\delta_0,\alpha_0)$ jointly affect the convergence of nPINNs.
\begin{table}[H]
\centering
\begin{center}
\begin{tabular}{l|c|rllll}
$\delta^*$, $\alpha^*$ & $\delta_0$, $\alpha_0$ & $N^{obs}$=& 10 & 20 & 50 & 100\\
\hline
1.4 & 1.0 & $\delta=$ & 0.774 & 0.783 & 1.404 & 1.406 \\
0.8 & 0.5 & $\alpha=$ & 0.004 & 0.041 & 0.802 & 0.803 \\\hline
14  & 10  & $\delta=$ & 21.01 & 15.07 &  13.22 & 14.44 \\
0.8 & 0.5 & $\alpha=$ & 0.827 & 0.806 & 0.795 & 0.802 \\ \hline
1.4 & 10  & $\delta=$ & 1.417 & 1.396 & 1.392 & 1.404 \\
1.8 & 0.5 & $\alpha=$ & 1.800 & 1.799 & 1.799 & 1.800 \\
\end{tabular}
\end{center}
\caption{Sensitivity of nPINNs to the number of observation points $N^{obs}$. The minimum number of observation points that we need to recover the true parameters depends is sensitive to the location of $(\delta^*,\alpha^*)$.}
\label{table}
\end{table}

\section{Application to turbulence modeling of Couette flow} \label{sec:turbulence}

\begin{figure}[H] 
\centering
\includegraphics[width=.7\textwidth]{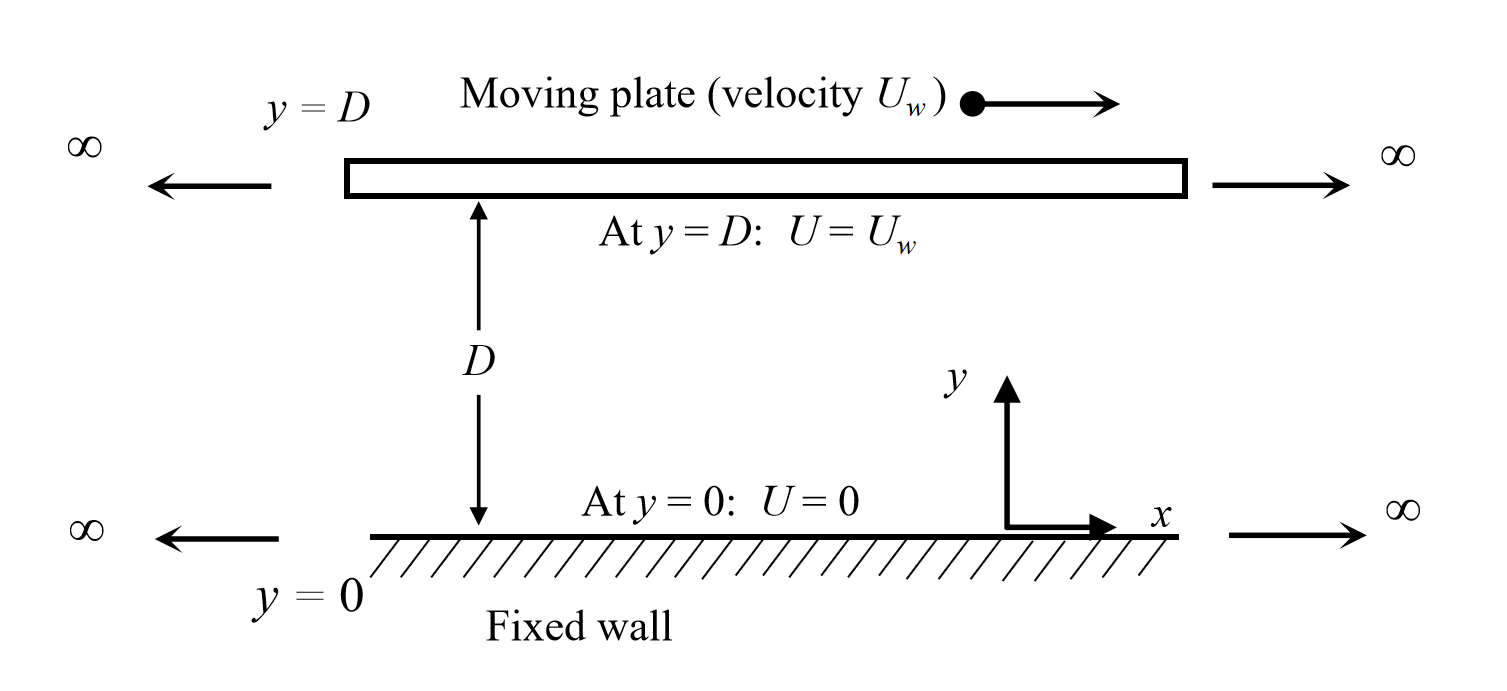}
\caption{Problem configuration for the Couette flow: the viscous fluid flows in the space between two infinitely large surfaces, one of which is moving tangentially relative to the other with velocity $U_w$.}
\label{Couette_flow}
\end{figure}

In this section we extend our algorithm to a more complex application of modeling wall turbulence using
the Reynolds-Averaged Navier-Stokes (RANS) equations. We propose a nonlocal model for Couette flow, i.e., the mean flow of a viscous fluid in the space between two infinitely large surfaces, one of which is moving tangentially relative to the other, as shown in Figure \ref{Couette_flow}. We introduce a new unified nonlocal operator that, compared to \eqref{eq:param-L}, features a different normalizing constant and a spatially variable decay rate $\alpha$, which allows us to describe more effectively the turbulent flow
behavior as a function of distance from the wall.
The classical local model is obtained by simplifying the RANS equations for the Couette flow set up; this yields the following one-dimensional equation \cite{kundu2012fluid,mehta2019discovering}
\begin{equation}\label{Couette-model}
\frac{d}{dy^+}\left(\frac{dU^+}{dy^+}-(\overline{uv})^+\right) = 0, \quad y^+ \in [0, 2\mbox{Re}_{\tau}].
\end{equation}
Here, $U^+$ and $y^+$ are dimensionless variables based on wall units and are defined as
\begin{equation}
    U^+ = \frac{U}{U_{\tau}}, \qquad y^+ = \frac{y U_{\tau}}{\nu},
\end{equation}
where $U$ is the mean value of the instantaneous velocity in the streamwise direction $x$, $U^2_\tau = \frac{\tau_w}{\rho}$ is the friction velocity, with $\rho$ the flow density and $\tau_w$ the wall shear stress, $y$ is the crossflow direction, and $\nu$ is the kinematic viscosity. The subscript $\tau$ refers to the shear stress. The dimensionless Reynolds stress $(\overline{uv})^+$ is defined as 
\begin{equation}
(\overline{uv})^+ = \frac{\overline{uv}}{U^2_{\tau}},
\end{equation}
where the overbar represents the average, and $u$ and $v$ are deviations from the mean values of the velocities in the $x$ and $y$ direction. Note that since the flow only occurs in the streamwise direction $x$, the mean value of the instantaneous velocity in the crossflow direction, $V$, is zero, i.e. $V\equiv 0$~\cite{kundu2012fluid}. The Reynolds number is defined as $\mbox{Re}_{\tau}=\frac{U_{\tau}\cdot (D/2)}{\nu}$. 
The term $\frac{dU^+}{dy^+}-(\overline{uv})^+$ in \eqref{Couette-model} represents the dimensionless total shear stress $\tau^+$, also defined as
\begin{equation}
\tau^+=\frac{\tau}{\rho U_{\tau}^2}.
\end{equation}

Equation \eqref{Couette-model} indicates that the total shear stress $\tau^+$ is constant with respect to $y^+$, i.e. $\tau^+=c$. This constant can be determined by evaluating the total shear stress near the wall, i.e., at $y^+\approx 0$. It is well known~\cite{kundu2012fluid} that in the viscous sublayer ($y^+<5$) the mean flow velocity $U^+\approx y^+$, the viscous stress $\frac{dU^+}{dy^+}$ dominates, and the Reynolds stress is negligible, i.e., $(\overline{uv})^+\approx 0$. Thus, we have $\tau^+(y^+=0)=\frac{dU^+(0)}{dy^+}-(\overline{uv})^+=1$. The total shear stress equation can be rewritten as
\begin{equation}\label{eq:classical-shearstress-eq}
    \frac{dU^+}{dy^+}-(\overline{uv})^+ = 1,\quad y^+ \in [0,2Re_{\tau}].
\end{equation}

We propose a different nonlocal model for the description of the Couette flow inspired by our unified nonlocal Laplace operator. Specifically, we model the total shear stress with a new nonlocal operator, $\widetilde\mcL^{\delta,\alpha(y^+)}$, similar to the one introduced in \eqref{eq:param-L}. We substitute \eqref{eq:classical-shearstress-eq} with
\begin{equation}\label{eq:nonlocal-shearstress-eq}
\widetilde{\mathcal{L}}^{\delta,\alpha(y^+)}U^+ = 1, \quad 
\delta>0, \; \alpha(y^+) \in (0,1), \; y^+\in[a,b],
\end{equation}
where
\begin{multline}\label{nonlocal-Couette}
    \widetilde\mcL^{\delta,\alpha(y^+)}U^+(y^+) := -C(y^+)\left[\int_{-\delta}^0\frac{U^+(y^+ + z)-U^+(y^+)}{|z|^{1+\alpha(y^+)}}dz 
    - \int_0^{\delta}\frac{U^+(y^+ + z)-U^+(y^+)}{|z|^{1+\alpha(y^+)}}dz \right.\\
    \left. + (-\delta_{max})^{-\alpha(y^+)}U^+(y^+ +\delta_{max}) - \delta_{min}^{-\alpha(y^+)}U^+(y^+ +\delta_{min}) \right] ,
\end{multline}
with $C(y^+)=\frac{\alpha(y^+)}{2\Gamma(1-\alpha(y^+))}$, $\delta_{max}=\max\{-\delta,a-y^+\}$, $\delta_{min}=\min\{\delta,b-y^+\}$, $a=0$, and $b=2Re_{\tau}$.

\paragraph{Properties of the variable-decay operator}
Similar to $\mcL^{\delta,\alpha}$, the discrete operator in \eqref{nonlocal-Couette} also spans a broad range of variations and includes well-known classical and fractional operators as limits. We first analyze the limit as $\delta \rightarrow \infty$; we have
\begin{equation}\label{Couette-limit}
\lim \limits_{\delta\rightarrow\infty}
\widetilde\mcL^{\delta,\alpha(y^+)}U^+(y^+) =
\frac{1}{2}\left(^*D^{\alpha(y^+)}_{a+}U^+(y^+)-
^*D^{\alpha(y^+)}_{b-}U^+(y^+)\right),
\end{equation}
where
\begin{equation}
 ^*D^{\alpha(y^+)}_{a+}U^+(y^+)=\frac{1}{\Gamma(1-\alpha(y^+))}\int_a^{y^+}(y^+-z)^{-\alpha(y^+)}(U^{+})'(z)dz,
\end{equation}
and
\begin{equation}
 ^*D^{\alpha(y^+)}_{b-}U^+(y^+)=-\frac{1}{\Gamma(1-\alpha(y^+))}\int_{y^+}^b(z-y^+)^{-\alpha(y^+)}(U^{+})'(z)dz,
\end{equation}
are the left- and right-sided variable-order Caputo fractional derivatives\footnote{Though straightforward, the proof of this result is long and not relevant for this work. Thus, we do not report it.}. Here, $(U^+)'(z)$ is the first derivative of $U^+(z)$ with respect to $z$. Thus, at the limit of infinite interactions $\widetilde\mcL^{\delta,\alpha(y^+)}$ reduces to a combination of Caputo fractional derivatives.

We also show that classical local derivatives can be seen as a specific instance of $\widetilde\mcL^{\delta,\alpha(y^+)}$ when simultaneously $\delta\to\infty$ and $\alpha(y^+)\rightarrow 1$. It is well-known~\cite{podlubny1998fractional} that the Caputo fractional derivative reduces to the integer-order derivative as $\alpha(y^+)\to 1$, i.e.
$\lim_{\alpha(y^+)\rightarrow 1} \, ^*D_{a+}^{\alpha(y^+)}U^+=\frac{dU^+}{dy^+}$ and $\lim_{\alpha(y^+)\rightarrow 1}\, ^*D_{b-}^{\alpha(y^+)}U^+=-\frac{dU^+}{dy^+}$.
Thus, by taking the limits as $\delta\to\infty$ and $\alpha(y^+)\to 1$ sequentially, we can show that
\begin{equation}\label{lim_alpha}
    \lim\limits_{\alpha(y^+)\rightarrow 1,\delta\rightarrow+\infty}\widetilde{L}^{\delta,\alpha(y^+)}U^+ = \frac{dU^+}{dy^+}.
\end{equation}
In this limit, our new model (\ref{eq:nonlocal-shearstress-eq}) reduces to the local model (\ref{eq:classical-shearstress-eq}) only in the viscous sublayer, inside which the Reynolds stress is almost zero. 

\subsection{Simulation results}
\begin{figure}[H] 
\centering
\includegraphics[width=.99\textwidth]{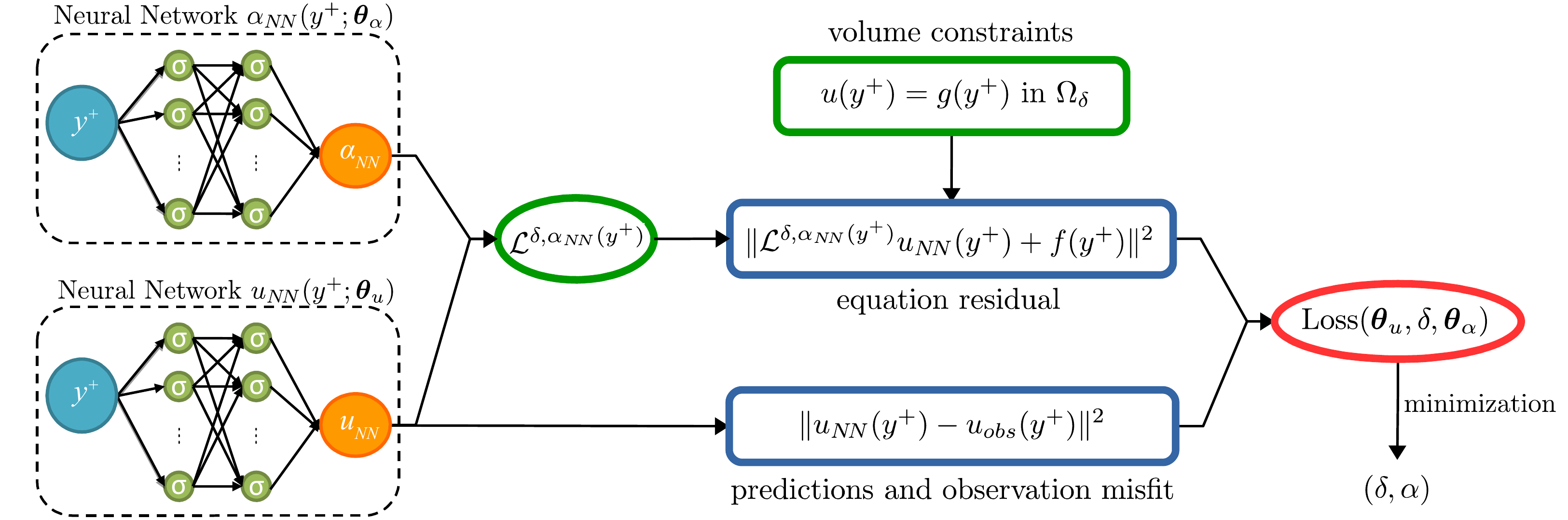}
\caption{Schematic of nPINNs for inferring the parameters $(\delta,\alpha(y^+))$ for the proposed nonlocal turbulence model.}
\label{nPINN_illustration2}
\end{figure}

To jointly estimate $\delta$ and $\alpha(y^+)$ using nPINNs, we describe the fractional order $\alpha(y^+)$ as a fully-connected NN, $\alpha(y^+)\cong\alpha_{N\!N}(y^+,\boldsymbol{\theta}_{\alpha})$. We distinguish the NNs corresponding to the solution ($U^+(y^+)\cong\unn(y^+,\thetab_u)$) and the fractional order by using different subscripts for their parameters. Compared to the nonlocal Poisson case, the parameter to be optimized, $\alpha$, is replaced by a set of parameters $\boldsymbol{\theta}_{\alpha}$, see Figure~\ref{nPINN_illustration2}, on the left. To guarantee that $\alpha_{N\!N}(y^+;\thetab_\alpha)\in(0,1)$, we consider a sigmoid function as activation function for the output layer of $\alpha_{N\!N}$. Figure \ref{nPINN_illustration2} illustrates the extension of the nPINNs algorithm to the operator described in the previous section.
We employ direct numerical simulation (DNS) data from~\cite{avsarkisov2014turbulent} as observations for the velocity $U^+$, we use the whole dataset. We consider $N$=600 residual points, uniformly distributed in log space of $y^+$, and three different Reynolds numbers $Re_{\tau}$=125, 180, and 250. Since the velocity profile features a large gradient near the wall and the moving plate, to guarantee accurate predictions, we consider deeper NNs compared to the one used for the nonlocal Poisson problem. Specifically, we set depth and width to 7 and 10, respectively, for both $\alpha_{N\!N}$ and $\unn$. The constant learning rate is $\eta$=1e-04 and the number of Adam iterations is one million.

\begin{figure}[H]
\centering
\subfloat[Profile of fractional order in wall units ($y^+$)]{
\includegraphics[width=.45\textwidth]{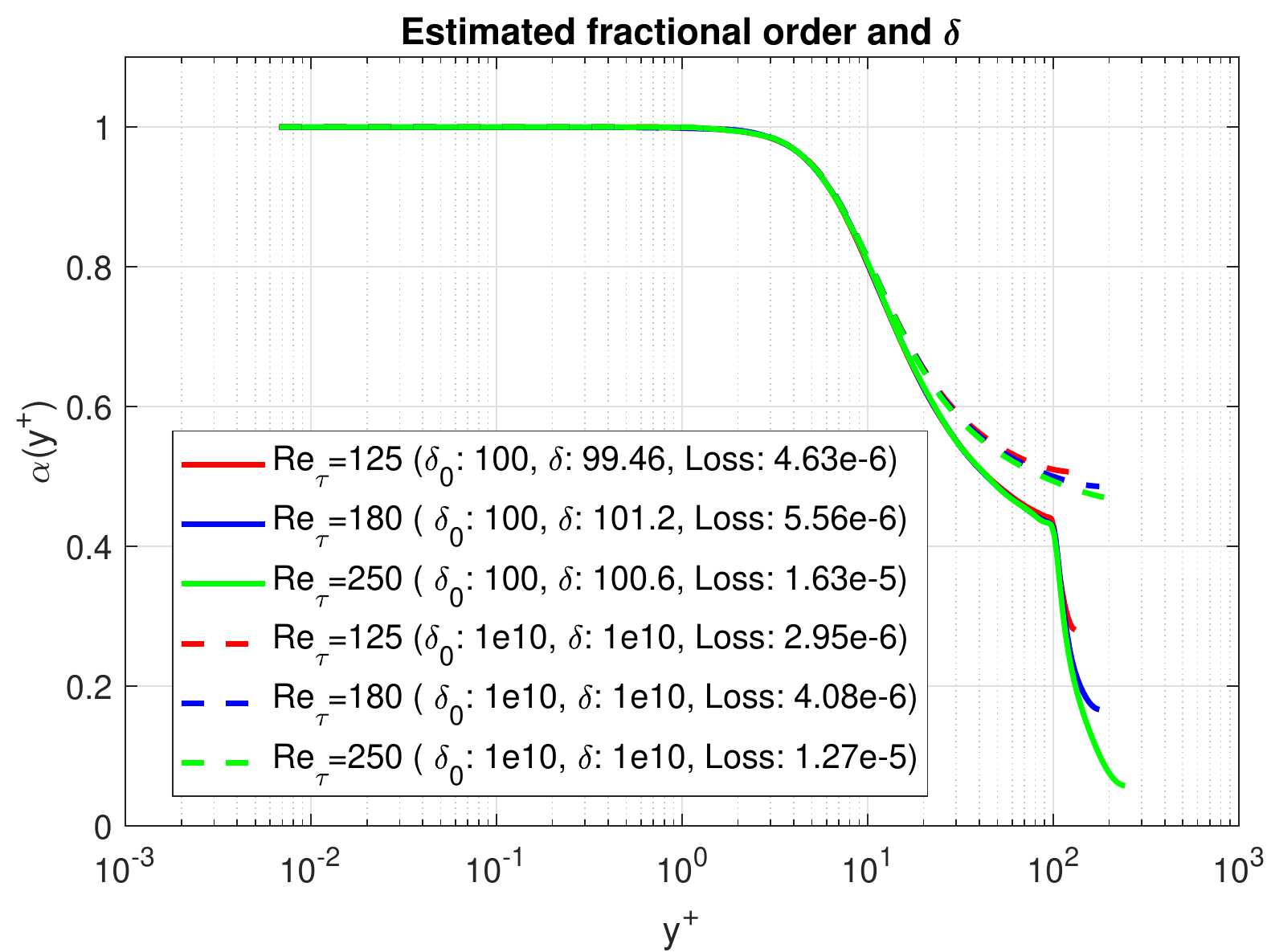}}
\subfloat[Profile of fractional order in scale $y^*$]{
\includegraphics[width=.45\textwidth]{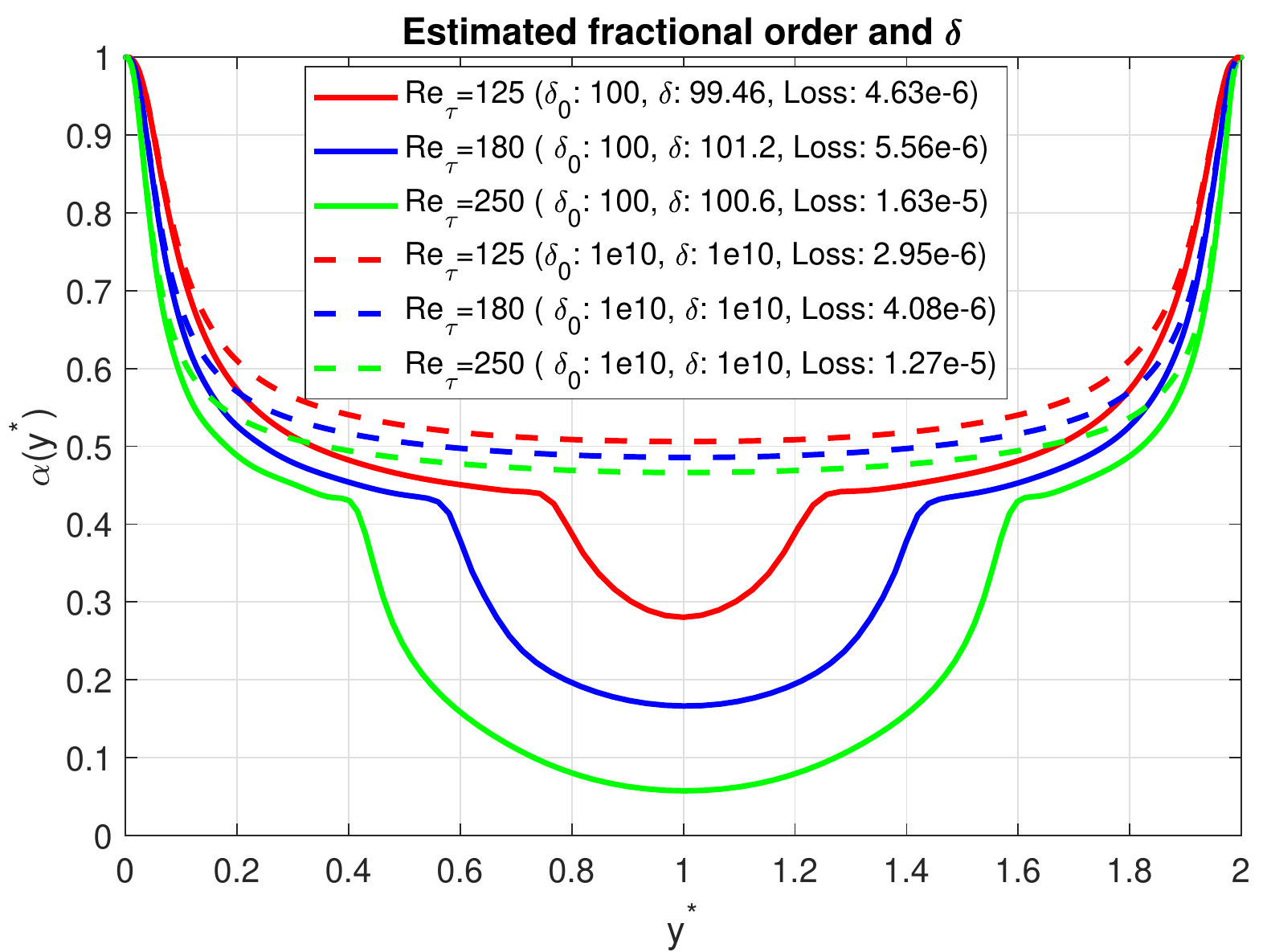}}\vfill
\subfloat[Computed total shear stress (true stress $\tau^+=1$)]{
\includegraphics[width=.45\textwidth]{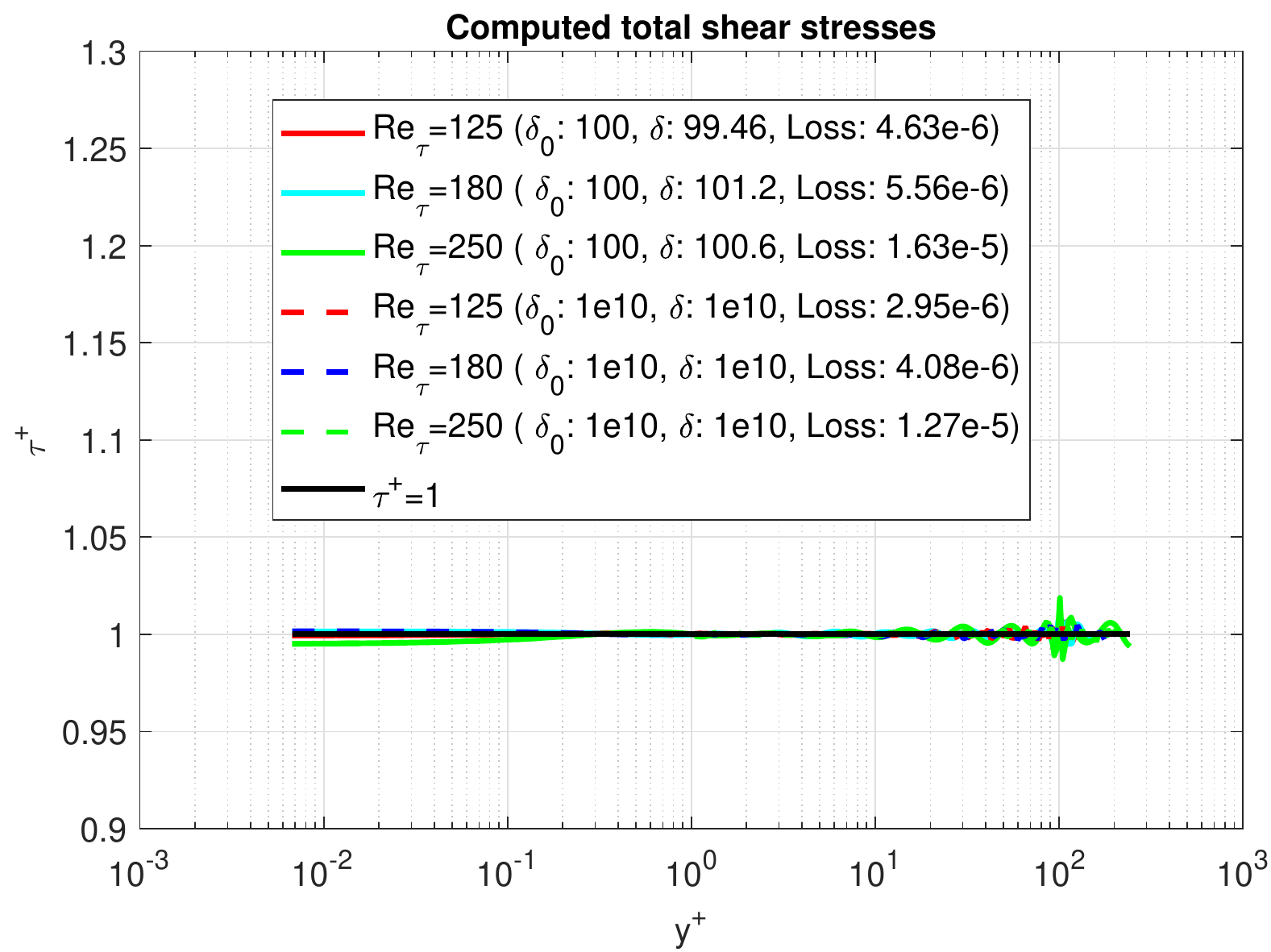}}
\caption{For the modified nonlocal operator for turbulence modeling of Couette flow: (a) $\alpha(y^+)$ in wall units learned by nPINNs; (b) $\alpha(y^*)$ in characteristic length scale; (c) total shear stresses $\tau^+$ learned by nPINNs, compared to the true total shear stresses $\tau^+=1$. In the legend, $\delta_0$ and $\delta$ are initial guess and estimate for the horizon in the nonlocal operator, respectively. ``Loss'' is the final loss after one million Adam optimization iterations. We observe the universal fractional order $\alpha(y^+)$ and universal horizon $\delta$ for varying Reynolds numbers when the initial guess for $\delta$ is fixed. For different initial guesses, we obtain two mimicking operators: $\tilde{\mathcal{L}}^{100,\alpha_{100}(y^+)}$ and $\tilde{\mathcal{L}}^{1{\rm e}10,\alpha_{\rm 1e10}(y^+)}$, where $\alpha_1(\cdot)$ and $\alpha_2(\cdot)$ denote two different fractional order profiles in the subplot (a).}
\label{Couette}
\end{figure}

Results in Figure \ref{Couette} correspond to two initial guesses for $\delta$, $\delta_0$=100 and $\delta_0$=1e10. In Figure \ref{Couette}(a), we report the fractional order profiles for different Reynolds numbers and initial guesses and the corresponding estimated $\delta$. We observe that the latter coincides with its initial guess, indicating that the value of the loss function is not sensitive to changes in the interaction radius. We also note that, independently of $Re_\tau$ and $\delta_0$, the estimated fractional order profiles are nearly on top of each other in the wall unit $y^+$. This implies that for fixed horizon $\delta$ there is a universal fractional order $\alpha(y^+)$ that reproduces the DNS data for different Reynolds numbers. From Figure \ref{Couette}(b) where $y$ is non-dimensionalized by using characteristic length and velocity scales, we observe that near the fixed wall and the moving plate the estimated fractional order is almost one. This agrees with the limit behavior \eqref{lim_alpha}, i.e. the turbulence effect is negligible in the viscous sublayer. As we move towards the centerline, the turbulence effect is intensified (the Reynolds stress dominates), and the decreasing fractional order indicates increasing nonlocal effects.
In Figure \ref{Couette}(c) we compare the computed total shear stress (i.e., the left-hand-side of the nonlocal model \eqref{eq:nonlocal-shearstress-eq}) with the true stress $\tau^+=1$; we observe that nPINNs accurately recovers the expected value. As in the previous section, we observe the operator mimicking phenomenon for two different initial guesses. For instance, as shown in Figure \ref{Couette}(a), $\tilde{\mathcal{L}}^{99.46,\alpha_{99.46}(y^+)}$ and $\tilde{\mathcal{L}}^{1{\rm e}10,\alpha_{\rm 1e10}(y^+)}$ correspond to comparable loss values (4.63e-6 and 2.95e-6, respectively) and reproduce the true stress in an equally accurate manner; however, the fractional order profiles $\alpha_1(\cdot)$ and $\alpha_2(\cdot)$ are noticeably different for $y^+>20$. Thus, the two operators are distinct, but their action on the velocity in the reduced RANS equation is the same. 

\begin{figure}[H]
\centering
\subfloat[Initial guess for $\delta$ is 100 in wall units]{
\includegraphics[width=.45\textwidth]{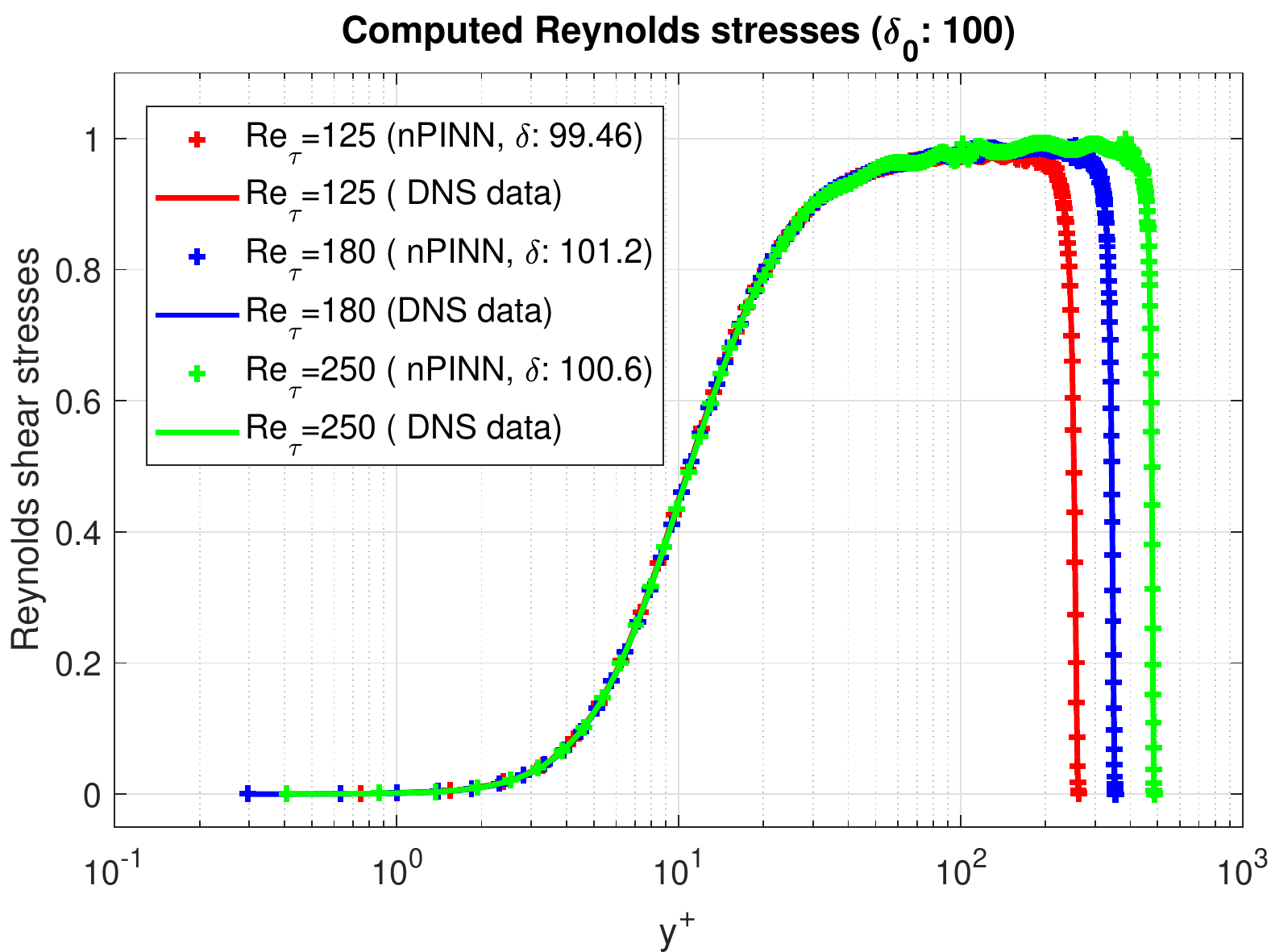}}
\subfloat[Initial guess for $\delta$ is $1\times 10^{10}$ in wall units]{
\includegraphics[width=.45\textwidth]{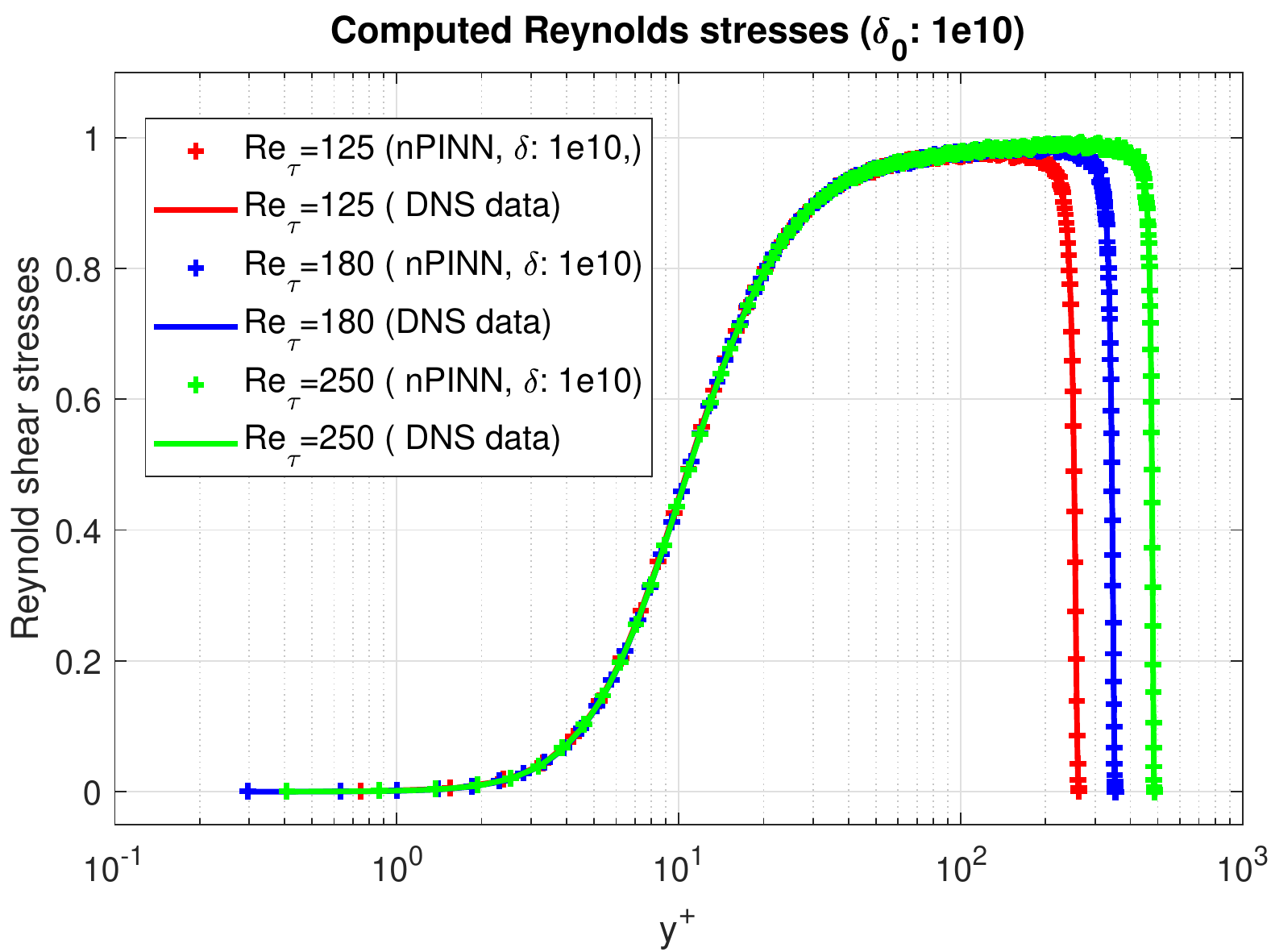}}\vfill
\caption{For the modified nonlocal operator for turbulence modeling of Couette flow: comparison of the Reynolds shear stresses computed by nPINNs and the stresses in DNS data for $\delta_0=100$ (a) and $\delta_0=1$e10 (b).}
\label{Reynolds-shear-stress}
\end{figure}

In Figure \ref{Reynolds-shear-stress} we additionally plot the computed Reynolds stresses and compare them to those reported in the DNS dataset for both initial guesses and different Reynolds numbers; the stress profiles are on top of each other.
Next, we further investigate the sensitivity of nPINNs to $\delta_0$, specifically, we consider the initial guesses $\delta_0\in \{1,\,10,\,100,\,1000,\,1{\rm e}10\}$ and analyze how they affect the learned parameters $\delta$ and $\alpha(y^+)$. 

\begin{figure}[H]
\centering
\subfloat[]{
\includegraphics[width=.45\textwidth]{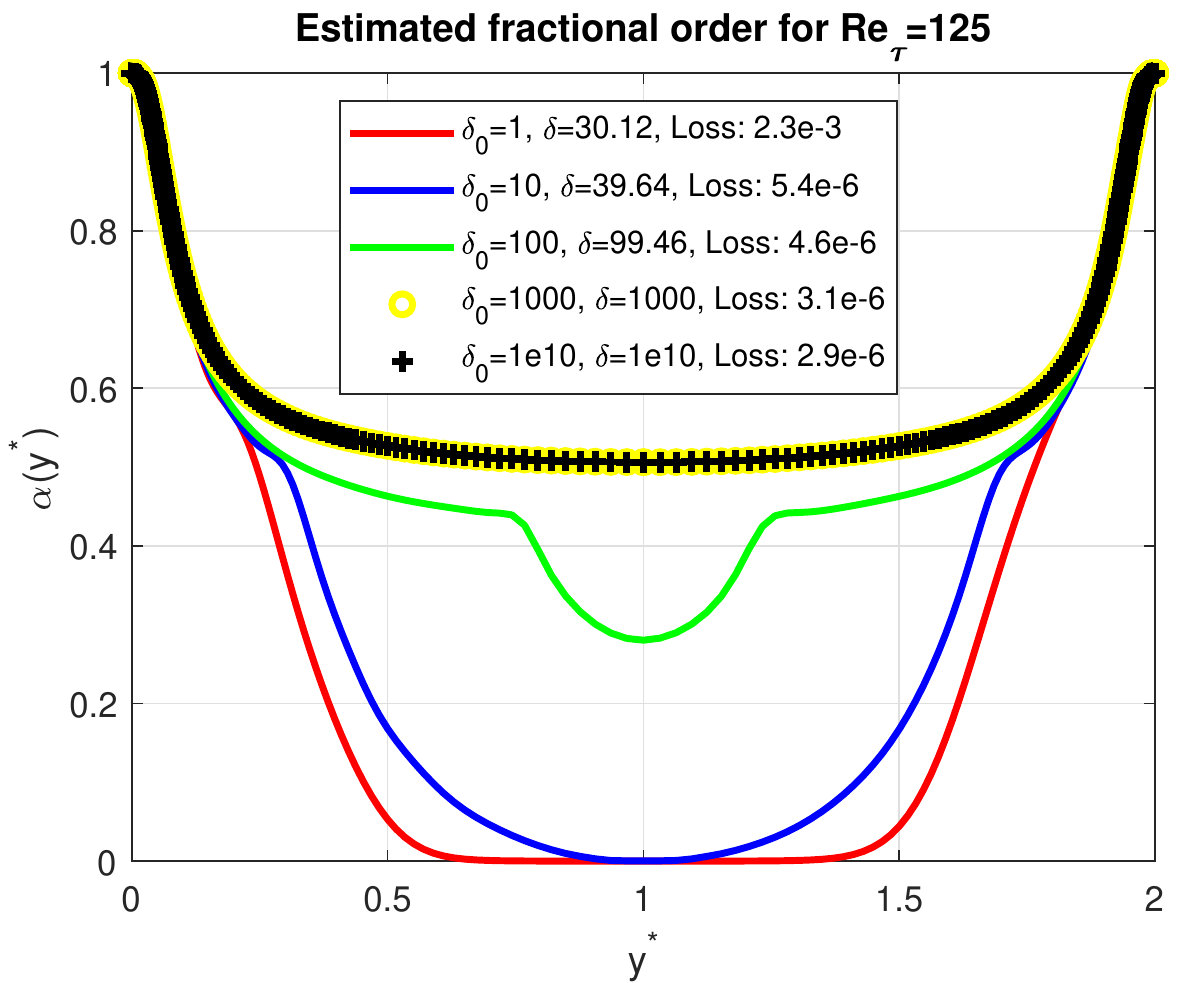}}
\subfloat[]{
\includegraphics[width=.45\textwidth]{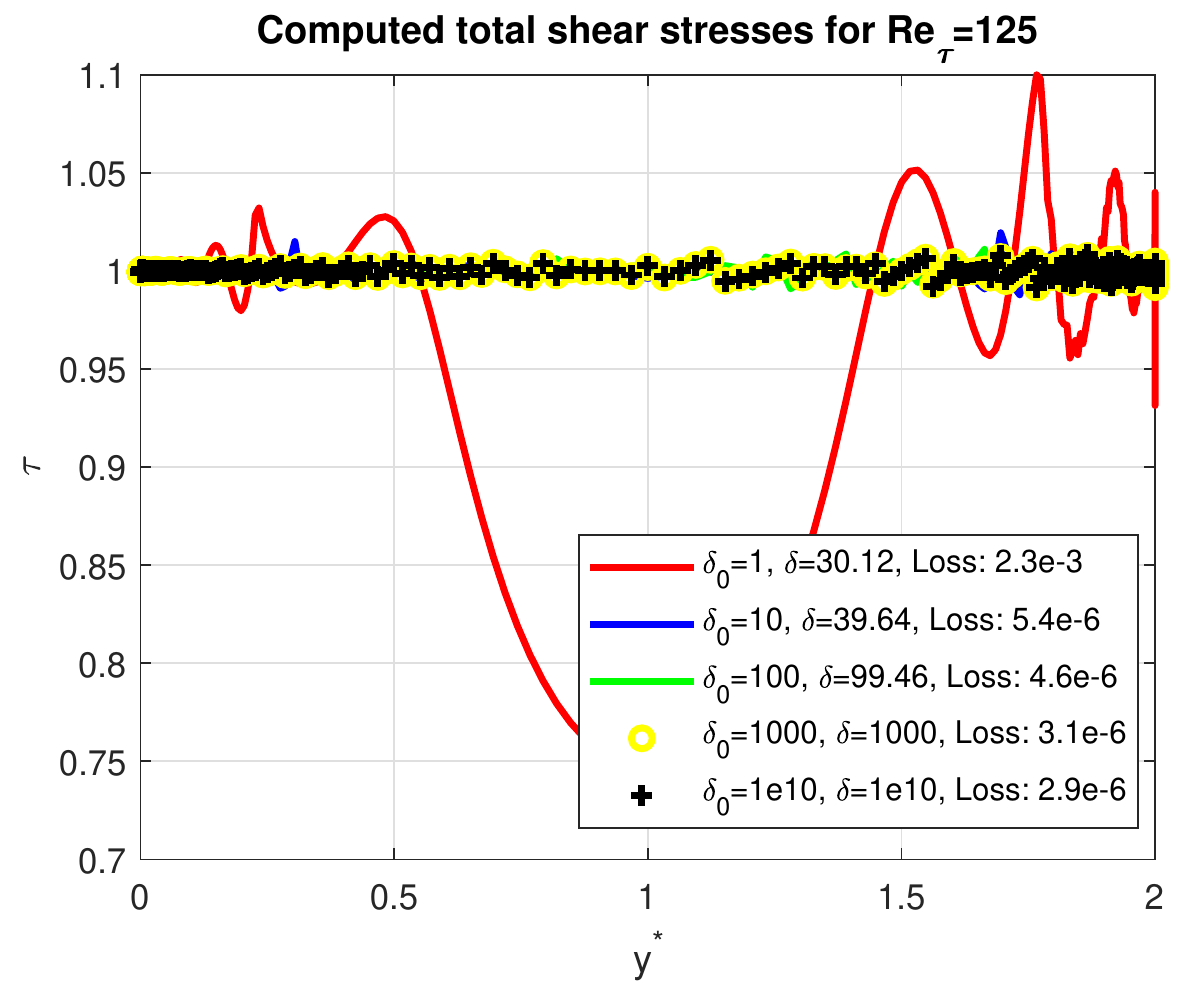}}\vfill
\subfloat[]{
\includegraphics[width=.45\textwidth]{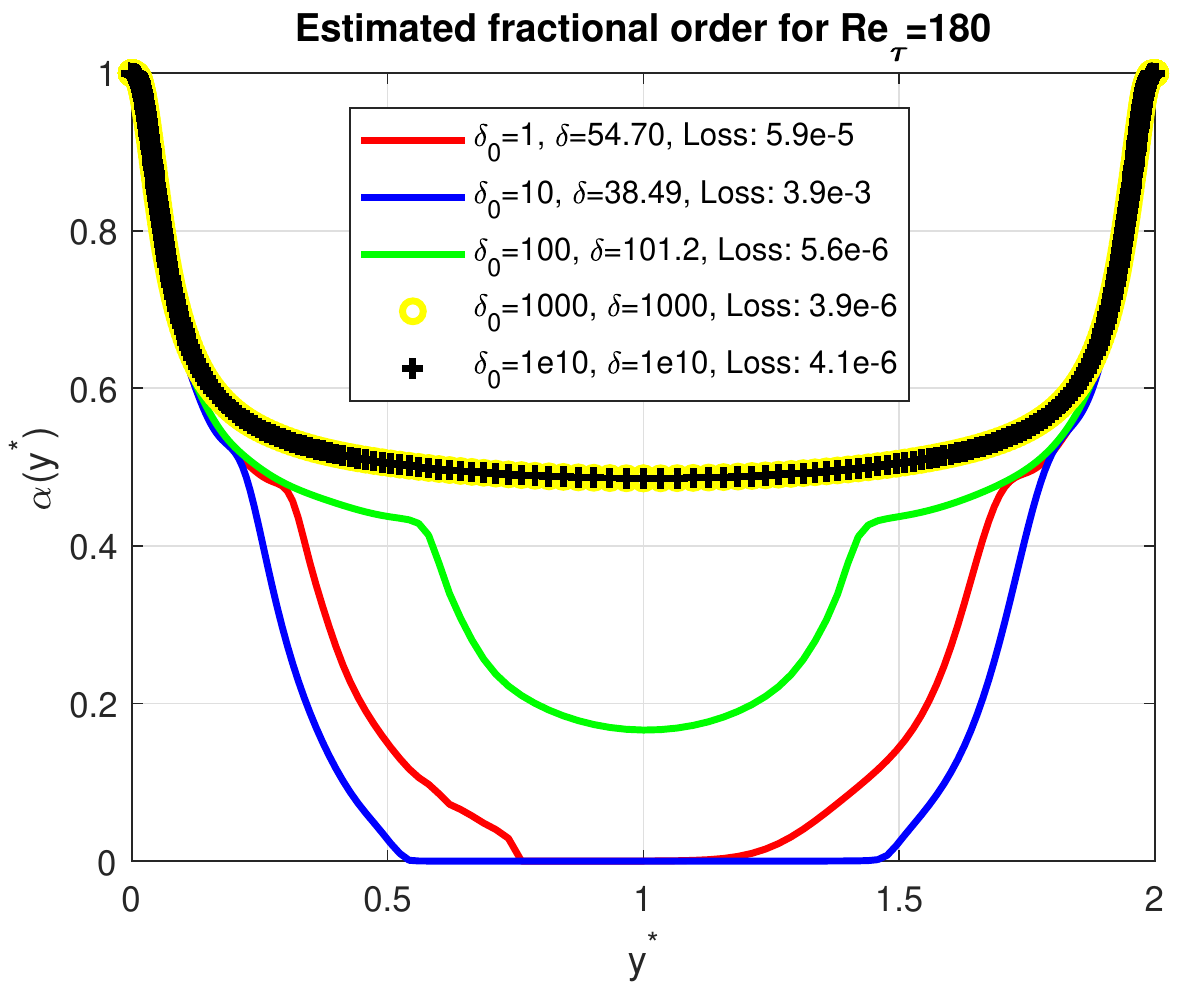}}
\subfloat[]{
\includegraphics[width=.45\textwidth]{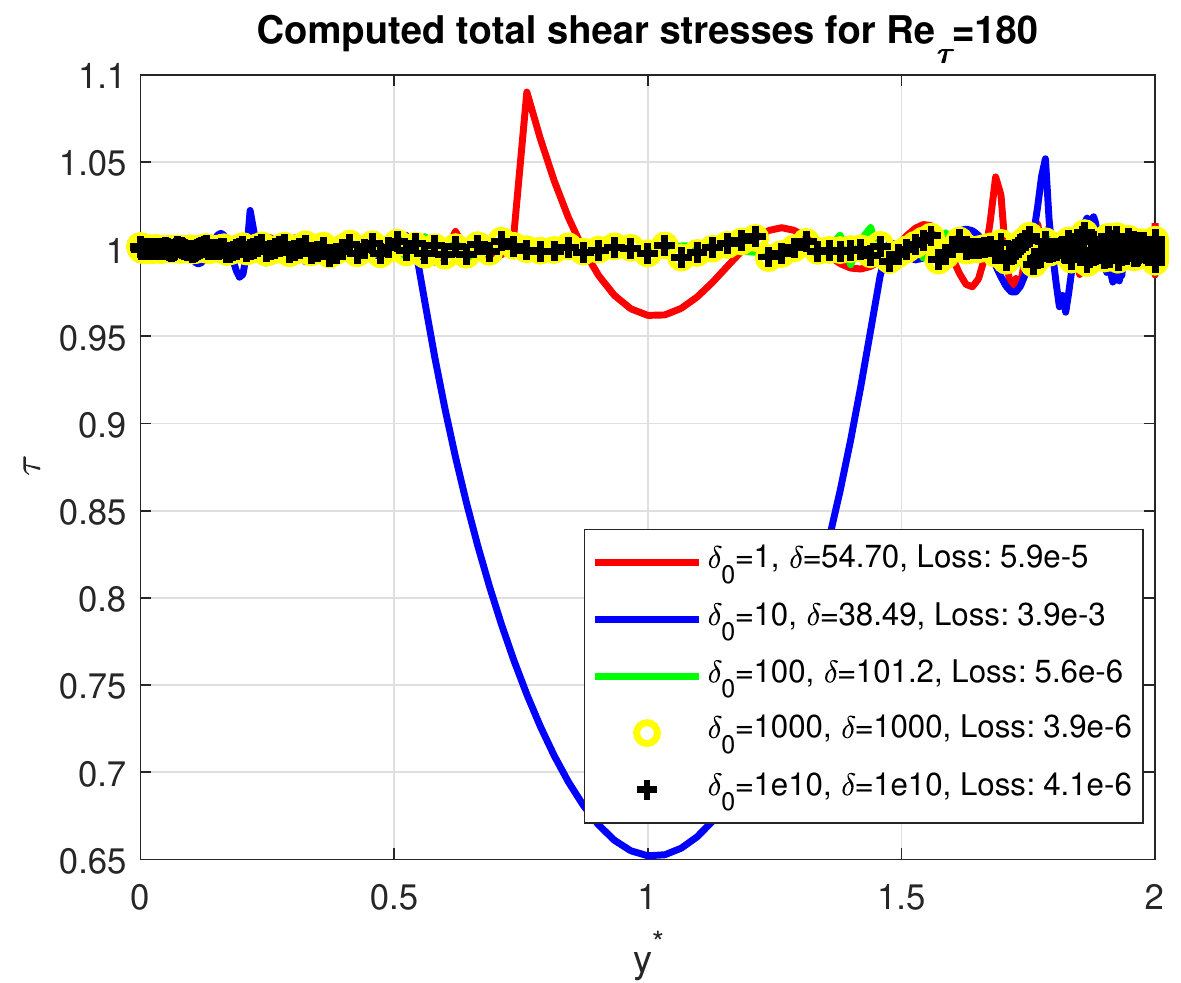}}\vfill
\subfloat[]{
\includegraphics[width=.45\textwidth]{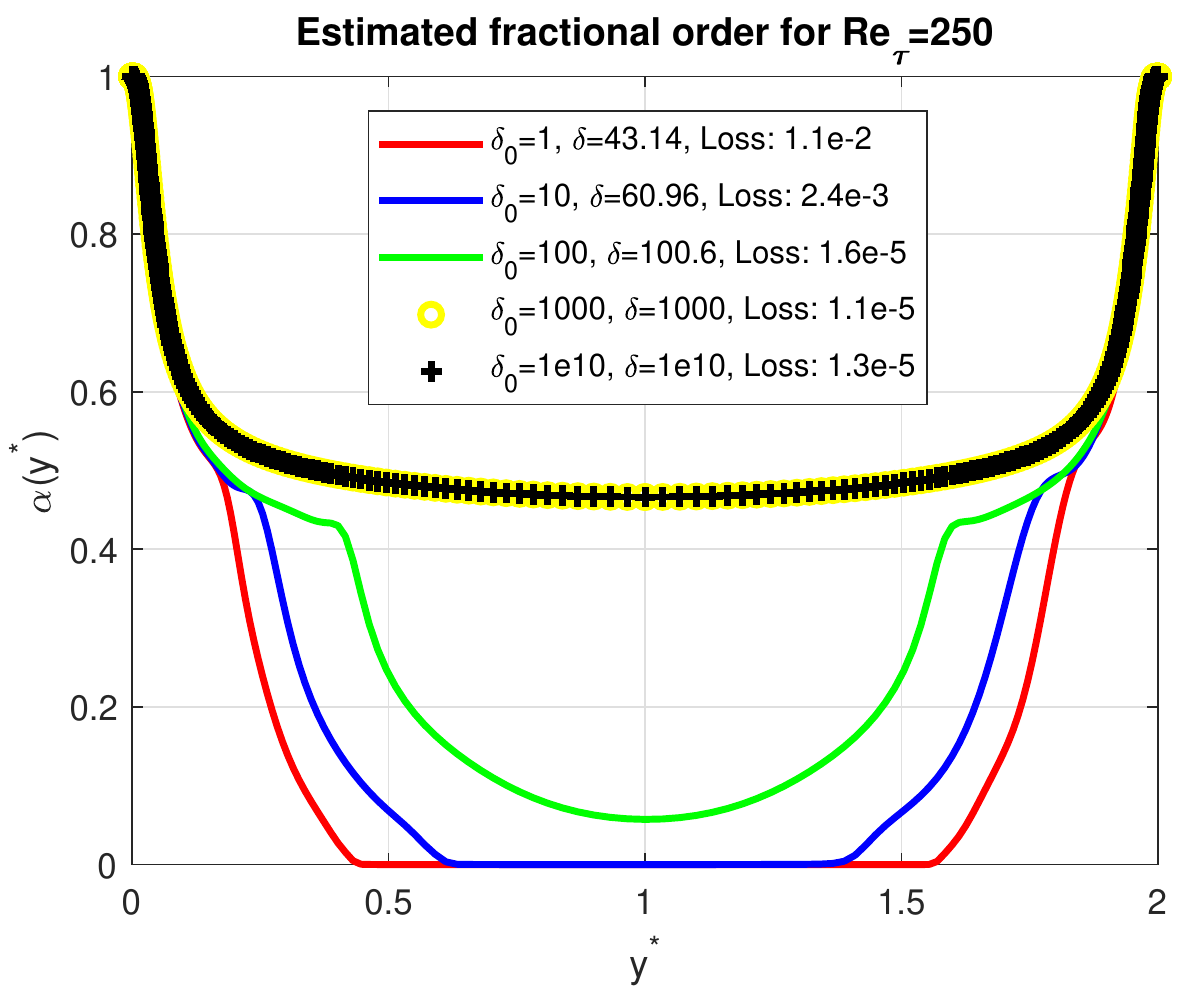}}
\subfloat[]{
\includegraphics[width=.45\textwidth]{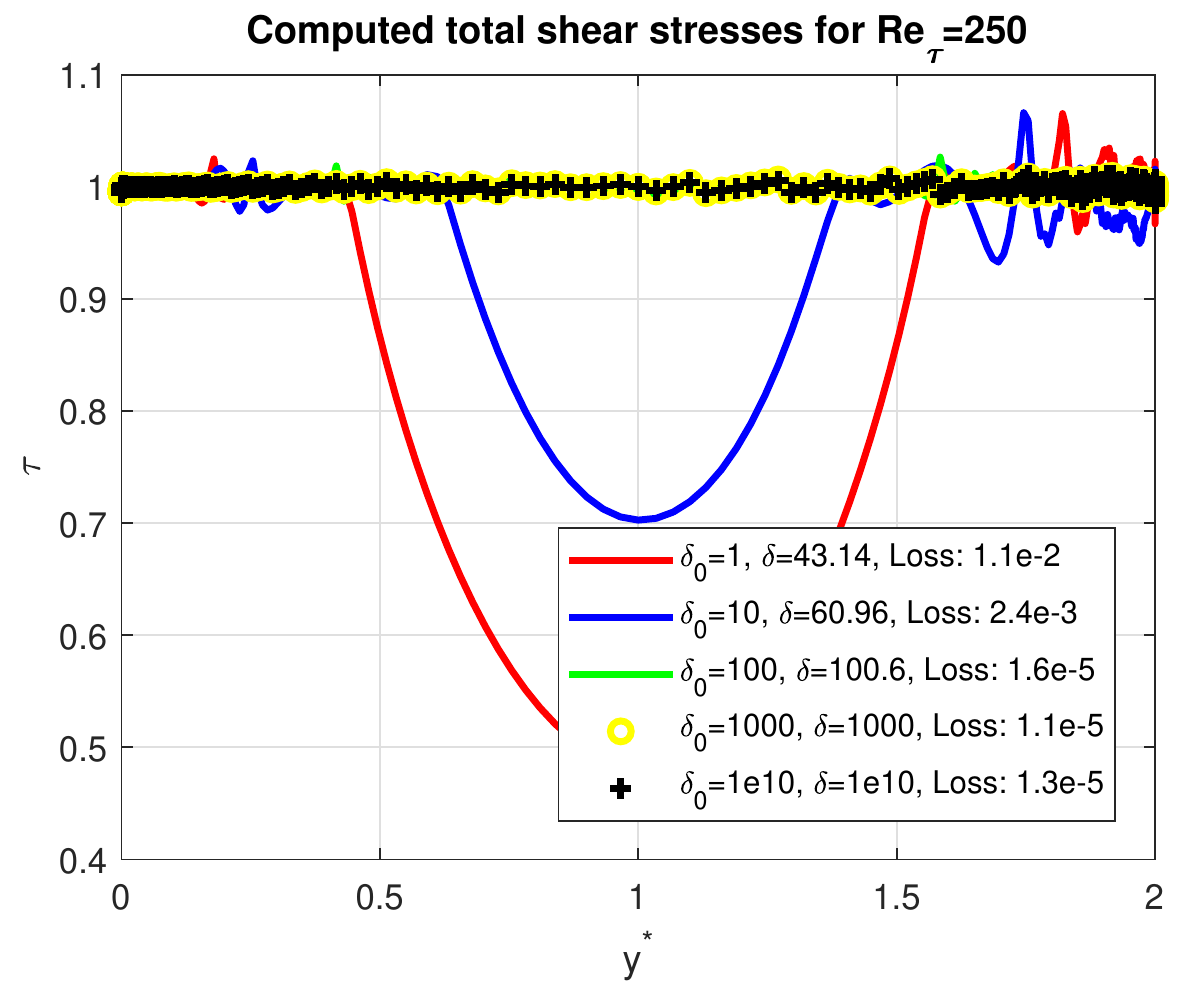}}
\caption{Influence of the initial guess for $\delta_0$ on the learned fractional order and associated total shear stress. The subplots in the first, second, and third rows correspond to $Re_{\tau}$=125, 180, and 250, respectively. The subplots in the first and second columns show the estimated fractional orders and the computed total shear stresses, respectively. The left column indicates that, for a fixed Reynolds number, the final loss value decreases monotonically as the initial guess $\delta_0$ goes from 1 to 1000, which indicates that a large horizon is physically meaningful. The right column indicates that by setting a threshold on the misfit $\Delta\tau^+$ between the computed shear stress and the true one we can exclude the blue and red curves, which correspond to bad local minima. For instance, by disregarding the results for which $\Delta\tau^+>0.01$, the only acceptable cases are the ones corresponding to $\delta_0$=1000 and 1e10.}
\label{initial_guess}
\end{figure}

In Figure \ref{initial_guess}, for Reynolds numbers $Re_{\tau}=125$, 180, and 250, we report the learned fractional orders and horizons, and the corresponding total shear stresses.
First, we observe that for different $\delta_0$ we obtain different fractional order and total shear stress profiles, some of which are non-physical (see profiles for $\delta_0$=1 and 10).
Second, for a fixed Reynolds number, the final loss value decreases monotonically as $\delta_0$ increases. This implies that a large horizon is more physically meaningful than a small one; however, we note that for $\delta \ge 100$ the loss values are of the same order. This suggests that there is a threshold for $\delta$ above which the nPINNs algorithm reaches the same accuracy.
To achieve a universal behavior of both $\delta$ and $\alpha(y^+)$, one could set a threshold, independent of $Re_{\tau}$, on the final loss value or on the misfit $\Delta\tau^+$ between the computed shear stress and the true one. For instance, if we require $\Delta\tau^+<0.01$, then, only the fractional order profiles associated to $\delta_0$=1000 and 1e10 are acceptable, see Figure \ref{initial_guess}. In this case, the resulting two nonlocal operators $\tilde{\mathcal{L}}^{1000,\alpha_{1000}(y^+)}$ and $\tilde{\mathcal{L}}^{1{\rm e}10,\alpha_{1{\rm e}10}(y^+)}$ mimic each other with almost the same fractional order profile but very different horizons. Note that by slightly increasing the aforementioned threshold from 0.01 to 0.013 for $Re_{\tau}$=125 and 180 and from 0.01 to 0.027 for $Re_{\tau}$=250, then, the results associated with $\delta_0=100$ are also acceptable and we have an additional mimicking operator.

\section{Summary}\label{sec:conclusion}

In this paper we introduced a universal nonlocal Laplace operator that bridges classical and fractional Laplacian operators, and developed a computational tool for inferring its parameters from data. Here, we summarize the main properties of this new computational modeling framework.

{\bf 1.} The unified nonlocal operator allows us to describe a broad spectrum of nonlocal elliptic operators by simply tuning the modeling parameters $\alpha$ and $\delta$. It is equivalent to the classical Laplacian at the limit of $\delta\to 0$ for $\alpha\in(0,2)$ and at the limit of $\alpha\to 2$ for $\delta=\infty$, and to the fractional Laplacian for $\delta=\infty$ and $\alpha\in(0,2)$. This property explains the occurrence of the operator mimicking phenomenon observed in Section \ref{sec:numerical-tests}, i.e. different pairs of parameters yield distinct operators that are equally effective in reproducing the training data. This is a consequence of the ill-posed nature of the identification problem, that resembles the non-uniqueness of the diffusivity identification problem in a PDE setting. Operator mimicking can be exploited in the numerical solution of nonlocal elliptic equations; as an example, given $(\delta_1,\alpha_1)$ and $(\delta_2,\alpha_2)$ for which the corresponding operators are equivalent, it is preferable to use the pair with smaller $\delta$ to reduce computational cost.

{\bf 2.} Our algorithm is flexible and requires minimal implementation effort. First, it handles forward and inverse problems in the same manner. For forward problems, the minimization is performed with respect to NN parameters only, whereas for inverse problems it is performed with respect to model parameters as well. Second, the discretization of the nonlocal operators is not tied to any specific method: available software can be used as a black box.

{\bf 3.} For the solution of forward problems, nPINNs are as accurate as other discretization methods such as, e.g. finite difference methods, and exhibit optimal error convergence rates with respect to the number of residual points (or discretization points). Furthermore, they can deal with rough solutions and are applicable in any dimension. 

{\bf 4.} The application to turbulence modeling shows that nPINNs can easily handle high-dimensional parameter spaces. In fact, the algorithm is learning the NN parameters of both the solution and the variable order. Furthermore, our results represent a preliminary step towards resolving the closure problem in turbulence modeling: the optimal parameters, $\alpha(y^+)$ and $\delta$, exhibit a universal behavior with respect to the Reynolds number.  In particular, the introduction of the truncation in the nonlocal operator allows us to identify an optimal interaction length. This fact has important consequences in terms of computational savings since computations with finite $\delta$ are significantly cheaper that with $\delta=\infty$.

{\bf 5.} Our computational framework provides the groundwork for an important open problem: the identification of the kernel function $\gamma$, i.e., the identification of the nonlocal operator itself. This is intrinsically a much more complex problem as it involves learning a functional form rather than model parameters and it is a high-dimensional identification problem. The identification of the variable decay rate in Section \ref{sec:turbulence} is a first step towards this challenging problem and suggests that one possible venue is to represent the kernel function with a NN. This is the subject of our future work, and the theoretical results of \cite{chen1995universal} on the universal approximation of functionals is encouraging. 

{\bf 6.} Tempered fractional operators do not belong to the set of operators spanned by $\mcL^{\delta,\alpha}$, however, truncation seems to have a similar effect as tempering, i.e. the second moment of truncated operators, such as ours, is finite, as for the tempered case. Also, the extension to generalized tempered operators only requires multiplication of the kernel function in \eqref{kernel} by the factor $\exp\{-\lambda\|\xb-\yb\|\}$, for some $\lambda>0$. Adapting nPINNs to this class of operators is straightforward and is part of our ongoing work.

\section{Acknowledgments}
This work was supported by the U.S. Department of Energy, Office of Science, Office of Advanced Scientific Computing Research under the Physics-Informed Learning Machines for Multiscale and Multiphysics Problems (PhILMs) project. GP and GK are also supported by the MURI/ARO at Brown University (W911NF-15-1-0562) and DARPA-AIRA (HR00111990025).
MD and MP are also supported by Sandia National Laboratories (SNL). SNL is a multimission laboratory managed and operated by National Technology and Engineering Solutions of Sandia, LLC., a wholly owned subsidiary of Honeywell International, Inc., for the U.S. Department of Energy's National Nuclear Security Administration under contract DE-NA-0003525. 

\smallskip 
The authors would like to thank Pavan Pranjivan Mehta (Brown University) for useful discussion on connections and possible extensions to tempered fractional derivatives.

\smallskip
This paper describes objective technical results and analysis. Any subjective views or opinions that might be expressed in the paper do not necessarily represent the views of the U.S. Department of Energy or the United States Government. Report number SAND2020-3980.

\begin{appendices}
\section{Proof of Lemma \ref{lem:limiting-behavior}}\label{sec:AppA}
Proofs for the two- and three-dimensional cases are based on transformations into polar and spherical coordinates respectively. We report steps of the proof for the two-dimensional case only, as the three-dimensional one can be treated in a similar manner. Without loss of generality, as done for the one-dimensional case, we assume that $u\equiv 0 \mbox{ for } \mathbf{x} \in \mathbb{R}^2\setminus\Omega$.

\smallskip\noindent 
I. The operator $\mathcal{L}^{\delta,\alpha}$ satisfies \eqref{lim1}.\\
Let $\yb=(y_1,y_2)$ and $\xb=(x_1,x_2)$; the transformation $y_1=x_1+r\cos\theta$ and $y_2=x_2+r\sin\theta$ yields
\begin{equation}
    \int_{B_{\delta}(\xb)}\frac{u(\yb)-u(\xb)}{|\yb-\xb|^{2+\alpha}}d\yb=\int_{0}^{\delta}\int_{0}^{2\pi}\frac{u(x_1+r\cos\theta,x_2+r\sin\theta)-u(x_1,x_2)}{r^{2+\alpha}}rd\theta dr.
\end{equation}
By Taylor series expansion we have
\begin{equation*}
   \begin{split}
    \int_{B_{\delta}(\xb)}\frac{u(\yb)-u(\xb)}{|\yb-\xb|^{2+\alpha}}d\yb&=\int_{0}^{\delta}\int_{0}^{2\pi}\frac{r\cos\theta\frac{\partial u}{\partial x_1}+r\sin\theta\frac{\partial u}{\partial x_2}}{r^{1+\alpha}}d\theta dr \\
    &+ \int_0^{\delta}\int_0^{2\pi}\frac{r^2\cos^2\theta\frac{\partial^2u}{\partial x_1^2}+2r^2\sin\theta\cos\theta\frac{\partial^2 u}{\partial x\partial y}+r^2\sin^2\theta\frac{\partial^2 u}{\partial x_2^2}}{2r^{1+\alpha}}d\theta dr \\
    & + \int_0^{\delta}\int_0^{2\pi}\frac{o(r^2)}{r^{1+\alpha}}d\theta dr.
   \end{split}    
\end{equation*}
The first integral on the right-hand side vanishes due to the fact that $\int_0^{2\pi}\cos\theta d\theta =\int_0^{2\pi}\sin\theta d\theta = 0$; by explicitly integrating in $\theta$ and $r$, the second integral equals $\frac{\pi\delta^{2-\alpha}\Delta u}{2(2-\alpha)}=\frac{\Delta u}{C'_{\delta,\alpha}}$; and the third integral is of order $\delta^\beta$ with $\beta>2-\alpha$. Since for $\delta\to 0$ the third term can be neglected, we have
\begin{equation*}
\lim\limits_{\delta\rightarrow0}C_{\delta,\alpha}\int_{B_{\delta}(\xb)}\frac{u(\yb)-u(\xb)}{|\yb-\xb|^{2+\alpha}}d\yb=\Delta u(\xb)+\lim\limits_{\delta\rightarrow0}\frac{C''_{\delta,\alpha}}{C'_{\delta,\alpha}}\Delta u(\xb)=\Delta u(\xb).
\end{equation*}

\smallskip\noindent
II. The operator $\mathcal{L}^{\delta,\alpha}$ satisfies \eqref{lim2}.\\
We consider the difference between the fractional Laplacian and the unified operator:
\begin{equation}
 \begin{split}
|-\mathcal{L}^{\delta,\alpha}u(\xb)-(-\Delta)^{\alpha/2}u(\xb)| & =\left|-C'_{\delta,\alpha}\int_{0}^{\delta}\int_0^{2\pi}\frac{u(\xb+r\thetab)-u(\xb)}{r^{2+\alpha}}rd\theta dr\right. \\
& \left. + C''_{\delta,\alpha}\int_{\delta}^{+\infty}\int_0^{2\pi}\frac{u(\xb+r\thetab)-u(\xb)}{r^{2+\alpha}}rd\theta dr \right|\\
& = \left|-C'_{\delta,\alpha}\int_{0}^{\delta}\int_0^{2\pi}\frac{u(\xb+r\thetab)-u(\xb)}{r^{2+\alpha}}rd\theta dr\right. \\
& \left.
-\frac{2\pi C''_{\delta,\alpha}}{\alpha}u(\xb)\delta^{-\alpha}\right|,
\end{split}    
\end{equation}
where $\thetab=(\cos\theta,\sin\theta)$ and where we used the fact that $u(\xb+r\thetab)=0$ for large $\delta$, see the assumption at the beginning of the proof. The thesis follows from the fact that, as $\delta\to\infty$, the first term on the right-hand side converges to zero with the asymptotic rate $\delta^{\alpha-2}$.

\section{Evaluating singular integrals}\label{sec:singular-integrals}
\label{sec:AppB}

\subsection{The one-dimensional case}
We describe the evaluation of the nonlocal operator $-\mathcal{L}^{\delta,\alpha}u(x_k;\boldsymbol{\mu})$ at the residual point $x_k\in\omg$ for given $\mub$ in the one-dimensional case. Note that we directly prescribe the nonlocal boundary condition, i.e. $u=g$ whenever we evaluate the $u$ in $\omgd$; instead, for $x_k\in\omg$, we set $u=\unn$. We introduce the constants $\rho\ll 1$ and $D=|\omg|$ and describe how to evaluate the integral according to the magnitude of $\delta$ with respect to $\rho$ and $D$. In order to avoid integration in the neighborhood of a discontinuity of $u$, we do not place any residual point in the following two regions: 

$$\Omega'_{\rho}=\bigcup\limits_{x\in\partial\Omega}B_{\rho}(x) 
\quad {\rm and} \quad
\Omega''_{\rho}=\bigcup\limits_{x\in \mathcal{S}}B_{\rho}(x),$$
where $\mathcal{S}$ denotes the set of points with a jump discontinuity. Thus, $\omg'_\rho$ consists of two layers of thickness $\rho$ surrounding the boundary of $\omg$ (one inside $\omg$ and one outside) and $\omg''_\rho$ is the set of all balls of radius $\rho$ surrounding discontinuity points.

\smallskip\noindent {\bf A.} $\delta\le\rho$\\
Using the Taylor expansion in (\ref{ss}), we have
\begin{equation}\label{append}
\int_{-\delta}^{\delta}\frac{u(x+z)-u(x)}{|z|^{1+\alpha}}dz =
\int_{-\delta}^{\delta}\frac{\unn(x+z)-\unn(x)}{|z|^{1+\alpha}}dz \approx \unn''(x)\frac{\delta^{2-\alpha}}{2-\alpha},\quad x\in\Omega\setminus(\Omega'_{\rho}\cup\Omega''_{\rho}).
\end{equation}
Note that in this case we only evaluate $u$ inside of $\omg$ so that $u=\unn$.

\smallskip\noindent {\bf B.} $\rho<\delta\le D$\\
We divide the integral into sub-integrals, isolating the singularity:
\begin{equation}\label{delta2}
\begin{split}
\int_{-\delta}^{\delta}\frac{u(x+z)-u(x)}{|z|^{1+\alpha}}dz & = \int_{-\rho}^{\rho}\frac{\unn(x+z)-\unn(x)}{|z|^{1+\alpha}}dz \\
& + \int_{-\delta}^{-\rho}\frac{\tilde{u}(x+z)-\unn(x)}{|z|^{1+\alpha}}dz + \int_{\rho}^{\delta}\frac{\tilde{u}(x+z)-\unn(x)}{|z|^{1+\alpha}}dz.
\end{split} 
\end{equation}
In the first sub-integral we used the fact that, for $z\in(-\rho,\rho)$, $x+z\in\omg$ $\tilde{u}(x+z)$ so that $u(z+x)=\unn(x+z)$. Thus, we approximate that sub-integral as in \eqref{append}, i.e.
\begin{equation*}
   \int_{-\rho}^{\rho}\frac{\unn(x+z)-\unn(x)}{|z|^{1+\alpha}}dz \approx u''_{N\!N}(x)\frac{\rho^{2-\alpha}}{2-\alpha}.
\end{equation*}
In the second and third sub-integral, we have that $|z|>\rho$; here, we let $\tilde{u}(x+z)$ be $\unn(x+z)$ for $x+z \in \Omega$ and $g(x+z)$ for $x+z \in \mathbb{R}\setminus\Omega$. We evaluate these integrals via composite Gauss quadrature, i.e.
\begin{equation}\label{GL_quad_1d}
\begin{split}
& \int_{-\delta}^{-\rho}\frac{\tilde{u}(x+z)-\unn(x)}{|z|^{1+\alpha}}dz + \int_{\rho}^{\delta}\frac{\tilde{u}(x+z)-\unn(x)}{|z|^{1+\alpha}}dz \\
&\approx \sum_{j=1}^m\sum_{k=1}^M w_{kj}\frac{\tilde{u}(x+z_{kj})-2\unn(x)+
\tilde{u}(x-z_{kj})}{z_{kj}^{1+\alpha}},
 \end{split}    
\end{equation}
where $\{z_{kj}\}_{k=1}^M$ and $\{w_{kj}\}_{k=1}^M$ are $M$ Gauss-Legendre quadrature points and weights in the $j$-th interval $(a_{j-1},a_j)$ ($j=1,2,\cdots,m$), respectively. The integration domain $(\rho,\delta)$ is divided into $m$ sub-intervals: $\rho=a_0<a_1<\cdots<a_{m-1}<a_m=\delta$. 

Note that there are two reasons to employ composite quadrature rules. First, for very small $\rho$, they mitigate the evaluation of nearly singular integrals. Second, when integrating possibly non-smooth or discontinuous functions, they provide a better accuracy by increasing the number of sub-intervals, $m$.

\smallskip\noindent {\bf C.} $\delta>D$\\ 
As done in {\bf B.} we split the integral as follows:
\begin{equation}\label{delta3}
\begin{split}
\int_{-\delta}^{\delta}\frac{\tilde{u}(x+z)-\tilde{u}(x)}{|z|^{1+\alpha}}dz & = \int_{-\rho}^{\rho}\frac{\unn(x+z)-\unn(x)}{|z|^{1+\alpha}}dz \\
& + \int_{-D}^{-\rho}\frac{\tilde{u}(x+z)-\unn(x)}{|z|^{1+\alpha}}dz + \int_{\rho}^{D}\frac{\tilde{u}(x+z)-\unn(x)}{|z|^{1+\alpha}}dz \\
& + \int_{-\delta}^{-D}\frac{g(x+z)-\unn(x)}{|z|^{1+\alpha}}dz + \int_D^{\delta}\frac{g(x+z)-\unn(x)}{|z|^{1+\alpha}}dz.
\end{split} 
\end{equation}
Here, the first sub-integral is again approximated by Taylor expansion as in \eqref{append}, whereas the second and the third sub-integrals are evaluated by using composite Gauss quadrature rules. By using the volume constraint, the fourth and fifth sub-integrals can be rewritten as:
\begin{equation}
\begin{split}
&\int_{-\delta}^{-D}\frac{g(x+z)-\unn(x)}{|z|^{1+\alpha}}dz + \int_D^{\delta}\frac{g(x+z)-\unn(x)}{|z|^{1+\alpha}}dz \\
& = \int_D^{\delta}\frac{g(x+z)+g(x-z)}{z^{1+\alpha}}dz +\frac{2\unn(x)(\delta^{-\alpha}-D^{-\alpha})}{\alpha},
\end{split}    
\end{equation}
where the first integral above can be approximated by using standard Gauss-Legendre quadrature.

We unify cases {\bf A}, {\bf B} and {\bf C} above in one formula. Let $a_0=\min\{\delta,\rho\}$, $a_m=\min\{\delta,D\}$, and $x\in\Omega\setminus(\Omega'_{\rho}\cup\Omega''_{\rho})$, we have
\begin{equation}\label{quad_1d}
\begin{split}
\int_{-\delta}^{\delta}\frac{\tilde{u}(x+z)-\tilde{u}(x)}{|z|^{1+\alpha}}dz & \approx \unn''(x)\frac{a_0^{2-\alpha}}{2-\alpha} \\
& +\sum_{j=1}^m\sum_{k=1}^M w_{kj}\frac{\tilde{u}(x+z_{kj})-2\unn(x)+\tilde{u}(x-z_{kj})}{z_{kj}^{1+\alpha}} \\
& + \left(\int_D^{\delta}\frac{g(x+z)+g(x-z)}{z^{1+\alpha}}dz +\frac{2\unn(x)(\delta^{-\alpha}-D^{-\alpha})}{\alpha}\right)H(\delta-D),
\end{split} 
\end{equation}
where $H(x)=1$ for $x>0$ and $H(x)=0$ for $x\le 0$.

\subsection{The multivariate case} \label{sec:AppC}
We extend the description of the previous section to the two- and three-dimensional setting.

\paragraph{Two-dimensional case} For $\thetab=(\cos\theta,\sin\theta)$, $a_0=\min\{\delta,\rho\}$, $a_{m_1}=\min\{\delta,D\}$, and $\xb\in\Omega\setminus(\Omega'_{\rho}\cup\Omega''_{\rho})$, we have
\begin{multline}\label{quad_2d}
   \int_{0}^{\delta}\int_0^{2\pi}\frac{\tilde{u}(\xb+r\thetab)-\tilde{u}(\xb)}{r^{2+\alpha}}rd\theta dr  \approx \Delta \unn(\xb)\frac{\pi a_0^{2-\alpha}}{2(2-\alpha)} \\
    +\sum_{j_1=1}^{m_1}\sum_{j_2=1}^{m_2}\sum_{k_1=1}^{M_1}\sum_{k_2=1}^{M_2} w^{j_1j_2}_{k_1}w^{j_1j_2}_{k_2}\frac{\tilde{u}(\xb+r^{j_1j_2}_{k_1}\thetab^{j_1j_2}_{k_2})-\unn(\xb)}{(r^{j_1j_2}_{k_1})^{1+\alpha}} \\
    + \left(\int_D^{\delta}\int_0^{2\pi}\frac{g(\xb+r\thetab)}{r^{1+\alpha}}d\theta dr +\frac{2\pi \unn(\xb)(\delta^{-\alpha}-D^{-\alpha})}{\alpha}\right)H(\delta-D),
 \end{multline} 
where $\thetab^{j_1j_2}_{k_2}=(\cos\theta^{j_1j_2}_{k_2},\sin\theta^{j_1j_2}_{k_2})$. The annulus in polar coordinate $(r,\theta)\in [a_0,a_{m_1}]\times [0,2\pi]$ is divided into $m_1\times m_2$ segments  $S_{j_1j_2}\overset{\Delta}{=}[a_{j_1-1},a_{j_1}]\times[b_{j_2-1}, b_{j_2}]$ for $j_1=0,1,\cdots,m_1$ and $j_2=0,1,\cdots,m_2$ where $a_0<a_1<\cdots<a_{m_1-1}<a_{m_1}$ and $0=b_0<b_1<\cdots<b_{m_2-1}<b_{m_2}=2\pi$. In each segment $S_{j_1j_2}$ we employ Gauss-Legendre quadrature with integration points $(r^{j_1j_2}_{k_1},\theta^{j_1j_2}_{k_2})$ and integration weights $w^{j_1j_2}_{k_1}w^{j_1j_2}_{k_2}$.

\paragraph{Three-dimensional case} 
For $\thetab=(\sin\phi\cos\theta,\sin\phi\sin\theta,\cos\phi)$, $a_0=\min\{\delta,\rho\}$,\\ $a_{m_1}=\min\{\delta,D\}$, and $\xb\in\Omega\setminus(\Omega'_{\rho}\cup\Omega''_{\rho})$, we have
\begin{multline}\label{quad_3d}
   \int_{0}^{\delta}\int_0^{2\pi}\int_0^\pi\frac{\tilde{u}(\xb+r\thetab)-\tilde{u}(\xb)}{r^{3+\alpha}}r^2\sin\phi d\phi d\theta dr  \approx \Delta \unn(\xb)\frac{2\pi a_0^{2-\alpha}}{3(2-\alpha)} \\
    +\sum_{j_1=1}^{m_1}\sum_{j_2=1}^{m_2}\sum_{j_3=1}^{m_3}\sum_{k_1=1}^{M_1}\sum_{k_2=1}^{M_2}\sum_{k_3=1}^{M_3} w^{j_1j_2j_3}_{k_1}w^{j_1j_2j_3}_{k_2}w^{j_1j_2j_3}_{k_3}\sin\phi^{j_1j_2j_3}_{k_3}\frac{\tilde{u}(\xb+r^{j_1j_2j_3}_{k_1}\thetab^{j_1j_2j_3}_{k_2k_3})-\unn(\xb)}{(r^{j_1j_2j_3}_{k_1})^{1+\alpha}} \\
    + \left(\int_D^{\delta}\int_0^{2\pi}\int_0^{\pi}\frac{g(\xb+r\thetab)\sin\phi}{r^{1+\alpha}}d\phi d\theta dr +\frac{4\pi \unn(\xb)(\delta^{-\alpha}-D^{-\alpha})}{\alpha}\right)H(\delta-D),
  \end{multline}
where $\thetab^{j_1j_2j_3}_{k_2k_3}=(\cos\theta^{j_1j_2j_3}_{k_2}\sin\phi^{j_1j_2j_3}_{k_3},\sin\theta^{j_1j_2j_3}_{k_2}\sin\phi^{j_1j_2j_3}_{k_3},\cos\phi^{j_1j_2j_3}_{k_3})$. The spherical shell in spherical coordinate $(r,\theta,\phi)\in [a_0,a_{m_1}]\times [0,2\pi] \times [0,\pi]$ is divided into $m_1\times m_2 \times m_3$ segments  $S_{j_1j_2j_3}\overset{\Delta}{=}[a_{j_1-1},a_{j_1}],[b_{j_2-1}\times b_{j_2}]\times[c_{j_3-1},c_{j_3}]$ for $j_1=0,1,\cdots,m_1$, $j_2=0,1,\cdots,m_2$, and $j_3=0,1,\cdots,m_3$ where $a_0<a_1<\cdots<a_{m_1-1}<a_{m_1}$, $0=b_0<b_1<\cdots<b_{m_2-1}<b_{m_2}=2\pi$, and $0=c_0<c_1<\cdots<c_{m_3-1}<c_{m_3}=\pi$. In each segment $S_{j_1j_2j_3}$ we employ Gauss-Legendre quadrature with integration points $(r^{j_1j_2j_3}_{k_1},\theta^{j_1j_2j_3}_{k_2},\phi^{j_1j_2j_3}_{k_3})$ and integration weights $w^{j_1j_2j_3}_{k_1}w^{j_1j_2j_3}_{k_2}w^{j_1j_2j_3}_{k_3}$.

\section{Quadrature accuracy}\label{sec:AppD}
We test the quadrature accuracy of the schemes (\ref{quad_1d}), (\ref{quad_2d}), and (\ref{quad_3d}). We consider the one-dimensional function
\begin{equation}\label{solu1}
   u_1(x)=\left\{
   \begin{array}{cc}
       x(1-x^2)^{1+\alpha/2},& \quad |x|\le 1 \\
       0, & \quad |x|>1,
   \end{array}  \right.
\end{equation}   
and the multi-dimensional functions with $\mathbf{x}\in \mathbb{R}^2$ or $\mathbb{R}^3$
\begin{equation}\label{solu2}
   u_2(\xb) (\text{or } u_3(\xb)) =\left\{
   \begin{array}{cc}
       (1-\|\xb\|_2^2)^{1+\alpha/2},& \quad \|\xb\|_2\le 1 \\
       0, & \quad \|\xb\|_2>1.
   \end{array}  \right.
\end{equation}  
The fractional Laplacian of the above functions can be computed analytically:
\begin{equation} 
   \begin{split}
    f_1(x) & = (-\Delta)^{\alpha/2}u_1(x)=\frac{\Gamma(\alpha+3)}{6}(3-(3+\alpha)x^2)x, \quad x\in\mathbb{R}, \\
    f_2(\xb) & =  (-\Delta)^{\alpha/2}u_2(\xb)= 2^{\alpha}\Gamma(\frac{\alpha}{2}+2)\Gamma(\frac{2+\alpha}{2})(1-(1+\frac{\alpha}{2})||\xb||_2^2), \quad \mathbf{x}\in \mathbb{R}^2, \\
    f_3(\xb) & = (-\Delta)^{\alpha/2}u_3(\xb)= 2^{\alpha}\Gamma(\frac{\alpha}{2}+2)\Gamma(\frac{3+\alpha}{2})\Gamma(1.5)^{-1}(1-(1+\frac{\alpha}{3})||\boldsymbol{x}||_2^2), \quad \mathbf{x}\in \mathbb{R}^3.
  \end{split}    
\end{equation}

In all these tests the decay exponent and the interaction radius are set to $\alpha=1.5$, and $\delta=10^{100}$, respectively. To show the quadrature accuracy, we define the relative error
\begin{equation}
    \epsilon_q^k = \frac{\sqrt{\sum_{j=1}^{100}((-\mathcal{L}^{\delta,\alpha}u_k(\xb^t_j)-(-\Delta)^{\alpha/2}u_k(\xb^t_j))^2}}{\sqrt{\sum_{j=1}^{100}((-\Delta)^{\alpha/2}u_k(\xb^t_j))^2}}\quad k=1,2,3.
\end{equation}
The test points for the one-dimensional case are taken in the following way: let $\{s_j\}_{j=1}^{100}\subset(0,1)$ be the first 100 points in the Sobol sequence; the test points are taken as the transformation $x^t_j=-1+\rho+(2-2\rho)*s_j$ since we exclude the neighborhood of the domain boundary $\Omega'_{\rho}$. 

For two-dimensional case, let $\{\mathbf{s}\}_j=\{s^1_j,s^2_j\}\subset (0,1)^2$ be the first 200 points in the Sobol the sequence; the test points are defined as $\xb^t_j=(s^1_j\cos s^2_j,s^1_j\sin s^2_j)$. 

\begin{figure}[H] 
\centering
\includegraphics[width=.95\textwidth]{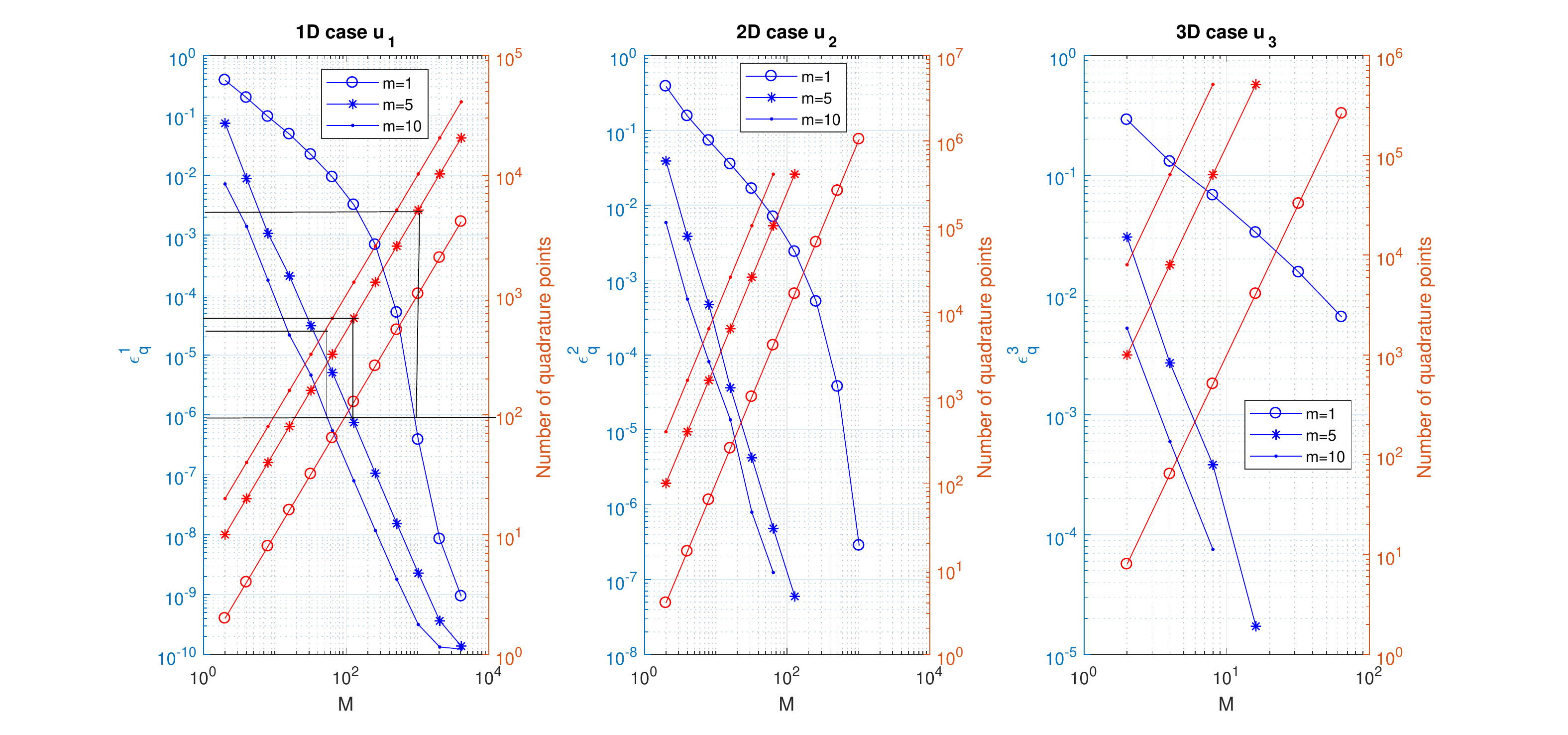}
\caption{Effects of number of sub-integrals ($m$) and number of Gauss-Legendre points in each sub-integral ($M$) when composite Gauss quadrature is performed. From the left $y-$axis (relative error), we see that quadrature accuracy drops with increasing $m$ and $M$. From the right $y-$axis (total number of quadrature points), we see that to achieve higher accuracy more quadrature points are needed. Note that the total number of quadrature points is $mM$, $m^2M^2$, and $m^3M^3$ for 1D, 2D, and 3D cases, respectively.}
\label{quad_acc}
\end{figure}

In all dimensions, the diameters of $\Omega$ is $D=2$; we set $\rho=10^{-5}$. The effects of the number of sub-integrals $m$ and the number of the Gauss-Legend quadrature points in each sub-integral, $M$, are shown in Figure \ref{quad_acc}.


\end{appendices}

\bibliographystyle{plain}

\end{document}